\documentclass[11pt,reqno]{article}

\usepackage{amsfonts}
\usepackage{amsmath}
\usepackage{amsthm}
\usepackage{amssymb}
\usepackage{mathrsfs}
\usepackage{amstext}
\usepackage{graphicx}
\usepackage{color}
\usepackage{url}
\font\msbm=msbm10

\numberwithin{equation}{section}

\theoremstyle{plain}
\newtheorem{Theorem}{Theorem}[section]
\newtheorem{lemma}[Theorem]{Lemma}
\newtheorem{corollary}[Theorem]{Corollary}
\newtheorem{example}[Theorem]{Example}
\newtheorem{proposition}[Theorem]{Proposition}
\newtheorem{definition}[Theorem]{Definition}

\newtheorem{remark}[Theorem]{Remark}
\def\mathbb#1{\hbox{\msbm{#1}}}

\newcommand{\re}{\mathbb{R}}\newcommand{\N}{\mathbb{N}}
\newcommand{\zz}{\mathbb{Z}}\newcommand{\C}{\mathbb{C}}

\newcommand{\com}{\mathbb{C}}
\newcommand{\Z}{{\zz}^d}
\newcommand{\n}{{\N}_0}
\newcommand{\R}{{\re}^d}

\newlength{\fixboxwidth}
\newcommand{\bproof}{\noindent {\bf Proof.}\, }
\newcommand{\eproof}{\hspace*{\fill} \rule{3mm}{3mm}}

\newcommand{\Co}{\mathsf{Co}}
\newcommand{\CoY}{\Co Y}

\newcommand{\x}{\mathbf{x}}
\newcommand{\y}{\mathbf{y}}
\newcommand{\z}{\mathbf{z}}

\newcommand{\supp}{{\rm supp \, }}

\newcommand{\cb}{{\mathcal B}}

\newcommand{\cp}{{\mathcal P}}

\newcommand{\cx}{{\mathcal X}}

\newcommand{\ch}{{\mathcal H}}

\newcommand{\cn}{{\mathcal N}}
\newcommand{\cf}{{\mathcal F}}

\newcommand{\cg}{{\mathcal G}}

\newcommand{\cu}{{\mathcal U}}
\newcommand{\osc}{{\operatorname{osc}}}

\newcommand{\dopp}[1]{\mathbb{#1}}

\newcommand{\esssup}[1]{\underset{#1}{\operatorname{ess\,sup}}}

\newcommand{\CoP}{\Co P^w_{B,q,a}}
\newcommand{\CoL}{\Co L^w_{B,q,a}}

\begin{document}
\title{Generalized coorbit space theory and inhomogeneous function spaces of Besov-Lizorkin-Triebel type}

\author{Holger Rauhut, Tino Ullrich\footnote{\tt rauhut@hcm.uni-bonn.de, tino.ullrich@hcm.uni-bonn.de}\\\newline\\
Hausdorff Center for Mathematics \& Institute for Numerical Simulation\\ 
Endenicher Allee 60, 53115 Bonn, Germany}

\maketitle

\begin{abstract} Coorbit space theory is an abstract approach to function spaces and their atomic decompositions. The original theory developed by Feichtinger and Gr{\"o}chenig in the late 1980ies heavily uses integrable representations of locally compact groups. 
Their theory covers, in particular, homogeneous Besov-Lizorkin-Triebel spaces, modulation spaces, Bergman spaces and the recent shearlet spaces.
However, inhomogeneous Besov-Lizorkin-Triebel spaces cannot be covered by their group theoretical approach.
Later it was recognized by Fornasier and the first named author \cite{fora05}
that one may replace coherent states related to the group representation by more general abstract continuous frames.
In the first part of the present paper we significantly extend this abstract generalized coorbit space theory
% developed by Fornasier and Rauhut \cite{fora05} 
to treat a wider variety of coorbit spaces. A unified 
approach towards atomic decompositions and Banach frames with new results for general coorbit spaces is presented. 
In the second part we apply the abstract setting to a specific framework and study coorbits of what we call Peetre 
spaces. They allow to recover inhomogeneous Besov-Lizorkin-Triebel spaces of various types of interest as coorbits. 
We obtain several old and new wavelet characterizations based on explicit smoothness, decay, and vanishing moment assumptions of the respective wavelet. %We emphasize that inhomogeneous spaces can not be treated with the classical theory of Feichtinger and Gr\"ochenig. 
As main examples we obtain results for weighted spaces (Muckenhoupt, doubling), general $2$-microlocal spaces, Besov-Lizorkin-Triebel-Morrey spaces, spaces of dominating mixed smoothness and even mixtures of the mentioned ones. 
Due to the generality of our approach, there are many more examples of interest where the abstract coorbit space theory is applicable. 
\end{abstract}

%\begin{tabbing}

\medskip


\newpage
\section{Introduction}
Coorbit space theory was originally developed by Feichtinger and Gr{\"o}chenig \cite{FeGr86, Gr88, Gr91}
with the aim to provide a unified approach for describing function spaces and their atomic decompositions, that is,
characterizations via (discrete) sequence spaces. Their theory uses locally compact groups together with an integrable group representation
as a key ingredient. The idea is to measure smoothness via properties of an abstract wavelet transform (the voice transform)
associated to the integrable group representation. More precisely, one asks whether the transform is contained in certain function 
spaces (usually $L_p$-spaces) on the index set of the transform, which is the underlying group.  
As main examples classical homogeneous Besov-Lizorkin-Triebel spaces \cite{Tr83,Tr88,Tr92} can be identified
as coorbit spaces \cite{T10}, and the abstract theory provides characterizations via wavelet frames. Also modulation spaces  and characterizations
via Gabor frames \cite{gr01,fe83-4}, Bergman spaces \cite{FeGr86}, and the more recent shearlet spaces \cite{dakustte09} can be treated via
classical coorbit space theory. In \cite{ra05-3} this theory was extended 
in order to treat also quasi-Banach function spaces. 

Later it was recognized that certain transforms and associated function spaces of interest 
do not fall into the classical group theoretical setting, and the theory was further generalized from groups to the setting
of homogeneous spaces, that is, quotients of groups via subgroups \cite{dastte04,dastte04-1,daforastte08}. 
Examples of spaces that fall into this setup are modulation spaces on the sphere \cite{dastte04}, as well as $\alpha$-modulation
spaces \cite{daforastte08}. The latter were originally introduced by Feichtinger and Gr{\"o}bner as ``intermediate'' spaces (but {\it not} interpolation spaces)
between modulation spaces and Besov spaces \cite{fegr85,gr92-2}. In another direction, the first named author developed a coorbit theory
in the setup of spaces of functions obeying symmetries such as radiality \cite{ra05,ra05-6}. 
Here, one takes the set of residue classes of the locally compact group
modulo a symmetry group leading to a hypergroup structure. In concrete setups, the theory provides then frames of radial wavelets (that is, each
frame element is a radial function) for radial homogeneous Besov-Lizorkin-Triebel spaces, as well as, radial Gabor frames for radial modulation
spaces. Coorbit space theory can then be used to show compactness of certain embeddings when restricting modulation spaces
to radial functions \cite{ra05-2}.

As the next step, the first named author together with Fornasier realized that group theory is not needed at all in order to develop a coorbit space
theory \cite{fora05}. The starting point is now an abstract continuous frame \cite{alanga93}, which induces an associated transform. Then one measures
``smoothness'' via the norm of the transform in suitable function spaces on the index set of the continuous frame. Under certain integrability and
continuity properties of the continuous frame, again discrete Banach frames for the associated 
coorbit spaces can be derived via sampling of the continuous frame. All the setups of coorbit space theory mentioned above fall into this generalization
(except that the theory for quasi-Banach spaces still needs to be extended). The advantage of the group theoretical setup is only that some 
of the required conditions are automatically satisfied, while in this general context they enter as additional assumptions, which means 
that they have to be checked in a concrete situation.

While the theory in \cite{fora05} essentially applies only to coorbit spaces
with respect to weighted Lebesgue spaces, we extend this abstract theory in the present paper in order to treat a wider variety of coorbit spaces. Our main motivation is to cover inhomogeneous Besov-Lizorkin-Triebel spaces and generalizations thereof. Those spaces indeed do not fit into any of the group theoretical approaches which were available before. In order to handle them in full generality, one needs to take coorbits with respect to more complicated spaces rather than only weighted Lebesgue spaces. Indeed, we will need (weighted) mixed $L_{p,q}$ spaces. We derive characterizations of such general coorbits via discrete Banach frames and atomic decomposition, i.e., characterizations using discrete sequence spaces. Such are very useful in order to study embeddings, s-numbers, interpolation properties etc., because the structure of sequence spaces is usually much easier to investigate. 

We further treat the identification as coorbits of inhomogeneous Besov-Lizorkin-Triebel type spaces in detail (Section \ref{Peetre}). The application of our general abstract coorbit space theory from Section \ref{abstrth} leads to concrete atomic decompositions and wavelet characterizations of the mentioned spaces. Such discretizations have a certain history. A remarkable breakthrough in the theory was achieved by Frazier, Jawerth \cite{FrJa90} with the invention of the $\varphi$-transform. They fixed the notion of smooth atoms and molecules as building blocks for classical function spaces. Afterwards many authors have dealt with wavelet characterizations of certain generalizations of Besov-Lizorkin-Triebel spaces in the past. To mention all the relevant contributions to the subject would go beyond the scope of this paper. We rather refer to the monograph \cite[Chapt.\ 2, 3]{Tr06}, the references given there and to our overview Section \ref{overview}. 
Our results on wavelet basis characterizations in this paper rely on the abstract discretization result in Theorem \ref{wbases2} below, which allows 
to use orthogonal and even biorthogonal wavelets as well as tight (discrete) wavelet frames.

We are able to come up with a suitable definition of weighted  Besov-Lizorkin-Triebel spaces and their wavelet characterizations when the weight is only assumed to be doubling. Muckenhoupt $\mathcal{A}_p$-weights fall into this class of weights, but there exist doubling weights for which
a proper notion of Besov-Lizorkin-Triebel spaces was more or less unavailable before, although there exits certain attempts, see for instance \cite{Bo07}. 
In addition, we treat 
generalized $2$-microlocal spaces, Morrey-Besov-Lizorkin-Triebel spaces, and Besov-Lizorkin-Triebel spaces of dominating mixed smoothness and their characterizations
via wavelet bases. The treatment of spaces with variable integrability, or more general, with parameters $p,q,s$ depending on $x$, will be considered in a subsequent contribution. As another main feature we also provide a better way to identify Lizorkin-Triebel type spaces as coorbits. 
So far, the (homogeneous) Lizorkin-Triebel spaces have been identified as coorbits of so called tent spaces \cite{Gr88,Gr91} on the $ax+b$ group. However, tent spaces \cite{CoMeSt85} are rather complicated objects. In this paper, we proceed by introducing a Peetre type maximal function, 
related to the one introduced in \cite{Pe75}, as well as corresponding
function spaces on the index set of the (continuous wavelet) transform. Then Lizorkin-Triebel spaces can also be identified as 
coorbits with respect to these new spaces, which we call Peetre spaces. This was recently 
accomplished for the homogeneous spaces \cite{T10}. It turns out that Peetre spaces are much easier to handle than tent spaces. 
%For this identification we restrict to the Banach space case since a suitable general coorbit space theory for the treatment of quasi-Banach spaces is currently under construction. The authors expect similar results also in this situation.

In the present paper we restrict our considerations to coorbit space theory for Banach spaces. 
While an extension to the setting of quasi-Banach spaces is available
for classical coorbit space theory \cite{ra05-3}, such extension is more technical for general coorbit spaces, and currently under development. We expect
similar results also in this situation. In order to be well prepared for applying this generalized coorbit space theory for quasi-Banach spaces once it is developed in detail, we state certain characterizations of generalized Besov-Lizorkin-Triebel type spaces also for the quasi-Banach space 
cases $p,q < 1$ -- although we do not need such cases in the present contribution.

We hope to convince the reader with this paper that the abstract coorbit space theory is a very powerful tool, 
and allows a unified treatment of function spaces. In contrast, the theory leading to atomic decompositions, wavelet characterizations
of several newly introduced function spaces is often developed from scratch. We believe, that most of these spaces 
can be interpreted as coorbit spaces in our setting. Once this is established,
then one has to follow an easy recipe checking only basic properties, in order to come up with corresponding discrete characterizations.
These can be widely applied for approximation issues, to prove certain embeddings, interpolation formulas, etc.  

While our main focus in this paper is on inhomogeneous Besov-Lizorkin-Triebel spaces, we expect that the principles of the abstract coorbit space theory
apply also to other setups. To be more precise, we expect that our theory can be used to introduce also inhomogeneous shearlet spaces, and their atomic
decompositions, and $\alpha$-modulation spaces with different $p,q$-indices (the paper \cite{daforastte08} only treats the case $p=q$).

The paper is structured as follows. After setting some basic notation we give an overview over the main results
and achievements of the paper in Section \ref{overview}. Section \ref{abstrth} is devoted to the extension
of the abstract generalized coorbit space theory from \cite{fora05}. In Section \ref{Peetre} we apply this abstract theory
to the specific situation of coorbits with respect to Peetre type spaces. 
We study several examples in Section 5 and give concrete discretizations for generalized inhomogeneous Besov-Lizorkin-Triebel spaces of various type in terms of wavelet bases with corresponding sufficient conditions for admissible wavelets. Appendix A contains some basic facts concerning orthonormal wavelet bases on $\re$
and $\R$, in particular, %the definition of 
orthonormal spline wavelets.\\ % which is of certain interest in this paper.\\
\newline
{\bf Acknowledgement.} The authors acknowledge support by the Hausdorff Center for Mathematics, University of Bonn. In addition, they would like to thank Stephan Dahlke, Hans Feichtinger, Yoshihiro Sawano, Martin Sch\"afer, and Hans Triebel for valuable discussions, critical reading of preliminary versions of this manuscript and for several hints
how to improve it.

\subsection{Notation}
  To begin with we introduce some basic notation. The symbols
$\re, \re_+,\com, \N, \n$ and $\zz$ denote the real numbers, positive real numbers, complex
numbers, natural numbers, natural numbers including $0$ and the
integers. Let us emphasize that $\R$ has the usual meaning and $d$ is reserved
for its dimension. The elements are denoted by
$x,y,z,...$ and $|x|$ is used for the Euclidean norm. We use
$|k|_1$ for the $\ell_1^d$-norm of a vector $k$. Sometimes the notation $\bar{a}$ is used to indicate 
that we deal with vectors $\bar{a}= (a_1,...,a_d)$ taken from $\R$. The notation $\bar{a}>b$, where $b\in \re$, 
means $a_i>b$ for every $i=1,...,d$. If $X$ is a (quasi-)Banach space and $f\in X$ we use $\|f|X\|$ or
simply $\|f\|$ for its (quasi-)norm. The class of linear continuous mappings from 
$X$ to $Y$ is denoted by $\mathcal{L}(X,Y)$ or simply $\mathcal{L}(X)$ if $X=Y$. Operator 
(quasi-)norms of $A \in \mathcal{L}(X,Y)$ are denoted by $\|A:X\to Y\|$, or simply by
$\|A\|$. As usual, the letter $c$ denotes a constant, which may vary
from line to line but is always independent of $f$, unless the
opposite is explicitly stated. We also use the notation
$a\lesssim b$ if there exists a constant $c>0$ (independent of the
context dependent relevant parameters) such that $a \le c\,b$. If
$a\lesssim b$ and $b \lesssim a$ we write $a \asymp b$\,. For a real number $t$, we denote $t_+ = \max\{t,0\}$ and
$t_-=\min\{t,0\}$. The ball in $\R$ with center $x \in \R$ and radius $r>0$ is denoted by $B(x,r) = \{y \in \R, |x-y| \leq r\}$, while $|B(x,r)|$ is its volume.

\subsection{Lebesgue spaces and tempered distributions}

For a measure space $(X,\mu)$ and a positive measurable weight function $w:X\to \re$, we define
%The space $L^1(X,\mu)$ is the collection of all complex valued $\mu$-measurable functions $F$ on $X$ with finite norm
%\begin{equation}\label{L1norm}
%    \|F|L^1\|:=\int\limits_{X}|F(x)|d\mu(x)\,.
%\end{equation}
the space $L^w_p(X,\mu)$, $1\leq p<\infty$, as usual by,
%is defined analogously via replacing
%\eqref{L1norm} by
$$
    \|F|L^w_p(X,\mu)\|:=\Big(\int_{X}|w(\x)F(\x)|^p\,d\mu(\x)\Big)^{1/p} < \infty\,.
$$
A function $F$ belongs to $L^w_{\infty}(X,\mu)$ if and only if $Fw$ is essentially bounded
with respect to the measure $\mu$. If $w\equiv 1$ we simply write $L_p(X,\mu)$ instead of $L_p^w(X,\mu)$.
Moreover, the space $L_1^{loc}(X,\mu)$ contains all functions $F$ for which
the integral over all subsets of finite measure $K\subset X$ is finite. If $X = \R$ and the measure $\mu$ is the Lebesgue measure $dx$ then we write $L_p(\R)$.
For a measurable weight function $v>0$, the space $L_p(\R,v)$,
$0< p\leq \infty$, is the collection of all functions $F$ such that 
\begin{equation}\label{wLp}
    \|F|L_p(\R,v)\|:=\Big(\int_{\R}|F(x)|^pv(x)\,dx\Big)^{1/p} < \infty\,,
\end{equation}
i.e., it coincides with $L_p(X,\mu)$ where $X=\R$ and $d\mu(x) = v(x)dx$\,.

As usual $\mathcal{S}(\re^d)$ is used for
the locally convex space of rapidly decreasing infinitely
differentiable functions on $\re^d$ where its topology is generated
by the family of semi-norms
$$
    \|\varphi\|_{k,\ell} = \sup\limits_{x\in \re^d,|\bar{\alpha}|_1 \leq \ell}
    |D^{\bar{\alpha}}\varphi(x)|(1+|x|)^k\quad,\quad \varphi\in \mathcal{S}(\re^d)\,,
$$
where $k,\ell\in \n$. The space $\mathcal{S}'(\re^d)$, the topological dual of $\mathcal{S}(\R)$, is also referred to as the set of tempered distributions on $\re^d$. Indeed, a linear mapping $f:\mathcal{S}(\re^d) \to \com$ belongs
to $\mathcal{S}'(\re^d)$ if and only if there exist numbers $k,\ell \in \n$
and a constant $c = c_f$ such that
\begin{equation}\nonumber
    |f(\varphi)| \leq c_f\sup\limits_{x\in \R,|\bar{\alpha}|_1 \leq \ell}
    |D^{\bar{\alpha}}\varphi(x)|(1+|x|)^k
\end{equation}
for all $\varphi\in \mathcal{S}(\re^d)$. The space $\mathcal{S}'(\re^d)$ is
equipped with the weak$^{\ast}$-topology.

The convolution $\varphi\ast \psi$ of two
integrable (square integrable) functions $\varphi, \psi$ is defined via the integral
\begin{equation}\label{conv}
    (\varphi \ast \psi)(x) = \int_{\R} \varphi(x-y)\psi(y)\,dy\,.
\end{equation}
If $\varphi,\psi \in \mathcal{S}(\R)$ then \eqref{conv} still belongs to
$\mathcal{S}(\R)$. The convolution can be extended to $\mathcal{S}(\R)\times \mathcal{S}'(\R)$ via
$(\varphi\ast f)(x) = f(\varphi(x-\cdot))$. It is a pointwise defined $C^{\infty}$-function in $\R$ of at most polynomial growth. 

The Fourier transform defined on both $\mathcal{S}(\R)$ and $\mathcal{S}'(\R$)
is given by %$(\cf f)(\varphi) := f(\cf \varphi)$, 
$\widehat{f}(\varphi) := f (\widehat{\varphi})$, 
where
$\,f\in \mathcal{S}'(\R), \varphi \in \mathcal{S}(\R)$, and
$$
%\cf \varphi(\xi) 
\widehat{\varphi}(\xi):= (2\pi)^{-d/2}\int_{\R} e^{-ix\cdot
\xi}\varphi(x)\,dx.
$$
The Fourier transform is a bijection (in both cases) and its inverse is
given by $\varphi^{\vee} = \widehat{\varphi}(-\cdot)$.
%$\cf^{-1}\varphi = \cf\varphi(-\cdot)$.

\subsection{The continuous wavelet transform}
\label{sectCWT}
Many considerations in this paper are based on decay results for the continuous wavelet transform $W_g f(x,t)$.
A general reference for this notion is provided by the monograph \cite[2.4]{Dau92}.
In \cite[App.\ A]{T10} the second named author provided decay results based on the following setting. 
For $x\in \R$ and $t>0$ we define the unitary dilation and translation operators $\mathcal{D}^{L_2}_{t}$
and $T_x$ by
$$
  \mathcal{D}^{L_2}_{t}g := t^{-d/2}g\Big(\frac{\cdot}{t}\Big)\quad\mbox{and}\quad T_xg := g(\cdot-x)\quad,\quad g\in L_2(\R)\,.
$$
The wavelet $g$ is said to be the analyzing vector for a function $f\in L_2(\R)$. The 
continuous wavelet transform $W_gf$ is then defined by
$$
    W_g f(x,t) = \langle T_x\mathcal{D}^{L_2}_t g, f\rangle\quad,\quad x\in \R, t>0\,,
$$ 
where the bracket $\langle \cdot, \cdot \rangle$ denotes the
inner product in $L_2(\R)$. 
We call $g$ an admissible wavelet if 
$$
    c_g:= \int_{\R} \frac{|\widehat{g}(\xi)|^2}{|\xi|^d}\,d\xi < \infty\,.
$$
If this is the case, then the family $\{T_x\mathcal{D}^{L_2}_t g\}_{t>0, x\in \R}$ represents a tight
continuous frame in $L_2(\re)$ where $C_1 = C_2 = c_g$ (see Subsection \ref{confr}). 

The decay of the function $|W_g f(x,t)|$ mainly depends on the number of vanishing moments of the wavelet $g$ as well as on the smoothness of $g$
and the function $f$ to be analyzed, as is made precise in the following definition. % To be more precise we give the following definition. 

\begin{definition}\label{basedef} Let $L+1\in \n$, $K>0$. We define the properties $(D)$, $(M_L)$ and $(S_K)$ for
a function $f \in L_2(\R)$ as follows.
\begin{enumerate}
     \item[$(D)$] For every $N\in \N$ there exists a constant
    $c_N$ such that
    $$
        |f(x)| \leq \frac{c_N}{(1+|x|)^N}\,.
    $$
    \item[$(M_L)$] All moments up to order $L$ vanish, i.e., 
    $$
        \int_{\R}x^{\alpha}f(x)\,dx = 0
    $$
    for all $\alpha\in \n^d$ such that $|\alpha|_1 \leq L$\,.
    \item[$(S_K)$] The function
    $$
        (1+|\xi|)^{K}|D^{\alpha}\widehat{f}(\xi)|
    $$
    belongs to $L_1(\R)$ for every multi-index $\alpha\in \n^d$.
\end{enumerate}
\end{definition}
\noindent
Property $(S_K)$ is rather technical. Suppose we have a function $f \in C^{K+d+1}(\R)$ for some $K\in \N$ 
such that $f$ itself and all its derivatives satisfy $(D)$. The latter holds, for instance, if $f$
is compactly supported. Then this function satisfies $(S_K)$ by elementary 
properties of the Fourier transform. Conversely, if a function $g \in L_2(\R)$ satisfies $(S_K)$ for some $K>0$
then we have $g\in C^{\lfloor K \rfloor}(\R)$. However, in case of 
certain wavelet functions $\psi$ where the Fourier transform $\cf \psi$ is given explicitly (See App.\ \ref{splines}) we can verify $(S_K)$ directly. Depending on these conditions we state certain decay results for the function $|W_gf(x,t)|$ in Lemma \ref{help1} below.
 
\section{Overview on main results}
\label{overview}
As suggested in \cite{fora05}
coorbit space theory can be generalized to settings without group structure, and thereby allows the treatment of even more
function spaces via coorbit space theory. We follow this path and develop
the theory even further. Our main application are inhomogeneous Besov-Lizorkin-Triebel spaces with several generalizations and their wavelet characterizations. 

\subsection{Abstract coorbit space theory}
\label{confr}
In this section we give a brief overview before going into details later in Section \ref{abstrth}.
Assume $\mathcal{H}$ to be a separable Hilbert space and $X$ be a locally compact Hausdorff
space endowed with a positive Radon measure $\mu$ with $\supp \mu = X$. A family
$\mathcal{F} = \{\psi_\x\}_{\x\in X}$ of vectors in $\mathcal{H}$
is called a continuous frame if there exist constants $0<C_1,C_2<\infty$ such that
\begin{equation}\label{stab}
    C_1\|f|\mathcal{H}\|^2 \leq \int_{X} |\langle f,\psi_\x \rangle|^2 d\mu(\x) \leq C_2\|f|\mathcal{H}\|^2\quad \mbox{for all} \quad
    f\in \mathcal{H}\,.
\end{equation}
For the sake of simplicity, we assume throughout this paper that $\|\psi_\x|\mathcal{H}\| \leq C$, $\x\in X$, and that 
the continuous frame is tight, i.e., $C_1 = C_2$. After a possible
re-normalization we may assume that $C_1=C_2 = 1$. We note,
however, that non-tight frames appear also in several relevant examples and the associated coorbit theory is worked out
in \cite{fora05} -- at least to a significant extent. (The generalizations in this paper can also be developed in the setting
of non-tight frames.)

Associated to a continuous frame we define the transform $V = V_{\cf}: \mathcal{H} \to L_2(X,\mu)$ by
$$
V_{\cf}f(\x) = \langle f,\psi_\x \rangle\, \quad,\quad f \in \ch, \x \in X,
$$
and its adjoint $V^{\ast}_{\cf}:L_2(X,\mu) \to \mathcal{H}$,
$$
   V^{\ast}_{\cf} F = \int_X %\limits_{X}
   F(\y)\psi_\y\,d\mu(\y)\,.
$$
%We see immediately that $S = V_{\mathcal{F}}^{\ast}V_{\mathcal{F}}$.
Since we assume the frame $\cf$ to be tight, %$\cf$ %, 
i.e., $C_1=C_2=1$ in \eqref{stab}, 
the operator $V_{\mathcal{F}}^{\ast}V_{\mathcal{F}}$
is the identity. 
%Let us assume without 
%restriction $C_1 = C_2$ in \eqref{stab},
%then $S$ is indeed the identity. 
Hence,
\begin{equation}\label{eq34}
   f = \int_{X} (V_{\cf}f)(\y)\psi_\y\,d\mu(\y)\quad\mbox{and}\quad
   V_{\cf}f(\x) = \int_{X} V_{\cf}f(\y)\langle \psi_\y,\psi_\x\rangle \,d\mu(\y)\,.
\end{equation}
It follows from the tightness of the frame that $V_{\cf}^{\ast} V_{\cf}$ is a multiple of the identity, so that the
transform $V_{\cf}$ is invertible (more precisely, has a left-inverse). The second identity in \eqref{eq34} is the 
crucial reproducing formula $F = R(F)$ on the image of $\mathcal{H}$ under $V_{\cf}$, where
\[
 R(\x,\y) = R_{\mathcal{F}}(\x,\y) = \langle \psi_\y, \psi_\x \rangle \quad,\quad \x,\y \in X\,,
\] 
is an integral kernel (operator). The idea of coorbit space theory is to measure ``smoothness'' of $f$ via properties, i.e., suitable norms of the transform $V_{\cf}f$. 
Under certain integrability properties of the kernel $R(\x,\y)$, see \eqref{eq35} in Paragraph \ref{absFS},
one can introduce a suitable space $\ch_v^1$ of test functions and its dual $(\ch_v^1)^{\sim}$ (which plays the role of the tempered distributions in
this abstract context, see \eqref{def:H1v}), and extend the definition of the transform $V_{\mathcal{F}}$ to $(\ch_v^1)^{\sim}$ in \eqref{extend:VF}. 
Then associated to a solid Banach space $Y$ of locally integrable functions on $X$ 
(see Definition \ref{def:coorbit}), one defines the coorbit space
\[
\CoY = \{f \in (\ch_v^1)^{\sim}~:~V_{\cf} f \in Y\}\,, \quad \|f|\CoY\| := \|V_{\cf} f| Y\|,
\]
provided that, additionally, the kernel $R_{\cf}$ acts continuously from $Y$ into $Y$ as an integral operator. The latter 
is expressed as $R_{\cf}$ being contained in an algebra $\mathcal{B}_{Y,m}$ of kernels, see \eqref{def:BY}.
Then $\CoY$ is a Banach space, and one can show that ``similar'' frames (in the sense that their cross
Gramian kernel satisfies suitable integrability properties) define the same coorbit spaces, see Lemma \ref{ind1}. 

A key feature of coorbit space theory is the discretization machinery, which provides discrete frames, and characterizations of coorbit spaces $\CoY$
via suitable sequence spaces $Y^\sharp$ and $Y^\flat$. This is, of course, very useful because many properties, such as embeddings, s-numbers etc.,
are much easier to analyze for sequence spaces. Here, the starting point is a suitable covering $\cu = \{U_i\}_{i \in I}$ of the space $X$, of compact
subsets $U_i \subset X$. One defines
the $\cu$-oscillation kernel
\begin{equation}\nonumber
      \osc_{\mathcal{U}}(\x,\y) := \sup\limits_{\z\in Q_\y}|\langle \varphi_\x, \varphi_\y-\varphi_\z\rangle|
      %\\
    %= \sup\limits_{z\in Q_y}|R_{\cf}(x,y)-R_{\cf}(x,z)|\,,
\end{equation}
where $Q_\y = \bigcup_{\y\in U_i} %\limits_{y\in U_i}
U_i$. This kernel can be viewed as a sort of modulus of continuity associated to the frame $\cf$ and the covering $\cu$. 
If $\osc$ together with its adjoint $\osc^{\ast}$ is also contained in the algebra $\cb_{Y,m}$, see \eqref{def:BY}, then
one obtains a discrete Banach frame and atomic decompositions 
by subsampling the continuous frame at points $\x_i \in U_i$, that is, $\cf_d = \{\varphi_{\x_i} \}_{i \in I}$, see Theorem \ref{discr2} for details.
In particular, the coorbit space $\CoY$ is discretized by the sequence space $Y^\flat$ with norm
\[
\| \{\lambda_i\}_{i \in I} | Y^\flat\| = \Big\| \sum_{i \in I} | \lambda_i | \chi_{U_i} | Y \Big\|.
\]
Another new important key result of the abstract theory 
is that orthogonal and biorthogonal basis expansions, as well as tight frame expansions, 
where the basis / frame elements are sampled from a continuous frame, extend
automatically from the Hilbert space $\mathcal{H}$ to coorbit spaces under certain natural conditions,
see Theorem \ref{wbases2}. In addition, these basis / frame expansions characterize the respective coorbit space. 
For the concrete setup of characterizing generalized Besov-Lizorkin-Triebel spaces these ``natural conditions'' reduce to certain moment, decay and smoothness conditions (Definition \ref{basedef}) on the used wavelet and dual wavelet. In Section \ref{examples} we give sufficient conditions for the orthonormal wavelet characterization of several common as well as new generalizations of the inhomogeneous  Besov-Lizorkin-Triebel spaces. 
%The case of biorthogonal wavelet bases is completely analogous and left to the reader. 

The described setup indeed generalizes the original group theoretical one due to 
Feichtinger and Gr{\"o}chenig, see \cite{fora05} for details. For convenience of the reader we briefly summarize 
the main innovations and advances with respect to previous results. 

\subsubsection*{Main contribution and novelty}
\begin{itemize}
\item Feichtinger and Gr{\"o}chenig used group representations as an essential ingredient in
their initial work on coorbit space theory \cite{FeGr86,FeGr89a,FeGr89b,Gr91,fegr92-1}. 
The formulation of the theory in the present paper gets completely rid of group theory and
uses general continuous frames instead. This general approach was initiated in \cite{fora05}. 
The present paper even removes certain strong restrictions on the spaces $Y$ to treat a wider variety of coorbit spaces. 
\item We provide characterizations of general coorbit spaces $\Co Y$ by (discrete) 
Banach frames or atomic decompositions under suitable conditions (Theorem \ref{discr2}). 
There is an easy and explicit connection of the 
corresponding sequence space $Y^{\sharp}$ to the function space $Y$. This discretization machinery may 
be useful in situations, where it is even hard to construct a related basis or frame
for the underlying Hilbert space $\mathcal{H}$.
\item In several cases, an orthonormal basis, a Riesz basis or a tight frame for the Hilbert space $\mathcal{H}$ 
that arises from samples of a continuous frame
can be constructed directly via methods outside coorbit space theory. Then under natural
conditions on the continuous frame, Theorem \ref{wbases2} below 
shows that the corresponding expansions and characterizations
automatically extend to the coorbit spaces. This represents one of the core results in the present paper. 
It generalizes a result from classical coorbit
theory (associated to group representations) due to Gr{\"o}chenig \cite{Gr88}. In contrast
to his approach, our proof is independent from the discretization machinery in Theorem \ref{discr2}.
\item Our extended coorbit theory allows
to identify a large class of function spaces as coorbits. Therefore, the abstract
discretization machinery is available to such function spaces. We emphasize that due to this unified
approach, the theory leading to atomic decompositions for several classes of
spaces does not have to be developed from scratch over and over again 
for each new class of function spaces. From this point of view there are numerous previous 
results on atomic decompositions, which are partly recovered as well as extended by our theory. 

With a similar intention Hedberg and Netrusov gave an axiomatic approach to function spaces of Besov-Lizorkin-Triebel type in their substantial paper \cite{HeNe07}. Their approach is different from ours but also leads to atomic decompositions in a unified way. In a certain sense our approach is more flexible since the abstract theory in Section \ref{abstrth} is also applicable to, e.g., the recent shearlet spaces \cite{dakustte09} as well as modulation spaces.
\end{itemize}

\subsection{Inhomogeneous Besov-Lizorkin-Triebel type spaces}  
\label{overBe}
In order to treat inhomogeneous spaces of Besov-Lizorkin-Triebel type, see \cite{Tr83,Tr92,Tr06} and 
the references given there, we introduce the index set $X = \R \times [(0,1) \cup \{\infty\}]$, where ``$\infty$'' denotes an isolated point, and  define the Radon measure
$\mu$ by
$$
    \int_{X} F(\x) d\mu(\x) = \int_{\R}\int_{0}^1 F(x,s) \frac{ds}{s^{d+1}}dx + \int_{\R} F(x,\infty) dx\,.
$$
The main ingredient is a solid Banach space $Y$ of functions on $X$. We use two general scales $P^w_{B,q,a}(X)$ and $L^w_{B,q,a}(X)$ of spaces on $X$. Here, we have $1\leq q\leq \infty$ and $a>0$. The parameter $B=B(\R)$ is a solid space of measurable functions on $\R$ in the sense of Definition \ref{condB} below, for instance, a weighted Lebesgue space.
%The letter ''P'' refers to Peetre's maximal function, see Definition \ref{defFS} below.
For a function $F: X \to \C$, the Peetre type maximal function $\mathcal{P}_aF$ defined on $X$ is given by
\begin{align}
    \mathcal{P}_a F(x,t) &:= \sup\limits_{z\in \R}\frac{|F(x+z,t)|}{(1+|z|/t)^a}\quad,\quad x\in \R, 0<t<1\,,\label{def:Peetre:max}\\
    \mathcal{P}_a F(x, \infty) &:= \sup\limits_{z\in \R} \frac{|F(x+z,\infty)|}{(1+|z|)^a} \quad , \quad x \in \R.\notag
\end{align}
The function $w: X \to \mathbb{R}_+$ is a weight function satisfying the technical growth conditions (W1) and (W2) in Definition \ref{def:weight}. 
Then the Peetre spaces and Lebesgue spaces are defined as
\begin{equation}\nonumber
  \begin{split}
     P^{w}_{B,q,a}(X) &:= \{F:X \to \C~:~\|F|P^{w}_{B,q,a}\| < \infty\},\\
     L^{w}_{B,q,a}(X) &:= \{F:X \to \C~:~\|F|L^{w}_{B,q,a}\| < \infty\},\\
  \end{split}   
\end{equation}
with respective norms
\begin{align}
   \|F|P^{w}_{B,q,a}\| & := \Big\|w(\cdot,\infty)\cp_a F(\cdot, \infty) |B(\R)\Big\| 
   + \Big\|\Big(\int_{0}^1 \Big[w(\cdot,t)\cp_a F(\cdot, t) \Big]^q\frac{dt}{t^{d+1}}\Big)^{1/q}|B(\R)\Big\|,\label{def:Peetre:Pspace}\\
    \|F|L^{w}_{B,q,a}\| & := \Big\|w(\cdot,\infty)\cp_a F(\cdot, \infty) |B(\R)\Big\| 
    + \Big(\int_{0}^1 \Big\|w(\cdot,t) \cp_a F(\cdot,t)|B(\R)\Big\|^q\frac{dt}{t^{d+1}}\Big)^{1/q}\,.\label{def:Peetre:Lspace}
\end{align}
%and 
%\begin{equation}\label{eq-10}
%  \begin{split} 
%   \|F|L^{w}_{B,q,a}\| := &\Big\|w(\cdot,\infty)\sup\limits_{z\in \R}\frac{|F(\cdot+z,\infty)|}{(1+|z|)^a}|B(\R)\Big\|\\
%    &+ \Big(\int_{0}^1 \Big\|w(\cdot,t)\sup\limits_{z\in \R}
%   \frac{|F(\cdot+z,t)|}{(1+|z|/t)^a}|B(\R)\Big\|^q\frac{dt}{t^{d+1}}\Big)^{1/q}\,.
%  \end{split} 
%\end{equation}
We give the definition of an admissible continuous
frame $\cf$ on $X$. 
  
\begin{definition}\label{admfr} A continuous (wavelet) frame $\cf = \{\varphi_\x\}_{\x\in X}$, $X = \R \times [(0,1) \cup \{\infty\}]$, is admissible if it is of the form 
\begin{equation}\label{eq30}
      \varphi_{(x,\infty)} = T_x \Phi_0\quad \mbox{and} \quad \varphi_{(x,t)} = T_x\mathcal{D}^{L_2}_t \Phi\,,
\end{equation}
where $\Phi$ denotes a radial function from $\mathcal{S}(\R)$ 
satisfying $\widehat{\Phi} > 0$ on $\{x:1/2< |x| < 2\}$ and 
$$
    \int_{\R} \frac{|\widehat{\Phi}(\xi)|^2}{|\xi|^d}\,d\xi = 1\,.
$$
We further assume that 
$\Phi$ has infinitely many vanishing moments (see Definition \ref{basedef}). %satisfies $(M_L)$ for all $L\in \N$, see Definition \ref{basedef}\,. 
This condition is satisfied, for instance, 
if $\widehat{\Phi}$ vanishes on $\{x:|x|<1/2\}$. The function $\Phi_0 \in \mathcal{S}(\R)$ is chosen such that 
$$
    |\widehat{\Phi}_0(\xi)|^2 + \int_{0}^1 |\widehat{\Phi}(t\xi)|^2\frac{dt}{t} = 1\,.
$$
\end{definition}
The functions $\Phi$ and $\Phi_0$ from Definition $\ref{admfr}$ satisfy $(D)$ and $(S_K)$ for every $K>0$. Additionally,
$\Phi$ satisfies $(M_L)$ for any $L\in \N$. Moreover, the continuous frame \eqref{eq30} represents a tight continuous frame in the sense of \eqref{stab}. Indeed, we apply Fubini's and Plancherel's theorem to get 
\begin{equation}
  %\begin{split}
    \|f|L_2(\R)\|^2 = \int_{\R}\Big(|\langle f,\varphi_{(x,\infty)}\rangle|^2+
    \int_{0}^1|\langle f,
    \varphi_{(x,t)}\rangle|^2\frac{dt}{t^{d+1}}\Big)\,dx = \int_{X} |\langle f,\varphi_{\x}\rangle|^2 d \mu(\x)\,.\notag
  % \end{split}
\end{equation}
The transform $V_{\cf}$ on $\mathcal{H}= L_2(\R)$ is then given by 
$V_{\cf}f(\x) = \langle f, \varphi_\x\rangle$, $\x\in X$.

With these ingredients at hand, the associated coorbit spaces are given as
\begin{align}
\CoP &:= \Co(P^{w}_{B,q,a}, \cf) = \{f \in \mathcal{S}'~:~ V_\cf f \in  P^{w}_{B,q,a}(X)\},\notag\\
\CoL &:= \Co(L^{w}_{B,q,a}, \cf) = \{f \in \mathcal{S}'~:~ V_\cf f \in L^{w}_{B,q,a}(X)\}.\notag
\end{align}
The spaces $\CoL$ can be interpreted as generalized Besov spaces, while the spaces $\CoP$ serve as generalized
Lizorkin-Triebel spaces. Below we use the abstract machinery of coorbit space theory to show that these are Banach spaces,
and we provide characterizations by wavelet bases, in particular, by orthonormal spline wavelets, see Appendix \ref{splines}. 
We will recover known and new spaces, as well as known and new wavelet characterizations.

We shortly give some examples.
\begin{itemize}
\item {\bf Classical inhomogeneous Besov and Lizorkin-Triebel spaces}. Here we take 
$1\leq p,q \leq\infty$, $\alpha\in \re$, $w_\alpha(x,t) = t^{-\alpha-d/2+d/q}$,
$B(\R) = L_p(\R)$ and $a > d/\min\{p,q\}$. 
Then
\begin{equation}\label{eq52}
\Co P^{w_{\alpha}}_{B,q,a} = F^{\alpha}_{p,q}(\R), \quad \Co L^{w_{\alpha}}_{B,q,a} = B^{\alpha}_{p,q}(\R),
\end{equation}
where $F^{\alpha}_{p,q}(\R)$ is the classical Lizorkin-Triebel space and $B^{\alpha}_{p,q}(\R)$ the classical Besov space, see \cite{Tr83}.
\item {\bf Weighted inhomogeneous Besov and Lizorkin-Triebel spaces.} For a doubling weight $v$ 
 and $B(\R) = L_p(\R,v)$ we obtain
\[
\Co P^{w_{\alpha}}_{B,q,a} = F^{\alpha}_{p,q}(\R,v), \quad \Co L^{w_{\alpha}}_{B,q,a} = B^{\alpha}_{p,q}(\R,v).
\]
Note, that there are doubling weights which do not belong to the Muckenhoupt class $\mathcal{A}_{\infty}$. We provide 
a reasonable definition of the respective spaces (Definition \ref{inhomdoubling}) and atomic decompositions also in this situation, see
Section \ref{doubling}. 
\item {\bf Generalized $\mathbf{2}$-microlocal spaces.}
The identities in \eqref{eq52} remain valid if we replace the spaces $F^{\alpha}_{p,q}(\R)$ and 
$B^{\alpha}_{p,q}(\R)$ by the generalized $2$-microlocal spaces $F^{w}_{p,q}(\R)$ and 
$B^{w}_{p,q}(\R)$ where $w \in \mathcal{W}^{\alpha_3}_{\alpha_1,\alpha_2}$ is an admissible $2$-microlocal weight, see Definition \ref{def:weight}.
\item {\bf Besov-Lizorkin-Triebel-Morrey type spaces.} By putting $B(\R) = M_{u,p}(\R)$, where the latter represents a
Morrey space, we obtain a counterpart of \eqref{eq52} also for Besov-Lizorkin-Triebel-Morrey spaces.
\end{itemize}
Furthermore, with a slightly different setup we also treat 
\begin{itemize}
\item {\bf Besov-Lizorkin-Triebel spaces of dominating mixed smoothness.} If $1\leq p,q\leq \infty$, $\bar{r}\in \R$, and 
$\bar{a} > 1/\min\{p,q\}$ we will show
$$
  \Co P^{\bar{r}+1/2-1/q}_{p,q,\bar{a}} = S^{\bar{r}}_{p,q}F(\R),\quad 
  \Co L^{\bar{r}+1/2-1/q}_{p,q,\bar{a}} = S^{\bar{r}}_{p,q}B(\R),
$$
where $S^{\bar{r}}_{p,q}F(\R)$, $S^{\bar{r}}_{p,q}B(\R)$ are Lizorkin-Triebel and Besov spaces of mixed dominating smoothness,
see e.g. \cite{ScTr87,Vy06}.
\end{itemize}
Let us summarize the innovations and main advances of our considerations with respect to the theory of Besov-Lizorkin-Triebel spaces. 

\subsubsection*{Main contribution and novelty}
\begin{itemize}

\item We work out in detail the application of the general abstract coorbit theory to inhomogeneous function
spaces of Besov-Lizorkin-Triebel type. We further create an easy recipe for finding concrete 
atomic decompositions which is applicable to numerous examples of well-known spaces on the one hand and new generalizations
on the other hand. 
 
In Section \ref{Peetre} we give a very general definition of the family of 
Besov-Lizorkin-Triebel type spaces as coorbits of Peetre type spaces, see \eqref{def:Peetre:Pspace} and
\eqref{def:Peetre:Lspace}.
These depend on a weight function $w$ on $X$ and a Banach space $B$ on $\R$. 
Indeed, our conditions on $w$ and $B$ (Definitions \ref{def:weight} and \ref{condB}) 
are rather general but, however, allow for introducing and analyzing corresponding coorbit spaces. 
To the best knowledge of the authors Besov-Lizorkin-Triebel spaces have not yet been introduced in this generality.
\item In the classical literature on coorbit spaces, tent spaces \cite{CoMeSt85} are used to identify
homogeneous Lizorkin-Triebel spaces as coorbits. While tent spaces are rather complicated objects, our newly
introduced Peetre type spaces are much easier to handle. Their structure \eqref{def:Peetre:Pspace}, \eqref{def:Peetre:Lspace} 
allows for the definition of inhomogeneous spaces. Indeed, 
combined with Proposition \ref{bkernel} this represents one of the core ideas in the present paper. 

\item The conditions on $w$ and the space $B(\R)$ in Definitions \ref{def:weight} and \ref{condB} 
involve parameters $\alpha_1,\alpha_2,\alpha_3,\delta_1,\delta_2,\gamma_1,\gamma_2$. 
We identify explicit conditions on the smoothness $K$ and number $L$ of vanishing moments 
(see Definition \eqref{basedef})
of wavelets in terms of these parameters, which allow to provide characterization
of the generalized Besov-Lizorkin-Triebel spaces via wavelet bases (Theorem \ref{wbasesP}).
While we state the result only for orthonormal wavelet bases it easily extends to biorthogonal
wavelets. The corresponding sequence spaces are studied in detail.

\item In Section \ref{examples} we identify several known generalizations of inhomogeneous 
Besov-Lizorkin-Triebel type spaces as coorbits and generalize even further, see Theorems \ref{mainex}, \ref{coo:doubling},
\ref{coo:morrey}, and \ref{coo:dom}. This requires 
some effort since the spaces are usually not given in terms of continuous characterizations, see \cite{T10}.
We provide these characterizations along the way, see Proposition \ref{contchar} and the paragraphs 
before Theorems \ref{coo:doubling}, \ref{coo:morrey}.
In particular, our analysis includes
classical Besov-Lizorkin-Triebel spaces,
2-microlocal Besov-Lizorkin-Triebel spaces with Muckenhoupt weights and
2-microlocal Besov-Lizorkin-Triebel-Morrey spaces. Moreover, we introduce 
Besov-Lizorkin-Triebel spaces with doubling weights, which are not necessarily Muckenhoupt. 
In the latter case, a ``classical'' definition is not available (but see \cite{Bo07}), 
and we emphasize that coorbit space theory
provides a natural approach for such spaces as well.

\item Special cases of our result concerning wavelet bases characterizations of Muckenhoupt weighted $2$-microlocal Besov-Lizorkin-Triebel spaces in Theorems \ref{wbasesmucken}, \ref{kempka1}, \ref{classbtl} already appeared in the literature. Indeed, see Theorem 3.10\ in \cite{HaPi08}, Theorem 4 in \cite{Ke10}, Theorem 1.20 in \cite{Tr08}, Theorem 3.5 in \cite{Tr06}, or Propositions 5.1, 5.2 in \cite{Ky03}. Our result concerning decompositions of Morrey type spaces, Theorem \ref{morreydec}, has a special case in \cite{Sa08} and in the recent monograph \cite[Thm 4.1]{YuSiYa10}. In the mentioned references the conditions on smoothness and cancellation (moment conditions) are often slightly less restrictive than ours for this particular case. But this fact might be compensated by the unifying nature of our approach.
However, compared to the conditions in \cite[Thm.\ 3.5]{Tr06} our restrictions in Theorem \ref{classbtl} are similar. Concerning characterizations of classical Besov spaces with orthonormal spline wavelets, see Appendix \ref{splines}, we refer to \cite{Bo95} for the optimal conditions with respect to the order $m$. 

\item For technical reasons several authors restrict to compactly supported atoms \cite[Sect.\ 3.1.3]{Tr06}, \cite[Sect.\ 1.2.2]{Tr08}, \cite[Sect.\ 2.2, 2.4]{Vy06}, \cite[Thm.\ 4.1]{YuSiYa10}, especially to wavelet decompositions using the well-known compactly supported but rather complicated Daubechies wavelet system \cite{Dau92}. In the literature more general atoms are called 
{\it molecules}. This term goes back to Frazier, Jawerth \cite[Thm.\ 3.5]{FrJa90}. Several authors \cite{Ky03, Ke10, Sa08, YuSiYa10, Bo07} used their techniques in order to generalize results in certain directions. In this sense our approach is already 
sufficiently general because we allow arbitrary orthonormal (biorthogonal) wavelets having sufficiently large smoothness, vanishing moments, and decay. 

\item By a slight variation of the setup of Section \ref{Peetre} we also identify inhomogeneous
Besov-Lizorkin-Triebel spaces of dominating mixed smoothness as coorbit spaces and derive
corresponding wavelet characterizations with explicit smoothness and moment conditions on
the wavelets (Theorem \ref{wbasestensor}). 
In contrast to most previous results \cite{Vy06, Tr09}, 
we are not restricted to compact support.  In particular, we obtained characterizations via orthonormal spline wavelets in Corollary \ref{corsplines}, which are comparable (with respect to the order of the splines) to the very recent results in the monograph \cite[Sect.\ 2.5]{Tr09}. Furthermore, since our arguments are based on the abstract Theorem \ref{wbases2}, our results extend in a straightforward way to discretization results using numerically convenient biorthogonal wavelets \cite{CoDaFe92}. 
\end{itemize}

\subsection{Further extensions and applications}
We conclude this section with a list of further possible extensions and applications of our work.

\begin{itemize}
\item The discrete wavelet characterizations derived in this paper allow to reduce many questions
on function spaces to related questions on the associated sequence spaces. For instance,
the study of embeddings or the computation of certain widths such as entropy, (non-)linear 
approximation, Kolmogorov, Gelfand,..., are much more straightforward by using the stated sequence 
space isomorphisms. This can be seen as the major application of our theory.

\item Our abstract approach would clearly allow to incorporate further extensions of 
Besov-Lizorkin-Triebel type spaces. For instance, one might think of coorbits with respect to 
(weighted) Lorentz spaces or (weighted) Orlicz spaces, or one may introduce weights also
in spaces of dominating mixed smoothness. Another recent development considers
variable exponents where $p,q$ are not constant but actually functions of the space variable.
The general theory would then provide also wavelet characterizations of such spaces. 

\item The abstract approach allows to handle also function spaces of different type than
Besov-Lizorkin-Triebel spaces, such as (inhomogeneous) shearlet spaces or
modulation spaces. For instance, it would be interesting to work out details for 
modulation spaces with Muckenhoupt weights \cite{sa08-1}.

\item The abstract coorbit space theory in the present stage applies only to Banach spaces.
An extension to quasi-Banach spaces, similar to the classical case in \cite{ra05-3}, 
is presently under investigation.
\end{itemize}

\section{General coorbit space theory}
\label{abstrth}
The classical coorbit space theory due to Feichtinger and Gr{\"o}chenig
\cite{FeGr86, Gr88, FeGr89a, FeGr89b, Gr91} can be generalized in various ways.
One possibility is to replace the locally compact group $G$ by a locally compact Hausdorff space $X$
without group structure equipped with a positive Radon measure $\mu$ that replaces the
Haar measure on the group \cite{fora05}. % This was done by Fornasier and Rauhut in \cite{fora05}. 
This section is intended to recall all the relevant background from \cite{fora05} and to extend the available abstract theory.

\subsection{Function spaces on $X$}
\label{absFS}
%We start by introducing the Lebesgue spaces on $X$.
%The space $L^1(X,\mu)$ is the collection of all complex valued $\mu$-measurable functions $F$ on $X$ with finite norm
%\begin{equation}\label{L1norm}
%    \|F|L^1\|:=\int_{X}|F(x)|d\mu(x)\,.
%\end{equation}
%The space $L^p(X,\mu)$, $1\leq p<\infty$,
%is defined analogously via replacing
%\eqref{L1norm} by
%$$
%    \Big(\int_{X}|F(x)|^p\,d\mu(x)\Big)^{1/p} < \infty\,.
%$$
%A function $F$ belongs to $L^{\infty}(X,\mu)$ if and only if it is essentially bounded
%with respect to the measure $\mu$.
%Moreover, the space $L^1_{loc}(X,\mu)$ contains all functions $F$ for which
%the integral over all compact subsets of $K\subset X$ is finite.
%
%Fix a weight function $v\geq 1$. The space $L^p_v(X,\mu)$,
%$1\leq p\leq \infty$, is the collection of all function $F$ such that $F\cdot v \in L^p(X,\mu)$.
In order to define the coorbit space with respect to a Banach space $Y$ of functions on $X$ we need to require
certain conditions on $Y$.

\begin{description}
 \item$(Y)$ The space $(Y,\|\cdot|Y\|)$ is a non-trivial Banach space of functions on $X$ that is
 contained in $L_1^{loc}(X,\mu)$ and satisfies the solidity condition, i.e.,
 if $F$ is measurable and $G \in Y$ such that $|F(\x)| \leq |G(\x)|$ a.e., then $F\in Y$ and $\|F|Y\| \leq \|G|Y\|$\,.
\end{description}
This property holds, for instance, for weighted $L^w_p(X,\mu)$-spaces.  % all the above defined function spaces.
The classical theory by Feichtinger and Gr{\"o}chenig \cite{FeGr86, Gr88, Gr91} heavily uses the group convolution.
Since the index space $X$ does not possess a group structure in general we have to find a proper
replacement for the convolution of functions on a group. Following \cite{fora05} we use
integral operators with kernels belonging to certain kernel algebras. Let
\begin{equation}\label{eq35}
  \mathcal{A}_1 := \{K:X \times X \to \C~:~ K \mbox{ is measurable and } \|K|\mathcal{A}_1\| < \infty\}\,,
\end{equation}
where
\begin{equation}\label{def:A:norm1}
    \|K|\mathcal{A}_1\| := \max\Big\{\esssup{\x\in X}\int_{X}|K(\x,\y)|d\mu(\y)~,~ \esssup{\y\in X} \int_{X}|K(\x,\y)|d\mu(\x)\Big\}\,.
\end{equation}
The sub-index $1$ indicates the unweighted case.
We further consider weight functions $v(\x)\geq 1$ on $X$.
The associated weight $m_v$ on $X \times X$ is given by
\begin{equation}%\nonumber
     m_v(\x,\y) := \max\Big\{\frac{v(\x)}{v(\y)}, \frac{v(\y)}{v(\x)} \Big\}\quad,\quad \x,\y\in X\,.\label{def:mv}
\end{equation}
For a weight $m$ on $X \times X$ the corresponding sub--algebra $\mathcal{A}_m \subset \mathcal{A}_1$ is defined as
$$
    \mathcal{A}_m := \{K:X\times X \to \com~:~Km \in \mathcal{A}_1\}
$$
endowed with the norm
$$
    \|K|\mathcal{A}_m\| := \|Km|\mathcal{A}_1\|\,.
$$
Later we will need that the kernel $R(\x,\y)$ from Subsection \ref{confr} and further related kernels
(see Subsection \ref{sectdiscr}) belong to $\mathcal{A}_m$ for a proper weight function $m$. In order to define
the coorbit of a given function space $Y$ we will further need that these particular kernels act boundedly
from $Y$ to $Y$, i.e., the mapping
$$
    K(F) = \int\limits_X K(\cdot,y)F(y)\,d\mu(y)
$$ 
is supposed to be bounded. It is easy to check that the condition $K \in \mathcal{A}_m$ is sufficient for $K$ to map $Y=L_p^{v}(X)$ into $L_p^{v}(X)$ boundedly. This, however, is not the case in general and has to be checked for particular spaces $Y$.  At this point we modify the setting in \cite{fora05} according to Remark 2 given there.
Associated to a space $Y$ satisfying $(Y)$ and a weight $m$ we introduce the subalgebra
$$
    \mathcal{B}_{Y,m} := \{K:X\times X \to \mathbb{C}\,:\,K\in \mathcal{A}_m\mbox{ and }K\mbox{ is bounded from }Y \mbox{ into } Y\}\,,
$$
where
\begin{equation}\label{def:BY}
    \|K|\mathcal{B}_{Y,m}\| := \max\{\|K|\mathcal{A}_m\|, \|K|Y\to Y\|\}
\end{equation}
defines its norm. 

\subsection{Associated sequence spaces}
Let us start with the definition of an admissible covering of the index space $X$.

\begin{definition} A family $\mathcal{U} = \{U_i\}_{i\in I}$ 
of subsets of $X$ is called admissible covering of $X$, if the following conditions are satisfied.
\begin{description}
\item(i) Each set $U_i, i\in I$, is relatively compact and has non-void interior.
\item(ii) It holds $X = \bigcup\limits_{i\in I} U_i$\,.
\item(iii) There exists some constant $N>0$ such that
\begin{equation}\label{def:N:overlap}
  \sup\limits_{j\in I}\sharp\{i\in I, U_i \cap U_j \neq \emptyset\} \leq N <\infty\,.
\end{equation}
\end{description}
Furthermore, we say that an admissible covering $\mathcal{U} = \{U_i\}_{i\in I}$ is moderate with respect to $\mu$, 
if it fulfills the following additional assumptions.
\begin{description}
  \item(iv) There exists some constant $D>0$ such that $\mu(U_i)\geq D$ for all $i\in I$.
  \item(v) There exists a constant $\tilde{C}$ such that
  $$
      \mu(U_i) \leq \tilde{C}\mu(U_j)\quad,\quad \mbox{for all }i,j \mbox{ such that } U_i\cap U_j \neq \emptyset\,.
  $$
\end{description}
\end{definition}
Based on this framework, we are now able to define sequence 
spaces associated to function spaces $Y$ on the set $X$ with 
respect to the covering $\mathcal{U}$.

\begin{definition}\label{absdefss} Let $\mathcal{U} = \{U_i\}_i$ be an admissible
covering of $X$ and let $Y$ be a Banach function space satisfying
$(Y)$, which contains all the characteristic functions $\chi_{U_i}$.
We define the sequence spaces $Y^{\flat}$ and $Y^{\sharp}$ associated to
$Y$ as
\begin{equation}\nonumber
  \begin{split}
    Y^{\flat} = Y^{\flat}(\mathcal{U}) &:= \Big\{\{\lambda_i\}_{i\in I}~:~
    \|\{\lambda_i\}_{i\in I}|Y^{\flat}\| := \Big\|\sum\limits_{i\in I}
    |\lambda_i|\chi_{U_i}|Y\Big\|<\infty
    \Big\}\,,\\
    Y^{\sharp} = Y^{\sharp}(\mathcal{U}) &:= \Big\{\{\lambda_i\}_{i\in I}~:~
    \|\{\lambda_i\}_{i\in I}|Y^{\sharp}\| := \Big\|\sum\limits_{i\in I}
    |\lambda_i|\mu(U_i)^{-1}\chi_{U_i}|Y\Big\|<\infty
    \Big\}\,.
  \end{split}
\end{equation}
\end{definition}

\begin{remark} Under certain conditions on the families $\mathcal{U} = \{U_i\}_{i\in I}$ and $\mathcal{V} = \{V_i\}_{i\in I}$
  over the same index set $I$, the sequence spaces
  $Y^{\sharp}(\mathcal{U})$ and $Y^{\sharp}(\mathcal{V})$ coincide (similar for $Y^\flat(\mathcal{U})$ and
  $Y^\flat(\mathcal{V})$), see Definition 7 and Lemma 6 in \cite{fora05}\,.
\end{remark}
The following lemma states useful properties of these sequence spaces.

\begin{lemma} Let $\mathcal{U} = \{U_i\}_i$ be an admissible
covering of $X$ and let $Y$ be a Banach function space satisfying
$(Y)$ which contains all the characteristic functions $\chi_{U_i}$.
\begin{description}
%Assume further
%\begin{description}
\item(i) If there exist constants $C,c>0$ such that $c\leq \mu(U_i) \leq C$
for all $i\in I$ then the spaces $Y^{\sharp}$ and $Y^{\flat}$ coincide in the
   sense of equivalent norms.
\item(ii) If for all $i\in I$ the relation $\|\chi_{U_i}|Y\| \lesssim v_i$ holds, 
where $v_i = \sup_{x\in U_i} v(x)$, then we have the continuous embeddings
   $$
        \ell_1^{v_i} \hookrightarrow Y^{\flat} \hookrightarrow Y^{\sharp}\,.
   $$
\end{description}
\end{lemma}

\bproof The statement in (i) is immediate. The first
embedding in (ii) is a consequence of the triangle inequality in $Y$,
indeed
$$
    \Big\|\sum\limits_{i\in I}|\lambda_i| \chi_{U_i}|Y\Big\| \leq
    \sum\limits_{i\in I}|\lambda_i|\|\chi_{U_i}|Y\| \lesssim
    \sum\limits_{i\in I}|\lambda_i|v_i\,.
$$
The second embedding is a consequence of the fact that $\mathcal{U}$ is an admissible 
covering.\eproof

\subsection{Coorbit spaces}
We now introduce properly our coorbit spaces and show some of their basic properties.
To this end we fix a space $Y$ satisfying $(Y)$, a weight $v \geq 1$, and a tight continuous frame 
$\cf = \{\psi_\x\}_{\x\in X} \subset \mathcal{H}$ which satisfies the following property $(F_{v,Y})$\,.
%a weight function $v(x)\geq 1$ and tight continuous
%frame $\mathcal{F} = \{\psi_x\}_{x\in X} \subset \mathcal{H}$
%satisfying
\begin{description}
  \item$(F_{v,Y})$ The image space $R_{\mathcal{F}}(Y)$ is continuously
 embedded into $L_{\infty}^{1/v}(X,\mu)$ and $R_{\mathcal{F}}$ belongs to the algebra $\mathcal{B}_{Y,m}$, where $m$ is the weight
 on $X \times X$ associated to $v$ via \eqref{def:mv}.
 %Let $m$ be the associated to $v$ on $X\times X$. The kernel . Moreover, 
 \end{description}
The embedding $R_{\mathcal{F}}(Y) \hookrightarrow L_{\infty}^{1/v}(X,\mu)$ might seem 
a bit strange at first glance. However, we will return to that point later on and 
reduce this question to conditions on the frame $\cf$ and the sequence space associated to $Y$.
The property $(\mathcal{F}_{v,Y})$ sets us in the position to define the coorbit space $\CoY = \Co (\mathcal{F},Y)$.
We first define the reservoir
\begin{equation}\label{def:H1v}
    \mathcal{H}^1_v = \{f\in \mathcal{H}~:~V_{\cf}f \in L_1^v(X,\mu)\}\,
\end{equation}
endowed with the norm
$$
   \|f|\mathcal{H}^1_v\| = \|V_{\cf}f|L_1^v\|\,.
$$
The space $\mathcal{H}^1_v$ is a Banach space, see \cite{fora05}. 
By $R_{\mathcal{F}} \in \mathcal{B}_{Y,m} \subset \mathcal{A}_m$ 
we see immediately that $\psi_\x \in \mathcal{H}^1_v$ for all $\x \in X$\,.
We denote by $(\mathcal{H}^1_v)^{\sim}$ the canonical anti-dual of $\mathcal{H}^1_v$. 
We may extend the transform $V$ to $(\mathcal{H}^1_v)^{\sim}$
by 
\begin{equation}\label{extend:VF}
      (V_{\cf}f)(\x) = f(\psi_\x)\quad,\quad \x\in X, f\in (\mathcal{H}^1_v)^{\sim}\,.
\end{equation}
The reproducing formula still holds true. If $F = V_{\cf}f$ for $f\in (\mathcal{H}^1_v)^{\sim}$ then 
$R_{\cf}(F) = F$. Conversely, if $F \in L_{\infty}^{1/v}$ satisfies the reproducing formula 
$F = R_{\cf}(F)$ then there exists an $f\in (\mathcal{H}^1_v)^{\sim}$ such that $F = V_{\cf}f$. For more details see \cite[Sect.\ 3]{fora05}.

Now we are able to give the crucial definition of the coorbit space $\Co Y$.

\begin{definition}\label{def:coorbit} Let $Y$ be a Banach function space on $X$ satisfying $(Y)$.
Let further $\cf = \{\psi_\x\}_{\x\in X}$ be a tight continuous frame on $X$ with property $(\mathcal{F}_{v,Y})$.
The coorbit $\Co(\mathcal{F},Y)$ of $Y$ with respect to $\mathcal{F}$ is given by
  $$
      \CoY = \Co(\cf, Y):= \{f\in (\mathcal{H}^1_v)^{\sim}~:~V_{\mathcal{F}} f \in Y\}\quad\mbox{with}\quad
      \|f|\CoY\| = \|V f|Y\|\,.
  $$
\end{definition}
For proofs of the following properties we refer to \cite{fora05}. As a consequence of property $(F_{v,Y})$ the space
$(\CoY,\|\cdot|\CoY\|)$ is a Banach space which is continuously embedded in $(\mathcal{H}^1_{v})^{\sim}$
and depends on the frame $\mathcal{F}$. Moreover, we have the identities $\Co L_1^v = \mathcal{H}^1_v$,
$\Co L_{\infty}^{1/v} = (\mathcal{H}^1_{v})^{\sim}$, and
$\Co L_2 = \mathcal{H}$\,.

Suppose that $w$ is another weight function such that $(F_{w,Y})$ is satisfied.
Let $m_w(\x,\y)$ be the associated weight on $X \times X$. If $m_w(\x,\y) \leq Cm_v(\x,\y)$ then the spaces $\CoY(v)$ and $\CoY(w)$ coincide and their norms are equivalent.

Finally, we shall focus on the essential question of the coincidence of the two
spaces $\Co(\cf,Y)$ and $\Co(\cg,Y)$, where $\cf$ and $\cg$ are two different
continuous frames. One way to answer the above question is the following proposition which is essentially
taken from \cite{fora05}. Since we start with tight continuous frames the situation simplifies slightly here.

\begin{lemma}\label{ind1} Let $Y$ be a Banach function space on $X$ satisfying property $(Y)$ and let $v$ be a weight function.
The tight continuous frames $\cg = \{g_\x\}_{\x\in X}$ and $\cf = \{f_\x\}_{\x\in X}$ on $\mathcal{H}$
are supposed to satisfy $(F_{v,Y})$. Moreover, we assume that the Gramian kernel
\begin{equation}\label{eq-28}
    G(\cf, \cg)(\x,\y) := \langle f_\y,g_\x\rangle \quad,\quad \x,\y\in X\,,
\end{equation}
belongs to the algebra $\mathcal{B}_{Y,m}$.
Then it holds
$$
    \Co(\cf,Y) = \Co(\cg, Y)\,
$$
in the sense of equivalent norms.
\end{lemma}
We close this paragraph with a result concerning the independence of the coorbit space $\Co (\mathcal{F},Y)$
on the used reservoir $(\mathcal{H}^1_v)^{\sim}$. We state a version of Theorem 4.5.13 in \cite{ra05-6}.

\begin{lemma}\label{ind2} Let $Y$ be a Banach function space on $X$ satisfying $(Y)$ and let $v \geq 1$ be a weight function. The definition of $\Co (\cf,Y)$ is independent of the reservoir $(\mathcal{H}^1_v)^{\sim}$ in
the following sense: Assume that $S \subset \mathcal{H}^1_v$ is a non-trivial locally convex vector
space and
 $\cf \subset S$ be a tight continuous frame satisfying $(\mathcal{F}_{v,Y})$.  
 Assume further that the reproducing formula $V_{\cf}f = R_{\cf}(V_{\cf}f)$ extends to all $f\in S^{\sim}$ (the topological anti-dual
of $S$) then 
$$
    \Co(\mathcal{F},Y) = \{f\in S^{\sim}~:~V_{\cf}f \in Y\}\,.
$$
\end{lemma}
\bproof Let $f\in \mathcal{S}^{\sim}$ such that $V_{\cf}f \in Y$. %f(x) = f(\psi_{x}) \in Y$. 
Since the reproducing formula extends to 
$\mathcal{S}^{\sim}$ we have $V_{\cf}f = R_{\cf}(V_{\cf}f)$ and hence $V_{\cf}f \in R_{\cf}(Y) \subset L_{\infty}^{1/v}(X,\mu)$ which gives $f\in (\mathcal{H}^1_v)^{\sim}$
by definition of the latter space.
\eproof

\subsection{Discretizations}
\label{sectdiscr}

Next we come to a main feature of coorbit space theory, the discretization machinery.
%We close this subsection with a further definition, which is essential for the discretization problem.
It is based on the following definition, which is
a slight modification of Definition 6 in \cite{fora05} according to Remark 5 there.

\begin{definition}\label{PropD} A tight continuous frame $\cf = \{\varphi_\x\}_{\x\in X}$ is said 
to possess property $D[\delta,m,Y]$ for a fixed $\delta>0$ and a weight $m:X\times X \to \re$ 
if there exists a moderate admissible covering
$\mathcal{U} = \mathcal{U}^{\delta} = \{U_i\}_{i\in I}$ of $X$ such
that
\[
 \sup_{i\in I}\sup_{\x,\y \in U_i} m(\x,\y) \leq
C_{m,\mathcal{U}}\;,
\]
if the kernel $R_{\cf}$ belongs to $\mathcal{B}_{Y,m}$, and if
$\osc_{\mathcal{U}}(\x,\y)$ and $\osc_{\mathcal{U}}^*(\x,\y)$ 
satisfy
$$
    \|\osc_{\mathcal{U}}|\mathcal{B}_{Y,m}\| < \delta\quad\mbox{ and }\quad \|
    \osc_{\mathcal{U}}^*|\mathcal{B}_{Y,m}\|<\delta\,.
$$
Here we put 
\begin{equation}\nonumber
      \osc_{\mathcal{U}}(\x,\y) := \sup\limits_{z\in Q_y}|\langle \varphi_\x,
      \varphi_\y-\varphi_\z\rangle|\\
    = \sup\limits_{\z\in Q_\y}|R_{\cf}(\x,\y)-R_{\cf}(\x,\z)|\,,
\end{equation}
$\osc_{\mathcal{U}}^*(\x,\y) = \osc_{\mathcal{U}}(\y,\x)$
and $Q_\y = \bigcup\limits_{y\in U_i}U_i$.
\end{definition}
The following lemma states conditions on the frame $\cf$ and the space $Y$ which ensure that 
at least the test functions in $\mathcal{H}^1_v$ are contained in $\Co Y$.
\begin{lemma} Let $Y$ be a Banach function space satisfying $(Y)$. Let further $v\geq 1$ be a weight function 
with the associated weight $m = m_v$ satisfying
$\sup_{i\in I}\sup_{\x,\y \in U_i} m(\x,\y) \leq C$ and put $v_i = \sup_{\x\in U_i}v(\x)$\,. The frame $\cf$ is supposed to satisfy $(F_{v,Y})$ as well as $D[1,1,Y]$ with 
corresponding covering $\mathcal{U} = \{U_i\}_{i\in I}$. If $\|\chi_{U_i}|Y\| \lesssim v_i$ then 
it holds $\|\varphi_\x|\Co Y\| \lesssim v(\x)$ and 
\begin{equation}\label{embs}
    \mathcal{H}^1_v \hookrightarrow \Co Y \hookrightarrow (\mathcal{H}^1_v)^{\sim}\,.
\end{equation}
\end{lemma}
\bproof For all $i \in I$ and $\x \in U_i$ we have %We have for all $i$ such that $x\in U_i$ 
\begin{equation}\nonumber
 \begin{split}
    \|\varphi_\x|\Co Y\| &= \|\langle\varphi_\x,\varphi_\y\rangle |Y\| 
    \leq \Big\|\mu(U_i)^{-1}\int_{X}\sup\limits_{\z \in Q_\x}|\langle \varphi_\z,\varphi_\y\rangle|\chi_{U_i}(\x)\,d\mu(\x)|Y\Big\|\\
    &\lesssim \|\osc+R_{\cf}|Y \to Y\|\cdot \mu(U_i)^{-1}\|\chi_{U_i}|Y\| \lesssim v_i \lesssim v(\x)\,.
  \end{split}
\end{equation}
The second embedding in \eqref{embs} follows from $R_{\cf}(Y) \subset L_{\infty}^{1/v}$. By Theorem 1 in
\cite{fora05} an element $f \in \mathcal{H}^1_v$ can be written as a
sum $f = \sum_{i\in I}|c_i|\varphi_{\x_i}$, where $I$ is a countable
subset and 
$$
    \|f|\mathcal{H}^1_v\| \asymp \inf\sum\limits_{i\in I}|c_i|v(\x_i)\,,
$$
where the infimum is taken over all representations of $f$ in the
above form. So let us take one of these representations and estimate
by using the triangle inequality
$$
    \|f|\CoY\| \leq \sum\limits_{i\in I} |c_i|\cdot\|\varphi_{\x_i}|\CoY\|
    \leq \sum\limits_{i\in I} |c_i|v(\x_i)\,.
$$
This concludes the proof.\eproof

We return to the question of ensuring $R_{\cf}(Y) \hookrightarrow L_{\infty}^{1/v}(X)$. 
The following lemma states a sufficient condition.
\begin{lemma}\label{Linftyemb} Let $Y$ be a Banach function space satisfying $(Y)$ and $v \geq 1$ be a weight function with associated weight $m$ satisfying
$\sup_{i\in I}\sup_{\x,\y \in U_i} m(\x,\y) \leq C$ and put $v_i = \sup_{\x\in U_i}v(\x)$\,.
If $\mathcal{U} = \{U_i\}_{i\in I}$ is a moderate admissible covering of $X$ and 
$\|\chi_{U_i}|Y\| \gtrsim 1/v_i\,$ then we have the continuous embedding $Y^{\sharp} \hookrightarrow (L_{\infty}^{1/v})^{\sharp}$. 
If the frame $\cf$ satisfies in addition 
$D[1,1,Y]$ with respect to this covering then we even have
$R_{\cf}(Y) \subset L_{\infty}^{1/v}(X,\mu)$.
\end{lemma}

\bproof Using the assumption $\|\chi_{U_i}|Y\| \gtrsim 1/v_i$ we get by the solidity of $Y$
$$
    \Big\|\sum\limits_{i \in I}|\lambda_i|\mu(U_i)^{-1}\chi_{U_i}|Y\Big\| \geq \mu(U_i)^{-1}|\lambda_i|\cdot\|\chi_{U_i}|Y\|
    \gtrsim \mu(U_i)^{-1}|\lambda_i|\frac{1}{v_i} 
$$
for all $i\in I$. This yields
$$
   \Big\|\sum\limits_{i \in I}|\lambda_i|\mu(U_i)^{-1}\chi_{U_i}|Y\Big\| \gtrsim 
   \sup\limits_{i\in I}\mu(U_i)^{-1}|\lambda_i|\frac{1}{v_i} \gtrsim \|\{\lambda_i\}_{i\in I}|(L^{1/v}_{\infty})^{\sharp}\|\,,
$$
where we applied \eqref{def:N:overlap} in the last step. This proves $Y^{\sharp} \hookrightarrow (L_{\infty}^{1/v})^{\sharp}$. To show $R_{\cf}(Y) \subset L_{\infty}^{1/v}(X,\mu)$ we start with $F \in Y$ and estimate as follows,
\begin{align*}%\nonumber
%   \begin{split}
      &\|R_{\cf}(F)|L_{\infty}^{1/v}\| = \sup\limits_{\x\in X}|R_{\cf}(F)(\x)|\frac{1}{v(\x)} \lesssim \sup\limits_{i\in I}\sup\limits_{\x\in U_i} \int_{X}|R_{\cf}(\x,\y)F(\y)|d\mu(\y)\frac{1}{v_i}\\
      &\lesssim \sup\limits_{i\in I}\Big\|\sup\limits_{\z\in U_i} \int_{X}|R_{\cf}(\z,\y)F(\y)|d\mu(\y)\chi_{U_i}(\x)|Y\Big\|
      \lesssim \sup\limits_{i\in I}\Big\|\int_{X}\sup\limits_{\z\in Q_\x}|R_{\cf}(\z,\y)F(\y)|d\mu(\y)\chi_{U_i}(\x)|Y\Big\|\\
      &\leq \Big\|\int_{X}(\osc^{\ast}(\x,\y)+|R_{\cf}(\x,\y)|)\cdot |F(\y)|d\mu(\y)|Y\Big\|\,.
%   \end{split}
\end{align*}
Property $D[1,1,Y]$ gives in particular 
the boundedness of the considered integral operator and we obtain
$$
  \|R_{\cf}(F)|L_{\infty}^{1/v}\| \leq c\|F|Y\|
$$
which concludes the proof.
\eproof

The following abstract discretization results for 
coorbit spaces is a slight generalization of Theorem 5 in
\cite{fora05}, see also Remark 5 there. We omit the proof since the necessary modifications
are straightforward.

\begin{Theorem}\label{discr2} Let $Y$ be a Banach space of functions on $X$ satisfying $(Y)$ and let $v \geq 1$ be a weight function
with associated weight $m$. Assume that $\cf = \{\varphi_\x\}_{\x\in X}$ is a tight continuous frame
satisfying $(F_{v,Y})$ and $D[\delta,m,Y]$ for some $\delta>0$ with corresponding 
moderate admissible covering ${\mathcal U}^{\delta}$ chosen in a way such that 
$$
    \delta(\|R|\mathcal{B}_{Y,m}\|+\max\{C_{m,\mathcal{U}^{\delta}}\|R|\mathcal{B}_{Y,m}\|,\|R|\mathcal{B}_{Y,m}\|+\delta\})\leq 1\,,
$$
where $C_{m,\mathcal{U}^{\delta}}$ is the constant from Definition \ref{PropD}\,. Choose points $\x_i \in U_i$.
Then the discrete system $\cf_d:= \{\varphi_{\x_i}\}_{i\in I}$ is both an atomic decomposition of $\CoY$ with corresponding sequence space $Y^{\sharp}$ as well as a Banach frame with corresponding sequence space $Y^\flat$\,. This means
that there exists a dual frame $\{e_i\}_{i\in I}$ such that for all $f\in \Co Y$
\begin{description}
    \item(a)
    $$
        \|f|\CoY\| \asymp \|\{\langle f,\varphi_{\x_i}\rangle\}_{i\in I}|Y^{\flat}\|\quad\mbox{and}\quad
        \|f|\CoY\| \asymp \|\{\langle f,e_i \rangle\}_{i\in I}|Y^{\sharp}\|\,.
    $$
    \item(b) If $f\in \CoY$ then the series
    $$
        f = \sum\limits_{i\in I}\langle f,e_i\rangle \varphi_{\x_i} = \sum\limits_{i\in I}\langle f,\varphi_{\x_i}\rangle e_i
    $$
    converge unconditionally in the norm of $\CoY$ 
    if the finite sequences are dense in $Y^{\sharp}$ and with
    weak$^\ast$--convergence induced by $(\mathcal{H}^1_v)^{\sim}$, in general.
\end{description}
\end{Theorem}
%
%Aufpassen mit Y^{\flat} = Y^{\sharp} bei der Einbettung von \ell_1^v \in Y^{sharp}
%Siehe Lemma \ref{propsequ}

In the sequel we are interested in (wavelet) bases for the spaces $\CoY$. 
%The main tools will be versions of Theorem C and D in
%\cite{Gr88}, see also Theorems 5.5, 5.6, 5.7 in the preprint version
%of \cite{ra05-3}\,. Since we want to obtain in addition characterizations in terms of atomic decompositions
%we give a complete different proof here.
In many situations, such as in wavelet analysis, one often has an orthonormal basis, biorthogonal basis or discrete tight frame for the Hilbert space
at disposal, which arises from sampling a continuous frame. (Of course, such an orthonormal basis has to be derived from different
principles than available in the abstract situation of coorbit space theory.) 
Then the next main discretization result, Theorem \ref{wbases2} below, provides simple conditions, which ensure that
the basis expansion extends to coorbit spaces, and characterizes them by means of associated sequence spaces.
Our result generalizes one of Gr{\"o}chenig in classical coorbit space theory, see \cite{Gr88} and also 
Theorem 5.7 in the preprint version of \cite{ra05-3}.
From an abstract viewpoint, extensions of basis expansions seem very natural. However, in classical function space theory usually much efforts
are carried out in order to provide such wavelet basis characterization. In contrast, our discretization result provides a general approach, which requires
to check only a single condition in a concrete setup.

%The goal is to discretize 
%$\CoY = \Co (\cf,Y)$ by a different frame $\cg = \{\psi_\x\}_{\x\in X}$. Later we will
%assume $\cg$ to be a wavelet frame. 

Before giving the precise statement of our result, we have to introduce some notation and state some auxiliary lemmas.
Given a continuous frame $\cf$ defining the coorbit space $\Co(\cf,Y)$ we would like to discretize
by a different frame $\cg = \{\psi_\x\}_{\x\in X}$. Essentially this reduces to
conditions on the Gramian kernel $G(\cf, \cg)(\x,\y)$ introduced
above. If $\mathcal{U} = \{U_i\}_{i\in I}$ denotes a moderate
admissible covering of $X$ and $\x_i\in U_i, i\in I$, then we define
the kernel
\begin{equation}\label{eq-26}
    K(\x,\y) = \sup\limits_{\z\in Q_\x} |G(\cf,\cg)(\z,\y)| = \sup\limits_{\z\in Q_\x} |\langle \varphi_\y, \psi_\z\rangle|\,,
\end{equation}
where $Q_\x = \bigcup\limits_{i~:~\x\in U_i} U_i$\,. Observe that $K(\x,\y)$ depends on $\cf, \cg$ and 
the covering $\mathcal{U}$\,.

\begin{lemma}\label{D} Let $Y$, $v$, $\mathcal{U} = \{U_i\}_{i \in I}$ be as above and $x_i \in U_i$, $i\in I$.
Let further $\cf = \{\varphi_\x\}_{\x\in X}$ be a tight continuous frame satisfying $(F_{v,Y})$,
and $\CoY = \Co(\cf, Y)$. Assume that 
$\cg = \{\psi_\x\}_{\x\in X} \subset \mathcal{H}^1_v$ is a further continuous frame such that the kernel $K$ in \eqref{eq-26} belongs to $\mathcal{B}_{Y,m}$. 
Then there exists a constant $C>0$ independent of $f$ such that
$$
    \|\langle f,\psi_{\x_i}\rangle \}_{i\in I}|Y^{\flat}\| \leq C\|f|\CoY(\cf,Y)\|\quad,\quad f\in \CoY\,.
$$
\end{lemma}
\bproof Since $\mathcal{F}$ is a tight continuous frame with frame constants one, we have $V_{\cf}^{\ast}V_{\cf} = \mbox{Id}$, see Subsection \ref{confr}.
We conclude that
 \begin{align}\nonumber
   %\begin{split}
     (V_{\cg}f)(\x_i) &= (V_{\cg}V_{\cf}^{\ast}V_{\cf}f)(\x_i)
                     = V_{\cg}\Big(\int_{X}V_{\cf}f(\y)\varphi_\y\,d\mu(\y)\Big)(\x_i)\\
                    & = \int_{X}V_{\cf}f(\y)\langle \varphi_\y, \psi_{\x_i} \rangle\,d\mu(\y).\notag
 %  \end{split}
\end{align}
This implies the relation $|(V_{\cg} f)(\x_i)| \leq K(|V_{\cf}f|)(\x_i)$\,. We define the function
$$
    H(x) = \sum\limits_{i \in I(\x)}(V_{\cg} f)(\x_i) \chi_{U_i}(\x)\,,
$$
where $I(\x) = \{j\in I:\x\in U_j\}$, and observe that by \eqref{def:N:overlap}
\begin{equation}\nonumber
  \begin{split}
      |H(\x)| &\leq  \sum\limits_{i\in I(\x)} \chi_{U_i}(\x) \int_{X}|(V_{\cf} f)(\y)|\cdot K(\x_i,\y)\,d\mu(\y)\\
      & \leq N \int_{X} |(V_{\cf} f)(\y)|\cdot K(\x,\y) \,d\mu(\y)\,.
     % & = \int_{X} |(V_{\cf} f)(\y)|\cdot K(\x,\y) \,d\mu(\y)\,.
  \end{split}
\end{equation}
Hence, $|H| \leq K(|(V_{\cf} f)|)$ and together with $(Y)$ and our assumption on $K$ we get finally
$$
   \|\{V_{\cg}f(\x_i)\}_{i\in I}|Y^{\flat}\| = \|H|Y\| \leq N \|K|Y \to Y\|\cdot \|V_{\cf}f|Y\| \leq C\|f|\CoY\|\,.
$$
\eproof

We need a further technical lemma.
\begin{lemma}\label{C} Let $Y$, $\mathcal{U}$, $v$, $\{x_i\}_{i\in I}$ and $m$ as above, such that $Y^{\sharp} \hookrightarrow (L_{\infty}^{1/v})^{\sharp}$. 
Let $\cf = \{\varphi_\x\}_{\x\in X}$ be a tight frame satisfying $(F_{v,Y})$, put $\CoY = \Co(\cf, Y)$, and assume $\cg = \{\psi_\x\}_{\x\in X} \subset \mathcal{H}^1_v$ to be a continuous 
frame such that also $K^*$, see \eqref{eq-26}, belongs to $\mathcal{B}_{Y,m}$. If $\{\lambda_i\}_{i\in I} \in Y^{\sharp}$ then the sum
$$
      f = \sum\limits_{i\in I} \lambda_i \psi_{\x_i}
$$
converges unconditionally in the weak$^{\ast}$-topology of $(\mathcal{H}^1_v)^{\sim}$ to an element $f \in \Co Y$
and there exists a constant $c>0$ such that 
\begin{equation}\label{eq-25}
    \|f|\CoY\| \leq c\|\{\lambda_i\}_i|Y^{\sharp}\|\,.
\end{equation}
If the finite sequences are dense in $Y^{\sharp}$ we even have unconditional convergence in the norm of $\Co Y$.
\end{lemma}

\bproof {\em Step 1:} We prove that $\sum_{i\in I}|\lambda_i|\cdot|\langle \psi_{\x_i},\varphi_\x\rangle|$ converges pointwise for every $\x\in X$ and that its pointwise limit function belongs to $L_{\infty}^{1/v}$. This implies that the sequence of partial sums of every rearrangement of $\sum_{i\in I} \lambda_i \psi_{\x_i}$ is uniformly bounded in $(\mathcal{H}_v^1)^{\sim}$. Since by Theorem 1 in \cite{fora05} $\mbox{span}\{\varphi_\x:\x\in X\}$ is dense in 
$\mathcal{H}^1_v$ we conclude with an analogous argument as used in \cite[Lem.\ 4.5.8]{ra05-6} that $\sum_{i\in I} \lambda_i \langle \psi_{\x_i},\varphi \rangle$ converges unconditionally for every $\varphi\in \mathcal{H}^1_v$. This defines the weak$^{\ast}$-limit of the 
expansion of $\sum_{i\in I} \lambda_i \psi_{\x_i}$. To show the necessary pointwise convergence we estimate as follows,
\begin{equation}\nonumber
  \begin{split}
   & \frac{1}{v(\x)}\sum\limits_{i}|\lambda_i\langle \psi_{\x_i},\varphi_\x\rangle| \leq \sum\limits_{i\in I} \frac{|\lambda_i|}{v_i}\cdot \frac{v_i}{v(\x)}|\langle \psi_{\x_i},\varphi_\x\rangle|\\
    &\leq \Big(\sup\limits_{i} \mu(U_i)^{-1}\frac{|\lambda_i|}{v_i}\Big)
    \int_{U_i} m(\x,\y)\sup\limits_{\z\in Q_\y}|\langle 
    \psi_\z,\varphi_\x\rangle|\,d\mu(\y)
    \lesssim \|\{\lambda_i\}_i|Y^{\sharp}\|\cdot \|K^{\ast}|\mathcal{A}_m\|\,.
  \end{split}
\end{equation}
In the last step we used the assumption $Y^{\sharp} \subset (L_{\infty}^{1/v})^{\sharp}$.

\medskip
\noindent
{\em Step 2}. We already know that $\sum_{i\in I} \lambda_i \psi_{\x_i} =: f\in (H^1_v)^{\sim}$. We claim that 
$f\in \Co Y$. Indeed,
\begin{align*}%\nonumber
  % \begin{split}
      &\|f|\Co(\cf,Y)\| = \|V_{\cf}f|Y\| = \Big\|\sum\limits_{i\in I}\lambda_i \langle \psi_{\x_i},\varphi_\x\rangle|Y\Big\|
      \leq \Big\|\sum\limits_{i\in I} |\lambda_i| \sup\limits_{\z\in U_i}|\langle \psi_{\z},\varphi_\x\rangle||Y\Big\|\\
      & = \Big\|\sum\limits_{i\in I} |\lambda_i| \mu(U_i)^{-1}
      \int_{X}\sup\limits_{\z\in Q_\y}|\langle \psi_{\z}, \varphi_\x\rangle|\chi_{U_i}(\y)\,d \mu(\y)|Y\Big\|\\
     % \leq \Big\|\sum\limits_{i\in I} \lambda_i \mu(U_i)^{-1}
      %\int_{X}\sup\limits_{z\in U_i}\langle \varphi_x, \psi_{z}\rangle\chi_{U_i}(y)\,d \mu(y)|Y\Big\|\\
      &\leq \Big\| \int_{X} K(\y,\x)\Big(\sum\limits_{i \in I}|\lambda_i| \mu(U_i)^{-1}\chi_{U_i}(\y)\Big)\,d\mu(\y)|Y\Big\|\,.
   %\end{split}
\end{align*}
By our assumption on $K^*$ we obtain consequently
$$
    \|f|\CoY\| \leq \|K^{\ast}|Y\to Y\|\cdot \Big\|\sum\limits_{i \in I}|\lambda_i| \mu(U_i)^{-1}\chi_{U_i}|Y\Big\|\,,
$$
which reduces to (\ref{eq-25}) using the definition of $Y^{\sharp}$\,. This type of argument also implies the convergence in $\Co Y$ if the finite sequences are dense in $Y^{\sharp}$.\eproof

Let now $\cg_r = \{\psi^{r}_\x\}_{\x\in X}$ and $\tilde{\cg}_r = \{\tilde{\psi}^r_{\x}\}_{\x\in X}$, $r=1,...,n$,
be continuous frames with associated Gramian kernels
$K_r(\x,\y)$ and $\tilde{K}_r(\x,\y)$ defined by (\ref{eq-26}) for a moderate
admissible covering $\mathcal{U} = \{U_i\}_{i\in I}$\,. 

%The following remarkable result generalizes Theorem 5.7 in the preprint version of \cite{ra05-3} which goes back to Gr\"ochenig \cite{Gr88}. However, our proof technique is different and does not make use of the discretization result in Theorem \ref{discr2}. It therefore represents an independent discretization result. 
Now we are prepared to state our next discretization result. In contrast to the proof of its predecessor in classical coorbit theory \cite{Gr88}, we note,
however, that it does not rely on our first discretization result Theorem \ref{discr2}.

\begin{Theorem}\label{wbases2} Let $Y$ be as above and $v$ and $m$ such that $Y^{\sharp} \subset (L_{\infty}^{1/v})^{\sharp}$. Let $\cf = \{\varphi_\x\}_{\x\in X}$ be a tight frame satisfying $(F_{v,Y})$ and put $\CoY = \Co(\cf, Y)$. The continuous frames $\cg_r = \{\psi^r_\x\}_{\x\in X}, \tilde{\cg}_r = \{\tilde{\psi}_\x\}_{\x\in X} \subset \mathcal{H}^1_v$ are such that the corresponding
kernels $K_r$ and $\tilde{K}_r^{\ast}$ belong to $\mathcal{B}_{Y,m}$. Moreover, assume that
\begin{equation}\label{waveletexp2}
    f = \sum\limits_{r=1}^n \sum\limits_{i\in I} \langle
     f, \psi^r_{\x_i}\rangle \tilde{\psi}^{r}_{\x_i}
\end{equation}
holds for all $f\in \mathcal{H}$ where $\x_i \in U_i$ (the same covering which is used for the Gramian kernels $K_r$ and $\tilde{K}_r$)\,.
Then the expansion (\ref{waveletexp2}) extends to all
$f\in \CoY$. Furthermore, $f\in (\mathcal{H}_v^1)^{\sim}$ belongs to $\CoY$ if and only
if $\{\langle f, \psi^r_{\x_i}\rangle\}_{i\in I}$ belongs to
$Y^{\sharp}$ for each $r=1,...,n$\,. Then we have 
\begin{equation}\label{normequiv}
  \|f|\Co Y\| \asymp \sum\limits_{r=1}^n \|\{\langle
     f, \psi^r_{\x_i}\rangle\}_{i\in I}|Y^{\sharp}\|\,.
\end{equation}
The convergence in \eqref{waveletexp2} is
in the norm of $\CoY$ if the finite sequences are dense in $Y^{\sharp}$. In
general, we have weak$^{\ast}$-convergence induced by $(\mathcal{H}^1_v)^{\sim}$. 
\end{Theorem}

\bproof By Lemmas \ref{C} and \ref{D} the expansion 
\begin{equation}\label{eq-74}
    \sum\limits_{r=1}^n \sum\limits_{i\in I} \langle
     f, \psi^r_{\x_i}\rangle \tilde{\psi}^{r}_{\x_i}
\end{equation}
converges in the weak$^{\ast}$-topology of $(\mathcal{H}^1_v)^{\sim}$ to an element $\tilde{f} \in \Co Y \subset (H_v^1)^{\sim}$ provided we assume 
that either $f\in \Co Y$ or $\{\langle f, \psi^r_{\x_i}\rangle\}_{i\in I}$ belongs to
$Y^{\sharp}$ for each $r=1,...,n$. If the finite sequences are dense in $Y^{\sharp}$ we even have convergence in $\Co Y$. 
It remains to show the identity $f = \tilde{f}$. % in this situation. 

{\em Step 1.} Let us start with an $\varphi\in \mathcal{H}^1_v$. We apply Lemma \ref{D} to the case $Y = L_1^v$ and $\mathcal{G} = \tilde{\mathcal{G}}_r$, $r=1,...,n$. The assumption $\tilde{K}_r \in \mathcal{A}_m$ implies then that $\tilde{K}_r$ maps $L_1^v$ boundedly into $L_1^v$. Therefore, 
Lemma \ref{D} yields that $\{\langle \varphi, \tilde{\psi}^r_{\x_i}\rangle\}_{i\in I}$ belongs to $\ell_1^{v_i}$ for all $r=1,...,n$. Lemma \ref{C}
gives then that the expansion 
\begin{equation}\label{eq-73}
      \sum\limits_{r=1}^n \sum\limits_{i\in I} \langle \varphi, \tilde{\psi}^r_{\x_i}\rangle \psi^r_{\x_i}
\end{equation}
converges in the norm of $\mathcal{H}^1_v$ to an element $g \in \mathcal{H}^1_v$ since the finites sequences are dense in $\ell_1^{v_i}$. 
Observe that our global assumption $v>1$ together with
$$
    \|h|\mathcal{H}\|^2 \leq \|\langle h,\varphi_\x\rangle|L_2(X)\|^2 \leq \|\langle h,\varphi_\x\rangle|L_{\infty}(X)\|\cdot
    \|\langle h,\varphi_\x\rangle|L_1(X)\| \leq \|\langle h,\varphi_\x\rangle|L_{\infty}(X)\|\cdot \|h|\mathcal{H}^1_v\|
$$
and 
$$
    \|\langle h,\varphi_\x\rangle|L_{\infty}(X)\| \leq \|h|\mathcal{H}\|\cdot \|\varphi_x|\mathcal{H}\| \lesssim \|h|\mathcal{H}\|\,,
$$
using $\|\varphi_\x|\mathcal{H}\| \leq C$, imply the continuous embedding $\mathcal{H}_1^v \hookrightarrow \mathcal{H}$.
Hence, \eqref{eq-73} converges also in $\mathcal{H}$ to $g$\,. On the other hand the identity in $\mathcal{H}$
$$    
    \eta = \sum\limits_{r=1}^n \sum\limits_{i\in I} \langle
     \eta, \psi^r_{\x_i}\rangle \tilde{\psi}^{r}_{\x_i}
$$
for arbitrary $\eta\in \mathcal{H}$ gives 
$$
     \langle \eta,\varphi\rangle  = \sum\limits_{r=1}^n \sum\limits_{i\in I} \langle
     \eta, \psi^r_{\x_i}\rangle \langle \tilde{\psi}^{r}_{\x_i},\varphi\rangle = 
     \langle \eta, \sum\limits_{r=1}^n \sum\limits_{i\in I}\langle \varphi, \tilde{\psi}^{r}_{\x_i}\rangle\psi^r_{\x_i}\rangle = \langle \eta,g\rangle\,.
$$
Choosing $\eta = \varphi-g$ gives $\varphi = g$\,.

{\em Step 2.} Using that \eqref{eq-73} converges to $\varphi$ in $\mathcal{H}^1_v$  and that $\tilde{f}$ is the weak$^\ast$--limit
of \eqref{eq-74}, we finally obtain
\begin{equation}\nonumber
  \begin{split}
    f(\varphi) &= f\Big(\sum\limits_{r=1}^n \sum\limits_{i\in I} \langle \varphi, \tilde{\psi}^r_{\x_i}\rangle \psi^r_{\x_i}\Big) = \sum\limits_{r=1}^n \sum\limits_{i\in I} \langle \tilde{\psi}^r_{\x_i},\varphi \rangle f(\psi^r_{\x_i})\\
    &= \sum\limits_{r=1}^n \sum\limits_{i\in I} \langle f, \psi^r_{\x_i}\rangle \langle \tilde{\psi}^r_{\x_i},\varphi \rangle 
     = \tilde{f}(\varphi)\,.
  \end{split}
\end{equation}
This implies $f = \tilde{f}$ since $\varphi$ was chosen arbitrarily. 
The norm equivalence in \eqref{normequiv} is a direct consequence of Lemmas \ref{D}, \ref{C}.\eproof

\section{Peetre type spaces and their coorbits}
\label{Peetre}

The generalized Besov-Lizorkin-Triebel spaces to be studied later in Section \ref{examples} are defined 
as coorbits of so-called Peetre spaces on the index set
$X = \R \times [(0,1) \cup \{\infty\}]$ equipped with the 
Radon measure $\mu$ given by
$$
    \int_{X} F(x) d\mu(x) = \int_{\R}\int_{0}^1 F(y,s) \frac{ds}{s^{d+1}}dy + \int_{\R} F(y,\infty) dy\,.
$$
We intend to define two general scales of Banach function spaces $P^w_{B,q,a}(X)$ and $L^w_{B,q,a}(X)$ on $X$. The parameter $B(\R)$ is a Banach 
space of measurable functions on $\R$, the parameter $w:X \to (0,\infty)$ represents a weight function on $X$, and 
$1\leq q\leq \infty$, $a>0$. The letter $P$
refers to Peetre's maximal function \eqref{def:Peetre:max} which is always involved in the definition of $P^w_{B,q,a}(X)$, see Definition \ref{defFS} below.
Let us start with reasonable restrictions on the parameters $w$ and $B(\R)$. We use the class $\mathcal{W}^{\alpha_3}_{\alpha_1,\alpha_2}$ of admissible
weights introduced by Kempka \cite{Ke09}\,.

\begin{definition}\label{def:weight} A weight function $w:X \to \re_+$ belongs 
to the class $\mathcal{W}^{\alpha_3}_{\alpha_1,\alpha_2}$
if and only if
there exist non-negative numbers $\alpha_1,\alpha_2,\alpha_3 \geq 0$ such that, 
for $\x=(x,t) \in X$, 
\begin{description}
 \item(W1) %for $t\in (0,1)$ and $x\in \R$
 $
    \left\{\begin{array}{lcl}
              \Big(\frac{t}{s}\Big)^{\alpha_1}w(x,s) \leq w(x,t) \leq \Big(\frac{t}{s}\Big)^{-\alpha_2}w(x,s)&:& 
              1 \geq s \geq t > 0\\\\
              t^{\alpha_1}w(x,\infty) \leq w(x,t) \leq t^{-\alpha_2}w(x,\infty)&:& s = \infty\,,0<t\leq 1,
          \end{array}\right.
 $
 \item(W2) 
 $
    w(x,t) \leq w(y,t)\left\{\begin{array}{rcl}
                                (1 + |x-y|/t)^{\alpha_3}&:& t\in (0,1)\\
                                (1+|x-y|)^{\alpha_3}&:& t=\infty
                             \end{array}\right.\, \quad\mbox{ for all }y \in \R.
 $
\end{description}
\end{definition}

\begin{example}  The main examples are weights of the form
$$
    w_{s,s'}(x,t) = \left\{\begin{array}{rcl}
                                t^{-s}\Big(1+\frac{|x-x_0|}{t}\Big)^{s'} &:& t\in (0,1)\\
                                (1+|x-x_0|)^{s'}&:& t=\infty
                    \end{array}\right.\,.
$$
where $s,s' \in \re$. The choice $s' = 0$ is most common.
\end{example}
\begin{remark} The above considered weights are continuous 
versions of weights appearing in the definition 
of certain $2$-microlocal function spaces of Besov-Lizorkin-Triebel type, see for instance \cite{Ke09,Ke10,Ke11}.
\end{remark}
The second ingredient is a Banach space $B(\R)$ of functions defined on $\R$. 
%Let us consider the following class
%of admissible Banach spaces.

\begin{definition}\label{condB} A solid Banach space $B = B(\R)$ of functions on $\R$ with norm $\|\cdot|B(\R)\|$ is called admissible if 
\begin{description} 
 \item(B1) the elements of $B(\R)$ are locally integrable functions with respect to the Lebesgue measure;
 \item(B2) there exist real numbers $\gamma_1 \leq \gamma_2$ and $\delta_2 \leq \delta_1$ with $\delta_1\geq 0$  such that for every $\alpha>0$
 there are constants $C_{\alpha},c_{\alpha}$ with
 $$
    c_{\alpha}t^{\gamma_2}\Big(1+\frac{|x|}{t}\Big)^{\delta_2}\leq\|\chi_{Q^{\alpha}_{(x,t)}}|B(\R)\| \leq C_{\alpha}t^{\gamma_1}\Big(1+\frac{|x|}{t}\Big)^{\delta_1}\quad,\quad x\in \R, t\in (0,1]\,,
 $$  
 where $Q^{\alpha}_{(x,t)} = x + t [-\alpha,\alpha]^d$ denotes a $d$-dimensional cube with center $x\in \R$.
 %$$
  %  Q^{\alpha}_{(x,t)}:= (x_1+t[-\alpha,\alpha])\times\cdots \times (x_d+t[-\alpha,\alpha])\,.
 %$$
\end{description}
\end{definition}

\begin{example}\label{Lp} If $B(\re^d) = L_p(\re^d)$ is the classical Lebesgue space then
  $$
      \|\chi_{Q^{\alpha}_{(x,t)}}|B(\re^d)\| = (2\alpha)^{d}t^{d/p}\,.
  $$
  Hence, the parameters in condition (B2) are given by 
  $C_{\alpha} = c_{\alpha} = (2\alpha)^d$, $\gamma_1=\gamma_2 = d/p$, and $\delta_1 = \delta_2 = 0$.
\end{example}

\subsection{Peetre type spaces on $X$}
\label{PeetSp}

Our key ingredient in recovering generalized Besov-Lizorkin-Triebel spaces are the following function spaces
on $X$ defined via the Peetre maximal function in \eqref{def:Peetre:max}.

\begin{definition}\label{defFS} Let $1\leq q\leq \infty$, $a>0$, and $w \in \mathcal{W}^{\alpha_3}_{\alpha_1,\alpha_2}$ 
be a weight function. Assume that $B(\R)$ is a solid Banach space of functions on $\R$ satisfying $(B1)$ and $(B2)$. Then we define by
\begin{equation}\nonumber
  \begin{split}
     P^{w}_{B,q,a}(X) &:= \{F:X \to \C~:~\|F|P^{w}_{B,q,a}\| < \infty\}\,,\\
     L^{w}_{B,q,a}(X) &:= \{F:X \to \C~:~\|F|L^{w}_{B,q,a}\| < \infty\}\\
  \end{split}   
\end{equation}
two scales of function spaces on $X$, where the norms are given by \eqref{def:Peetre:Pspace} and \eqref{def:Peetre:Lspace}.
%
%put
%\begin{equation}
%  \begin{split} 
%   \|F|P^{w}_{B,q,a}\| := &\Big\|w(\cdot,\infty)\sup\limits_{z\in \R}\frac{|F(\cdot+z,\infty)|}{(1+|z|)^a}|B(\R)\Big\|\\
%    &+ \Big\|\Big(\int_{0}^1 \Big[w(\cdot,t)\sup\limits_{z\in \R}
%   \frac{|F(\cdot+z,t)|}{(1+|z|/t)^a}\Big]^q\frac{dt}{t^{d+1}}\Big)^{1/q}|B(\R)\Big\|
%  \end{split} 
%\end{equation}
%and 
%\begin{equation}\label{eq-10}
%  \begin{split} 
%   \|F|L^{w}_{B,q,a}\| := &\Big\|w(\cdot,\infty)\sup\limits_{z\in \R}\frac{|F(\cdot+z,\infty)|}{(1+|z|)^a}|B(\R)\Big\|\\
%    &+ \Big(\int_{0}^1 \Big\|w(\cdot,t)\sup\limits_{z\in \R}
%   \frac{|F(\cdot+z,t)|}{(1+|z|/t)^a}|B(\R)\Big\|^q\frac{dt}{t^{d+1}}\Big)^{1/q}\,.
%  \end{split} 
%\end{equation}
\end{definition}

\begin{remark}\label{Tx} Assume that in addition the space $B(\R)$ is %uniform 
uniformly translation invariant, i.e., the translation operators
defined by $T_x g = g(\cdot -x)$ are uniformly bounded from $B(\R)$ to $B(\R)$, % In other words, we assume
\begin{equation}\label{eq31}
    \sup\limits_{x\in \R} \|T_x:B(\R)\to B(\R)\| < \infty\,.
\end{equation}
Moreover, we assume that 
$$
    w(x,t) = \tilde{w}(t)\quad,\quad (x,t)\in X\,.
$$
Under this stronger condition we define the scale of spaces
$L^w_{B,q}(X)$, $1\leq q\leq \infty$, by
\begin{equation}\nonumber
   \|F|L^{w}_{B,q}\| := |\tilde{w}(\infty)|\cdot\|F(\cdot,\infty)|B(\R)\|
    + \Big(\int_{0}^1 |\tilde{w}(t)|^q\|F(\cdot,t)|B(\R)\|^q\frac{dt}{t^{d+1}}\Big)^{1/q}.
\end{equation}
These spaces can then also be taken in replacement of $L^w_{B,q,a}(X)$.
An important class of examples of uniformly translation invariant spaces are the unweighted classical Lebesgue space $L_p(\R)$. 
\end{remark}
In the following we prove assertions on the boundedness of certain integral operators between 
these spaces. Recall that for a function $G:X \to \com$ the action of a kernel
$K$ on $G$ is defined by
$$
   K(G)(x,t) = \int_{\R} K((x,t),(y,\infty))G(y,\infty)\,dy +
   \int_{\R} \int_{0}^1K((x,t),(y,s))G(y,s)\,\frac{ds}{s^{d+1}}\,dy\,.
$$
Condition \eqref{Kb} below will be satisfied for kernels associated to continuous wavelet transforms to be studied later.

\begin{proposition}\label{bkernel} Assume that $K((x,t),(y,s))$ denotes a kernel function on $X \times X$
such that
\begin{equation}\label{Kb}
 K((x,t),(y,s)) \leq \left\{\begin{array}{rcl}
        G_1\Big(\frac{y-x}{t},\frac{s}{t}\Big)&:&t,s \in (0,1)\;,\\
        G_2\Big(\frac{y-x}{t},\frac{1}{t})&:& t\in (0,1), s=\infty\;,\\
        G_3(y-x,s)&:& t=\infty, s\in (0,1)\;,\\
        G_4(y-x)&:&t=s=\infty\;
 \end{array}\right.
\end{equation}
for some functions $G_1,G_2,G_3,G_4$.
Let $1 \leq q \leq \infty$, $a>0$, and $w \in \mathcal{W}^{\alpha_3}_{\alpha_1,\alpha_2}$. 
Assume $B(\R)$ is a solid Banach function space satisfying (B1) and (B2) and 
suppose that the following quantities are finite, 
\begin{equation}\label{bddcond}
\begin{split}
     M_1&:=\int_{0}^{\infty} \int_{\R} |G_1(y,r)|(1+|y|)^{a}
     r^{d/q}\max\{1,r^{-a}\}\max\{r^{-\alpha_1},r^{\alpha_2}\}dy\,\frac{dr}{r^{d+1}}\,,\\
     M_2&:= \int_{1}^{\infty} \int_{\R} t^{\alpha_2+d/q}(1+|y|)^{a}
     \sup\limits_{t/2\leq t'\leq t}|G_2(y,t')|\,dy \frac{dt}{t^{d+1}}\,,\\
     M_3&:= \int_{0}^1 \int_{\R}r^{-(\alpha_1+2a+d/q'-d)}\sup\limits_{r/2\leq r' \leq r}|G_3(y,r')|(1+|y|)^a
    \,dy\frac{dr}{r^{d+1}}\,,\\
     M_4&:=\int_{\R}|G_4(y)|(1+|y|)^a\,dy\,,
\end{split}
\end{equation}
where $q'$ is such that $1/q + 1/q'=1$.
%are supposed to be finite. 
Then %we have 
\begin{align}
   \|K:P^w_{B,q,a}(X) \to P^{w}_{B,q,a}(X)\|
    &\lesssim M_1+M_2+M_3+M_4\notag\\
    \mbox{and} \qquad  
    \|K:L^{w}_{B,q,a}(X) \to L^{w}_{B,q,a}(X)\|
    &\lesssim M_1+M_2+M_3+M_4\,.\notag
\end{align}
\end{proposition}
\bproof We prove the assertion only for the space $P^w_{B,q,a}$. For $L^w_{B,q,a}$ the calculation is simpler and the modifications are straightforward. 
We first observe that, for a function $F \in P^w_{B,q,a}(X)$,
\begin{equation}\nonumber
   \begin{split}
        &\|K(F)|P^w_{B,q,a}\| \leq \\
        &~\Big\|\Big(\int_{0}^1 \Big[\sup\limits_{z\in \R}
        \frac{w(\cdot,t)}{(1+|z|/t)^a}\int_{\R}
        \int_{0}^1 |K((\cdot+z,t),(y,r))F(y,r)|\frac{dr}{r^{d+1}}\,dy\Big]^q\frac{dt}{t^{d+1}}\Big)^{1/q}|B(\R)\Big\|\\
        &+\Big\|\Big(\int_{0}^1 \Big[\sup\limits_{z\in \R}
        \Big[\frac{w(\cdot,t)}{(1+|z|/t)^a}\int_{\R}
        K((\cdot+z,t),(y,\infty))F(y,\infty)\,dy\Big]^q\frac{dt}{t^{d+1}}\Big)^{1/q}|B(\R)\Big\|\\
        &+\Big\|\sup\limits_{z\in \R}
        \frac{w(\cdot,\infty)}{(1+|z|)^a}\int_{0}^1\int_{\R}
        |K((\cdot+z,\infty),(y,r))F(y,r)|\,dy\frac{dr}{r^{d+1}}|B(\R)\Big\|\\
        &+\Big\|\sup\limits_{z\in \R}
        \frac{w(\cdot,\infty)}{(1+|z|)^a}\int_{\R}
        |K((\cdot+z,\infty),(y,\infty))F(y,\infty)|\,dy|B(\R)\Big\|\,.
   \end{split}
\end{equation}
We denote the summands appearing on the right-hand side by
$S_1,S_2,S_3,S_4$. Let us first treat $S_4$. We have
\begin{align*}
   S_4 &\leq \Big\|\sup\limits_{z\in \R}
   \frac{w(\cdot,\infty)}{(1+|z|)^a}\int_{\R}
   |G_4(y-(\cdot+z))F(y,\infty)|\,dy|B(\R)\Big\|\\
   &\leq \Big\|\int_{\R}
   |G_4(y)|w(\cdot,\infty)\sup\limits_{z\in \R}
   \frac{|F(\cdot + y+z,\infty)|}{(1+|z|)^a}\,dy|B(\R)\Big\|\\
   & \leq \int_{\R}
   |G_4(y)|(1+|y|)^a\,dy \Big\|w(\cdot,\infty)\sup\limits_{z\in \Z}
   \frac{|F(\cdot+z,\infty)|}{(1+|z|)^a}|B(\R)\Big\| = M_4\|F|P^{w}_{B,a,q}\|\,.
\end{align*}
%The last term can be estimated by
%$$
%   \Big\|w(x,\infty)\sup\limits_{z\in \Z}
%   \frac{|F(x+z,\infty)|}{(1+|z|)^a}|B(\R)\Big\|\,,
%$$
%which gives
%$$
%   S_4 \lesssim M_4\|F|P^{w,a,q}_B\|\,,
%$$
%where
%$$
%  M_4 = \int_{\R}
%   |G_4(y)|(1+|y|)^a\,dy
%$$
%%
Similarly, we obtain
\begin{align*}%\nonumber
%\begin{split}
    S_2 &\leq \Big\|\Big(\int_{0}^1 \Big[
   \sup\limits_{z\in \R}
   \frac{w(\cdot,t)}{(1+|z|/t)^{a}}%\times\\
  % &~~~~~~~~~~\times
  \int_{\R}
   \Big|G_2\Big(\frac{y-(\cdot+z)}{t},\frac{1}{t}\Big)\cdot F(y,\infty)\Big|\,dy\Big]^q\frac{dt}{t^{d+1}}\Big)^{1/q}|B(\R)\Big\|\\
   &\lesssim\Big\|w(\cdot,\infty)\sup\limits_{z\in \R}
   \frac{|F(\cdot+z,\infty)|}{(1+|z|)^a}|B(\R)\Big\|\Big(\int_{0}^1 \Big[t^{-\alpha_2}
   \int_{\R}\Big|G_2\Big(\frac{y}{t},\frac{1}{t}\Big)(1+|y|)^a\Big|\,dy\Big]^q\frac{dt}{t^{d+1}}\Big)^{1/q}\\
   &\lesssim \Big(\int_{1}^{\infty} \Big[t^{\alpha_2-d+2d/q}
   \int_{\R}|G_2(y,t)|(1+|y|)^a\,dy\Big]^q\frac{dt}{t^{d+1}}\Big)^{1/q} \|F|P^w_{B,q,a}\| \leq M_2 \|F|P^{w}_{B,q,a}\| \;.
%\end{split}
\end{align*}
%
%where we put
%$$
%   M'_2 = \Big(\int_{1}^{\infty} \Big[t^{\alpha_2-d+2d/q}
%   \int_{\R}|G_2(y,t)|(1+|y|)^a\,dy\Big]^q\frac{dt}{t^{d+1}}\Big)^{1/q}
%$$
%which satisfies $M'_2 \leq M_2$\,.
%%
The next step is to estimate
\begin{align*}%\nonumber
%  \begin{split}
    S_1&\lesssim \Big\|\Big(\int_{0}^{1} \Big[
    w(\cdot,t)\int_{0}^{1/t}\int_{\R}
    |G_1(y,r)|\sup\limits_{z\in \R}\frac{|F(\cdot+z+ty,rt)|}{(1+|z|/t)^a}\,dy\frac{dr}{r^{d+1}}\Big]^q\frac{dt}{t^{d+1}}
    \Big)^{1/q}|B(\R)\Big\|\\
    &\lesssim \Big\|\Big(\int_{0}^{1} \Big[w(\cdot,t)
    \int_{0}^{1/t}\int_{\R}
    |G_1(y,r)|(1+|y|)^a
    %\times\\
    %&~~~~~~\times
    \sup\limits_{z\in \R}\frac{|F(\cdot+z,rt)|}{(1+|z|/t)^a}\,dy\frac{dr}{r^{d+1}}\Big]^q\frac{dt}{t^{d+1}}
    \Big)^{1/q}|B(\R)\Big\|\,.
  %\end{split}
\end{align*}
Minkowski's inequality and a change of variable in the integral over $t$ gives
\begin{equation}\nonumber
  \begin{split}
     S_1 &\lesssim \int_{0}^{\infty}\int_{\R}|G_1(y,r)|(1+|y|)^a\max\{1,r^{-a}\}\max\{r^{-\alpha_1},r^{\alpha_2}\}
     r^{d/q}\times\\
     &~~~~~~\times\Big\|\Big(\int_{0}^{1}
     \Big[w(\cdot,t)\sup\limits_{z\in \R}\frac{|F(\cdot+z,t)|}{(1+|z|/t)^a}\Big]^q\,\frac{dt}{t^{d+1}}\Big)^{1/q}
     |B(\R)\Big\|\,dy\frac{dr}{r^{d+1}}\\
     &\lesssim M_1 \|F|P^{w}_{B,q,a}\|\,.
  \end{split}
\end{equation}
%where we put 
%$$
%    M_1=\int_{0}^{\infty}\int_{\R}|G_1(y,r)|(1+|y|)^a\max\{1,r^{-a}\}
%    \max\{r^{-\alpha_1},r^{\alpha_2}\}r^{d/q}\,dy\frac{dr}{r^{d+1}}\,.
%$$
It remains to estimate $S_3$. Using (W1) we get
\begin{equation}\label{eq-6}
  \begin{split}
    S_3 &\leq \Big\|\sup\limits_{z\in \R}
    \frac{w(\cdot,\infty)}{(1+|z|)^a}\int_{0}^1\int_{\R}
    |G_3(y-(\cdot+z),r)F(y,r)|\,dy\frac{dr}{r^{d+1}}|B(\R)\Big\|\\
    &\leq \Big\|\int_{0}^1\int_{\R}
    w(\cdot,r)r^{-\alpha_1}|G_3(y,r)|\sup\limits_{z\in \R}
    \frac{|F(\cdot + y+z,r)|}{(1+|z|)^a}\,dy\frac{dr}{r^{d+1}}|B(\R)\Big\|\,.
\end{split}
\end{equation}
For $r \in (0,1)$ we can estimate the supremum above by %taking $r\in (0,1)$ into account by
\begin{align*}%\nonumber
 %\begin{split}
    &\sup\limits_{z\in \R}
    \frac{|F(x+y+z,r)|}{(1+|z|)^a} = \sup\limits_{w\in \R}
    \frac{|F(x+w,r)|}{(1+|w-y|)^a}\cdot\frac{(1+|w|/r)^a}{(1+|w|/r)^a}\\
    &\leq \sup\limits_{w\in \R}
    \frac{|F(x+w,r)|}{(1+|w|/r)^a} \cdot r^{-a}(1+|y|/r)^a
    \leq \sup\limits_{w\in \R}
    \frac{|F(x+w,r)|}{(1+|w|/r)^a}r^{-2a}(1+|y|)^a\,.
  %\end{split}
\end{align*}
Pluggin this into (\ref{eq-6}) and using H{\"o}lder's inequality with $1/q+1/q' = 1$ we finally get
$$
  S_3 \lesssim \Big(\int_{0}^1\Big[r^{-(\alpha_1+2a)}
  \int_{\R}G_3(y,r)(1+|y|)^a\,dy\Big]^{q'}\,\frac{dr}{r^{d+1}}\Big)^{1/q'}
 \|F|P^{w}_{B,a,q}\| \leq M_3 \|F|P^{w}_{B,a,q}\|\;.
$$
This concludes the proof.
%where 
%$$
%    M'_3 = $$
%and $M'_3 \leq M_3$\,.
\eproof
\begin{remark} According to Remark \ref{Tx} the conditions in \eqref{bddcond} are simpler in the translation invariant case. 
The parameter $a$ is then not required. 
\end{remark}

We need a similar statement in order to guarantee that $K$ belongs to $\mathcal{A}_{m_v}$, where $m_v$ is the associated
weight to $v:X \to \re$ given by
\begin{equation}\label{mw}
    m_v((x,t),(y,s)) := \max\Big\{\frac{v(x,t)}{v(y,s)}, \frac{v(y,s)}{v(x,t)} \Big\}
\end{equation}
for the special choice 
\begin{equation}\label{weight}
    v(x,t) := \left\{\begin{array}{rcl}
                        t^{-\gamma}(1+|x|/t)^\eta&:&t\in (0,1]\\
                        (1+|x|)^{\eta}&:& t=\infty
                    \end{array}\right.,
\end{equation}
where $\eta,\gamma\geq 0$.
Recall that we define $K^{\ast}(\x,\y) = K(\y,\x)$.

\begin{proposition}\label{bddalgebra} Let $K$ be a kernel function on $X \times X$ such that $K$ and $K^{\ast}$ satisfy \eqref{Kb}
with functions $G_i$ and $G_i^{\ast}$, $i=1,...,4$, respectively. Let further $v$ and $m_v$ be given by \eqref{weight}
and \eqref{mw}. If the quantities
\begin{equation}\label{Am}
\begin{split}
     S_1&:=\int_{0}^{\infty} \int_{\R} |G_1(y,t)|\max\{t,t^{-1}\}^{|\eta|+|\gamma|}(1+|y|)^{|\eta|}dy\,\frac{dt}{t^{d+1}}\,,\\
     S_2&:= \esssup{t>1}\,t^{|\eta|+|\gamma|-d}\int_{\R} |G_2(y,t)|(1+|y|)^{|\eta|}dy\,,\\
     S_3&:= \int_{0}^1 \int_{\R}|G_3(y,t)|t^{-(|\eta|+|\gamma|)}(1+|y|)^{|\eta|}dy\,\frac{dt}{t^{d+1}}\,,\\
     S_4&:=\int_{\R}|G_4(y)|(1+|y|)^{|\eta|}\,dy\,,
\end{split}
\end{equation}
and the corresponding ones for $K^{\ast}$ in terms of the function $G_i^{\ast}$ are finite then we have
$K,K^{\ast} \in \mathcal{A}_{m_v}$.

\end{proposition}

\bproof A straightforward computation shows  
\begin{equation}\label{eq59}
    m_v((x,t),(y,s)) \leq \max\Big\{\frac{s}{t},\frac{t}{s}\Big\}^{|\gamma|+|\eta|}\Big(1+\frac{|x-y|}{t}\Big)^{|\eta|}
\end{equation}
(obvious modification in case $s=\infty$ or $t=\infty$). According to \eqref{def:A:norm1} we have to show that 
$$
    \|K|\mathcal{A}_{m_v}\| := \max\Big\{\esssup{\x\in X}\int_{Y}m_v(\x,\y)|K(\x,\y)|d\mu(\y)~,~ \esssup{\y\in X} \int_{X}m_v(\x,\y)|K(\x,\y)|d\mu(\x)\Big\}\,
$$
(similar for $K^{\ast}$) is finite. Combining \eqref{Kb}, \eqref{Am}, and \eqref{eq59} finishes the proof. 
\eproof

\subsection{Associated sequence spaces}
\label{seq}

As the next step we study the structure of the sequence spaces $P^{w}_{B,q,a}(X)^{\sharp}$
and $L^{w}_{B,q,a}(X)^{\sharp}$\, associated to Peetre type spaces. % from the preceding subsection.
We will use the following covering of the space $X$. 
For $\alpha>0$ and $\beta>1$ we consider
the family $\mathcal{U}^{\alpha,\beta} = \{U_{j,k}\}_{j\in \n, k\in \zz^d}$ of subsets
\begin{align}\nonumber
    U_{0,k} &= Q_{0,k}\times\{\infty\}
    \quad,\quad k\in \Z\,,\notag\\
U_{j,k} &= Q_{j,k} \times [\beta^{-j},\beta^{-j+1})\quad , \quad j \in \N, k \in \zz^d\,,\nonumber
\end{align}
where $Q_{j,k} = \alpha k + \alpha \beta^{-j}[0,1]^d$.
%$$
%    U_{j+1,k} = Q_{j+1,k}\times [\beta^{-(j+1)},\beta^{-j}]\quad,\quad j\in \n\,,k\in \Z\,,
%$$
%where
%$$
%   Q_{j,k} = [\alpha k_1\beta^{-j},\alpha (k_1+1)\beta^{-j}]\times \cdots \times[\alpha k_d\beta^{-j},\alpha (k_d+1)\beta^{-j}]\quad,\quad j\in \n.
%$$
We will use the notation
$$
    \chi_{j,k}(x) = \left\{\begin{array}{rcl}
                           1 &:& x\in U_{j,k}\;,\\
                           0 &:& \mbox{otherwise}\;.
                        \end{array}\right.
$$
Clearly, we have $X \subset \bigcup_{j\in \n, k\in \Z} U_{j,k}$ and $\mathcal{U}^{\alpha,\beta}$ is a moderate 
admissible covering of $X$.
We now investigate properties of the sequence spaces $(\CoP)^{\sharp}$ and $(\CoL)^{\sharp}$, recall 
Definition \ref{absdefss}. 

\begin{lemma}\label{chainembappl} Let $1 \leq q \leq \infty$, $a>0$, and
$w \in \mathcal{W}^{\alpha_3}_{\alpha_1,\alpha_2}$. Let  $B(\R)$ be a solid Banach space satisfying (B1) and (B2). 
Then we have 
\begin{equation}\label{eq-33}
        v_{j,k}^{-1} \lesssim \|\chi_{U_{j,k}}|P^{w}_{B,q,a}\|\quad,\quad j\in \n, k\in \Z\,,
\end{equation}
where $v_{j,k} = \sup\limits_{(x,t) \in U_{j,k}} v_{w,B,q}(x,t)$  with 
\begin{equation}\label{rem1}
v_{w,B,q}(x,t) = \left\{\begin{array}{rcl}
            t^{-|\alpha_1+\gamma_2-d/q|}(1+|x|/t)^{|\alpha_3-\delta_2|}&:&0<t\leq 1\\
            (1+|x|)^{|\alpha_3-\delta_2|}&:&t=\infty\,.
         \end{array}\right.
\end{equation}
for $(x,t) \in X$. The same holds for $L^w_{B,q,a}(X)$ in replacement of $P^w_{B,q,a}(X)$.
\end{lemma}
\bproof Assume $j\geq 1$, $k\in \Z$ and choose $(x,t) \in U_{j,k}$ arbitrarily. By $(W2)$ we have
\begin{equation}\label{f3}
  \begin{split}  \nonumber
    \|\chi_{U_{j,k}}|P^{w}_{B,q,a}\| &\asymp t^{-d/q}w(x,t)\Big\|\sup\limits_{z\in \R}\frac{\chi_{j,k}(\cdot+z)}{(1+|z|/t)^a}|B(\R)\Big\|\\
    &\gtrsim t^{-d/q+\alpha_1}(1+|x|/t)^{-\alpha_3}\|\chi_{j,k}|B(\R)\|\\
    &\gtrsim t^{-d/q+\alpha_1+\gamma_2}(1+|x|/t)^{\delta_2-\alpha_3}\,.
  \end{split}  
\end{equation}
Hence, choosing $v_{w,B,q}(x,t)$ as above gives \eqref{eq-33}. Note that $v_{w,B,q} \geq 1$. In case $j=0$ the modifications are straightforward. 
\eproof

Our next result provides equivalent norms of the sequence spaces associated to the Peetre type function spaces on $X$.

\begin{Theorem} Let $1 \leq q \leq \infty$, $B(\R)$, $w$ as in Lemma \ref{chainembappl}, and $a>0$. Then 
\begin{equation}\label{sequ_1}
  \begin{split}  
    &\|\{\lambda_{j,k}\}_{j,k}|(P^{w}_{B,q,a})^{\sharp}\|\\ 
    &\asymp \Big\|
    \sup\limits_{z\in \R}\frac{1}{(1+|z|)^a}\sum\limits_{k\in \Z}w(\alpha k,\infty)|
    \lambda_{0,k}|\chi_{0,k}(\cdot+z)|B(\R)\Big\|\\
    &~~~+\Big\|\Big(\sum\limits_{j\in \n}\Big[\sup\limits_{z\in \R}\frac{\beta^{dj/q}}{(1+\beta^j|z|)^a}\sum\limits_{k\in \Z}w(\alpha k \beta^{-j},\beta^{-j})|\lambda_{j,k}|\chi_{j,k}(\cdot+z)\Big]^q\Big)^{1/q}|B(\R)\Big\|
  \end{split}  
\end{equation}
and
\begin{equation}\label{sequ_2}
  \begin{split}  
    &\|\{\lambda_{j,k}\}_{j,k}|(L^{w}_{B,q,a})^{\sharp}\|\\ 
    &\asymp \Big\|
    \sup\limits_{z\in \R}\frac{1}{(1+|z|)^a}\sum\limits_{k\in \Z}w(\alpha k,\infty)|
    \lambda_{0,k}|\chi_{0,k}(\cdot+z)|B(\R)\Big\|\\
    &~~~+\Big(\sum\limits_{j\in \n}\Big\|\sup\limits_{z\in \R}\frac{\beta^{dj/q}}{(1+\beta^j|z|)^a}\sum\limits_{k\in \Z}w(\alpha k \beta^{-j},\beta^{-j})|\lambda_{j,k}|\chi_{j,k}(\cdot+z)|B(\R)\Big\|^q\Big)^{1/q}\,.
  \end{split}  
\end{equation}
Additionally, we have $(L^{w}_{B,q,a})^{\sharp} = (L^{w}_{B,q,a})^{\flat}$ and 
$(P^{w}_{B,q,a})^{\sharp} = (P^{w}_{B,q,a})^{\flat}$, respectively.
 
\end{Theorem}
\bproof According to Definition \ref{absdefss} the statement is a result of a straightforward computation taking (W2) into account.
\eproof

If we have additional knowledge on the space $B(\R)$, then the structure of 
the sequence spaces $(P^{w}_{B,q,a})^{\sharp}$ 
and $(L^{w}_{B,q,a})^{\sharp}$
simplifies significantly. Indeed, under some additional conditions (see below) they coincide
with the spaces  $p_{B,q}^w$ and $\ell_{B,q}^w$ of which the norms are given by % i.e., they coincide with the spaces 
\begin{equation}\nonumber
  \begin{split}  
    \|\{\lambda_{j,k}\}_{j,k}|p^{w}_{B,q}\|
    &= \Big\|
    \sum\limits_{k\in \Z}w(\alpha k,\infty)|
    \lambda_{0,k}|\chi_{0,k}|B(\R)\Big\|\\
    &~~~+\Big\|\Big(\sum\limits_{j\in \n}\beta^{dj}\Big[\sum\limits_{k\in \Z}w(\alpha k \beta^{-j},\beta^{-j})|\lambda_{j,k}|\chi_{j,k}\Big]^q\Big)^{1/q}|B(\R)\Big\|\,,
  \end{split}  
\end{equation}
and 
\begin{equation}\nonumber
  \begin{split}  
    \|\{\lambda_{j,k}\}_{j,k}|\ell^{w}_{B,q}\|
    &\asymp \Big\|
    \sum\limits_{k\in \Z}w(\alpha k,\infty)|
    \lambda_{0,k}|\chi_{0,k}|B(\R)\Big\|\\
    &~~~+\Big(\sum\limits_{j\in \n}\beta^{dj}\Big\|
    \sum\limits_{k\in \Z}w(\alpha k \beta^{-j},\beta^{-j})|\lambda_{j,k}|\chi_{j,k}|B(\R)\Big\|^q\Big)^{1/q}\,,
  \end{split}  
\end{equation}
respectively, and get therefore independence of $a$. Before giving a precise statement we first introduce the Hardy--Littlewood maximal function $M_r f$, $r>0$. It is defined for $f\in L_1^{loc}(\R)$ via
\begin{equation}\nonumber
  (M_rf)(x) = \sup\limits_{x\in Q} \Big(\frac{1}{|Q|}\int_{Q}\,|f(y)|^r\,dy\Big)^{1/r}\quad,\quad x\in\R\,,
\end{equation}
where the $\sup$ runs over all rectangles $Q$ containing $x$ with sides parallel to the coordinate axes. 
The following majorant property of the Hardy--Littlewood maximal function is taken from \cite[II.3]{StWe71}.

\begin{lemma}\label{majprop} Let $f\in L_1^{loc}(\R)$ and $\varphi \in L_1(\R)$ where $\varphi(x) = \psi(|x|)$ with 
a nonnegative decreasing function $\psi:[0,\infty)\to \re$. Then we have 
$$
   |(f \ast \varphi)(x)| \leq (M_1f)(x)\|\varphi|L_1(\R)\|
$$
for all $x\in \R$\,.
\end{lemma}
\bproof A proof can be found in \cite[II.3]{StWe71}, page 59. 
\eproof
\newline\\Let us further define the space $B(\ell_q,\R)$ as the space of all sequences of measurable functions $\{f_k\}_{k\in I}$ on $\R$ such that 
$$
    \|\{f_k\}_{k\in I}|B(\ell_q,\R)\| := \Big\|\Big(\sum\limits_{k\in I}|f_k|^q\Big)^{1/q}|B(\R)\Big\| < \infty.
$$
\begin{corollary}\label{corsequ} Let $1\leq q\leq \infty$, $a>0$, and $B(\R)$, $w$ as above.\\\newline
{\em (i)} If for some $r>0$ with $ar>d$ the Hardy--Littlewood maximal operator 
$M_r$ is bounded on $B(\R)$ and on $B(\ell_q,\R)$ then $(L^{w}_{B,q,a})^{\sharp}=\ell^w_{B,q}$ and $(P^{w}_{B,q,a})^{\sharp} = p^w_{B,q}$,
respectively. \\\newline
{\em (ii)} If $B(\R)$ is uniformly translation invariant, see \eqref{eq31}, then  $(L^w_{B,q,a})^{\sharp} = (L^w_{B,q})^{\sharp} = \ell^w_{B,q}$ provided
$a>d$\,.
\end{corollary}
\bproof For $(x,t) \in U_{j,k}$, we have
\begin{equation}\label{eq-0}
      \sup\limits_{z}\frac{|\chi_{j,k}(w+z)|}{  (1+\beta^j|z|)^{ar}} \lesssim
      \frac{1}{(1+|w-x|/t)^{ar}}\lesssim \Big(\chi_{j,k} \ast \frac{t^{-d}}{(1+|\cdot|/t)^{ar}}\Big)(w)\,.
\end{equation}
Indeed, the first estimate is obvious. Let us establish the second one
\begin{align}
\nonumber
  %\begin{split}
    & \Big(\chi_{j,k} \ast \frac{1}{(1+|\cdot|/t)^{ar}}\Big)(w) =
     \int_{Q_{j,k}}
     \frac{1}{(1+|w-y|/t)^{ar}} \,dy
     \gtrsim \int_{|y|\leq ct}
    \frac{1}{(1+|w-x-y|/t)^{ar}}\,dy\notag\\
    &= \int_{|y|\leq ct}
    \frac{1}{(1+|w-x|/t+|y|/t)^{ar}}\,dy
    \gtrsim t^d\int_0^{1}
    \frac{s^{d-1}}{(1+|w-x|/t+s)^{ar}}\,ds
    \gtrsim \frac{t^d}{(1+|w-x|/t)^{ar}}\,.\nonumber
  % \end{split}
\end{align}
Because of $ar>d$ the functions $g_j = \beta^{jd}(1+\beta^j|\cdot|)^{-ar}$ belong to $L_1(\R)$ and the $L_1(\R)-$norms are uniformly bounded
in $j$.

(i) We use Lemma \ref{majprop} in order to estimate the convolution on the right-hand side of \eqref{eq-0} by the Hardy--Littlewood maximal function and obtain
$$
   \sup\limits_{z}\frac{|\chi_{j,k}(x+z)|}{(1+\beta^j|z|)^{ar}} \lesssim M_1(\chi_{j,k})(x)\quad,\quad x\in \R\,.
$$
Hence, we can rewrite \eqref{sequ_1} as
\begin{equation}\nonumber
  \begin{split}  
    &\|\{\lambda_{j,k}\}_{j,k}|(P^{w}_{B,q,a})^{\sharp}\|
    \lesssim \Big\|
    M_r\Big[\sum\limits_{k\in \Z}w(\alpha k,\infty)|\lambda_{0,k}|\chi_{0,k}\Big](\cdot)|B(\R)\Big\|\\
    &~~~+\Big\|\Big(\sum\limits_{j\in \n}\Big[M_r\sum\limits_{k\in \Z}\beta^{dj/q}w(\alpha k \beta^{-j},\beta^{-j})|\lambda_{j,k}|\chi_{j,k}
    \Big]^q\Big)^{1/q}|B(\R)\Big\|\,.
  \end{split}  
\end{equation}
Since by assumption $M_r$ is bounded on $B(\ell_q,\R)$ we obtain the desired upper estimate. The corresponding estimate from below is trivial. 
The proof of the coincidence $(L^w_{B,q,a})^{\sharp} = \ell^w_{B,q}$ is similar. %\\\newline

For the proof of (ii) we do not need the Hardy--Littlewood maximal function. We use \eqref{eq-0} with $r=1$ and simply Minkowski's inequality. 
This yields
\begin{equation}\nonumber
 \begin{split}
    &\Big(\sum\limits_{j\in \n}\Big\|\sup\limits_{z\in \R}\frac{\beta^{dj/q}}{(1+\beta^j|z|)^a}\sum\limits_{k\in \Z}w(\alpha k \beta^{-j},\beta^{-j})|\lambda_{j,k}|\chi_{j,k}(\cdot+z)|B(\R)\Big\|^q\Big)^{1/q}    \\
    &\lesssim \Big(\sum\limits_{j\in \n}\Big\|\Big(\sum\limits_{k\in \Z}\beta^{jd/q}w(\alpha k \beta^{-j},\beta^{-j})|\lambda_{j,k}|\chi_{j,k}(\cdot)\Big)\ast g_j(\cdot)|B(\R)\Big\|^q\Big)^{1/q}\\
    &\lesssim \Big(\sum\limits_{j\in \n}\Big[\int_{\R}g_j(y)\Big\|\sum\limits_{k\in \Z}\beta^{jd/q}w(\alpha k \beta^{-j},\beta^{-j})|\lambda_{j,k}|\chi_{j,k}(\cdot-y)|B(\R)\Big\|\,dy\Big]^q\Big)^{1/q}\\
    &\lesssim \Big(\sum\limits_{j\in \n}\beta^{jd}\Big\|\sum\limits_{k\in \Z}w(\alpha k \beta^{-j},\beta^{-j})|\lambda_{j,k}|
    \chi_{j,k}(\cdot-y)|B(\R)\Big\|^q\Big)^{1/q}\,.
  \end{split}
\end{equation}
The same argument works for the first summand in \eqref{sequ_2}. The estimate from below is trivial. \eproof

\begin{remark} The main examples for spaces $B(\R)$ satisfying the assumptions in Corollary \ref{corsequ} are ordinary Lebesgue spaces $L_p(\R)$, $1 \leq p \leq \infty$, Muckenhoupt weighted Lebesgue space $L_p(\R,v)$, $1\leq p\leq\infty$, and Morrey spaces $M_{u,p}(\R)$, $1\leq p\leq u\leq \infty$, defined in Subsection \ref{Morrey}. 
\end{remark}

\begin{remark} Corollary \ref{corsequ}(ii) remains valid if we weaken condition \eqref{eq31} in the following sense,
$$
    \|T_x:B(\R)\to B(\R)\| < (1+|x|)^{\eta}
$$
for some $\eta>0$. One has to adjust the parameter $a$ in this case. 
This setting applies to certain weighted $L_p$-spaces $B(\R) = L_p(\R,\omega)$ with polynomial 
weight $\omega(y) = (1+|y|)^{\kappa}$.
\end{remark}

\subsection{The coorbits of $P^{w}_{B,q,a}(X)$ and $L^{w}_{B,q,a}(X)$}
\label{coorbits}
Now we apply the abstract coorbit space theory from Section \ref{abstrth} to our concrete setup. %The situation is as follows. 
We put $\mathcal{H} = L_2(\R)$ and fix an admissible continuous wavelet frame $\cf$ in the sense of Definition \ref{admfr}.
According to the abstract theory in Section \ref{abstrth}
the operator $R_{\cf}$ is then given by
$$
    R_{\cf}((x,t),(y,s)) = \langle\varphi_{(x,t)},
    \varphi_{(y,s)}\rangle\quad,\quad (x,t),(y,s) \in X\,.
$$
The relevant properties of this kernel and the kernels below depend on smoothness and decay conditions of the wavelets, see Definition \ref{basedef}.
The next result plays a crucial role and is proved in \cite[Lem.\ A.3]{T10}\,. 
Similar results which are stated in a different language can be found for instance in \cite[Lem.\ B1, B2]{FrJa90}, \cite[Lem.\ 1]{Ry99a}, and \cite[Lem.\ 1.2.8, 1.2.9]{HeNe07}.

\begin{lemma}\label{help1} Let $L\in \n$, $K>0$, and $g,f, f_0 \in L_2(\R)$.
   \begin{description}
   \item(i) Let $g$ satisfy $(D)$, $(M_{L-1})$ and let $f_0$ satisfy $(D)$, $(S_{K})$. Then for every 
   $N\in \N$ there exists a constant $C_N$ such that the
   estimate 
   $$
       |(W_g f_0)(x,t)| \leq C_N \frac{t^{\min\{L,K\}+d/2}}{(1+|x|)^N}
   $$
   holds true for $x\in \R$ and $0<t<1$\,.
   \item(ii) Let $g,f$ satisfy $(D)$, $(M_{L-1})$ and $(S_K)$. For every $N\in \N$ there exists a constant $C_N$ such that the
   estimate 
   $$
       |(W_g f)(x,t)| \leq C_N
       \frac{t^{\min\{L,K\}+d/2}}{(1+t)^{2\min\{L,K\}+d}}\Big(1+\frac{|x|}{1+t}\Big)^{-N}
   $$
   holds true for $x\in \R$ and $0<t<\infty$\,.
   \end{description}
\end{lemma}

Based on this lemma, we can show that the kernel $R_{\cf}$ acts continuously on $P^{w}_{B,q,a}(X)$.

\begin{lemma}\label{R(Y)} Let $\cf$ be an admissible continuous wavelet frame. Then the 
operator $R_\cf$ belongs to $\mathcal{B}_{Y,m}$ for 
$Y = P^w_{B,q,a}(X)$ or $Y = L^w_{B,q,a}(X)$ and every $v,m$ given by 
\eqref{weight} and \eqref{mw}. Moreover, it holds $R_{\cf}(Y) \subset L_{\infty}^{1/v_{w,B,q}}(X,\mu)$
where $v_{w,B,q}$ is defined in \eqref{rem1}.
\end{lemma}
\bproof We use that 
$$
    R_{\cf}((x,t),(y,s)) = \left\{\begin{array}{rcl}
                               W_{\Phi}\Phi\Big(\frac{y-x}{t},\frac{s}{t}\Big)&:&t,s \in (0,1]\;,\\
                               W_{\Phi}\Phi_0(y-x,s)&:& t=\infty, s\in (0,1]\;,\\
                               W_{\Phi_0}\Phi\Big(\frac{y-x}{t},\frac{1}{t}\Big)&:& t\in (0,1], s=\infty\;,\\
                               W_{\Phi_0}\Phi_0(y-x,1)&:&t=s=\infty\;,
                            \end{array}\right.
$$
where the operator $W$ denotes the continuous wavelet transform, see Subsection \ref{sectCWT}. 
Together with Propositions \ref{bkernel}, \ref{bddalgebra} in combination with Lemma \ref{help1}, this yields
that $R_{\cf}$ belongs to $\mathcal{B}_{Y,m}$\,. The embedding $R(Y) \subset L_{\infty}^{1/v_{w,B,q}}(X,\mu)$
follows from the abstract result in Lemma \ref{Linftyemb} and the choice of the weight $v_{w,B,q}$ in Lemma \ref{chainembappl}\,. 
To prove that $\cf$ satisfies the 
property $D[1,1,Y]$ we refer to Subsection \ref{subsectdiscr} below and Proposition \ref{bddkernels}.\eproof

Now we are ready to define the coorbits $\CoP$ and $\CoL$. 

\begin{definition}\label{defcoorbit} Let $1 \leq q \leq \infty$, $B(\R)$ and $w$ as above, $\cf$ be an admissible 
continuous frame in the sense of Definition \ref{admfr}, and $a>0$. We define
\begin{equation}\nonumber
  \begin{split}  
    \CoP &= \Co(P^w_{B,q,a},\cf):= \{f\in (\mathcal{H}^1_{v_{w,B,q}})^{\sim}~:~V_{\mathcal{\cf}}f \in P^w_{B,q,a}(X)\}\,,\\
    \CoL &= \Co(L^w_{B,q,a},\cf):=\{f \in (\mathcal{H}^1_{v_{w,B,q}})^{\sim}~:~V_{\mathcal{\cf}}f \in L^w_{B,q,a}(X)\}\,.
  \end{split}  
\end{equation}
\end{definition}

Based on the abstract theory we immediately obtain the following basic properties
of the introduced coorbit spaces.

\begin{Theorem}\label{props} Let $1 \leq q \leq \infty$, $a>0$, $w\in \mathcal{W}^{\alpha_3}_{\alpha_1,\alpha_2}$, $\cf$ be
an admissible frame, and let $B(\R)$ satisfy (B1) and (B2). Then we have the following  properties. 
\begin{description}
 \item(a) If $a>0$ then the spaces $\CoL$ and $\CoP$ are Banach spaces. 
 \item(b) A function $F \in P^w_{B,q,a}$ (or $L^w_{B,q,a}$) is of the form $V_{\cf} f$ for some $f\in \CoP$ (or $\CoL$) if and only if $F = R_{\cf}(F)$.
 \item(c) The spaces $\CoP$ and $\CoL$ do not depend on the frame $\cf$ in the sense that a different admissible frame 
 in the sense of Definition \ref{admfr} leads to the same space. Furthermore, if we use a weight of the form
 (\ref{weight}) satisfying $v(x) \ge v_{w,B,q}(x)$ 
 then the corresponding spaces coincide as well. We also have %even have the coincidence
 $$
    \CoP = \{f\in \mathcal{S}'(\R)~:~V_{\mathcal{F}}f \in P^w_{B,q,a}\}\,,
 $$
and similarly for $\CoL$\,.
\end{description}
\end{Theorem}

\bproof Assertions (a), (b) follow from Proposition 2(a),(b) in \cite{fora05} and Lemma \ref{R(Y)}\,. The assertion in (c)
is a consequence of the abstract independence results in Lemmas \ref{ind1}, \ref{ind2} together with Proposition \ref{bkernel}.

\subsection{Discretizations}
\label{subsectdiscr}

In the following we use a covering 
$\mathcal{U} = \mathcal{U}^{\alpha,\beta} = \{U_{j,k}\}_{j,k}$ as introduced in 
Subsection \ref{seq}\,. 
\begin{definition} The oscillation kernels
$\osc_{\alpha,\beta}$ and
$\osc^{\ast}_{\alpha,\beta}$ are given as
follows
$$
    \osc_{\alpha,\beta}((x,t),(y,s)) = \sup\limits_{(z,r) \in Q_{(y,s)}}
    |R_{\cf}((x,t),(y,s))-R_{\cf}((x,t),(z,r))|\,,
$$
where $Q_{(y,s)} = \bigcup\limits_{(j,k):(y,s)\in U_{j,k}} U_{j,k}$, and 
$\osc_{\alpha,\beta}^*$ is its adjoint.
%$$
%    \osc^{\ast}_{\alpha,\beta}((x,t),(y,s)) = \sup\limits_{(z,r) \in Q_{(x,t)}}
%    |R_{\cf}((z,r),(y,s))-R_{\cf}((x,t),(y,s))|\,.
%$$
\end{definition}

Next, we show that the norms of these kernels can be made arbitrarily small by choosing a sufficiently fine covering.

\begin{proposition}\label{bddkernels} Let $\cf = \{\varphi_{\x}\}_{\x\in X}$ be an admissible continuous frame in the sense of Definition
\ref{admfr}.
\begin{description}
     \item(i) Let $\alpha_0 >0$ and $\beta_0>1$ be arbitrary. The kernels $\osc_{\alpha,\beta}$ and $\osc^{\ast}_{\alpha, \beta}$
     with $0<\alpha\leq \alpha_0$ and $1<\beta\leq \beta_0$ are uniformly bounded operators on
     $P^{w}_{B,q,a}(X)$ %to $P^{w}_{B,q,a}(X)$ 
     and belong to $\mathcal{A}_{m_v}$ for every weight $v$ of the form \eqref{weight}.
     \item(ii) If $\alpha\downarrow 0$ and $\beta\downarrow 1$ then
     $$
        \|\osc_{\alpha,\beta}:P^{w}_{B,q,a} \to P^{w}_{B,q,a}\| \to 0\quad,\quad
        \|\osc^{\ast}_{\alpha,\beta}:P^{w}_{B,q,a} \to
        P^{w}_{B,q,a}\| \to 0\,.
     $$
  \end{description}
\end{proposition}

\bproof Because of the particular structure of $\cf$, see Definition \ref{admfr}, we guarantee that 
$\Phi$ satisfies $(D), (M_L)$ and $(S_L)$ and that $\Phi_0$ satisfies
$(D)$ and $(S_L)$ for all $L>0$. Putting $G_1(y,s) =
W_{\Phi}\Phi(y,s)$, $G_2(y,s) = (W_{\Phi_0}\Phi)(y,s)$,
$G_3(y,s) = (W_{\Phi}\Phi_0)(y,s)$ and $G_4(y) =
(W_{\Phi_0}\Phi_0)(y,1)$ then Lemma \ref{help1} yields the following
estimates for every $L>0$ and every $N\in \N$
\begin{equation}\label{f17}
    |G_i(y,s)| \leq C_N \frac{s^{\alpha_i}}{(1+s)^{\beta_i}}\Big(1+\frac{|y|}{1+s}\Big)^{-N}\quad,\quad y\in \R,s\in (0,1], i=1,2,3\,,
\end{equation}
where $\alpha_1 = L+d/2, \beta_1 = 2L+d$, $\alpha_2 = 0, \beta_2 = L+d/2$ and $\alpha_3 = L+d/2, \beta_3 = 0$. Moreover, we
have 
$$
    |G_4(y)| \leq C_N(1+|y|)^{-N}\,.
$$
For $K=\mbox{osc}$ we choose the set $U =
[-\alpha,\alpha]^d\times [\beta^{-1},\beta]$ and $U_0 =
[-\alpha,\alpha]^d$ and use the functions from $\eqref{Kb}$,
\begin{equation}\nonumber
 \begin{split}
    G_1^{\sharp}(y,s) &= \sup\limits_{(z,r) \in
    (y,s)\cdot U}|G_1(y,s)-G_1(z,r)|\;,\\
    G_2^{\sharp}(y,s) &= \sup\limits_{(z,s) \in (y,s)\cdot [U_0 \times \{1\}]}
    |G_2(y,s)-G_2(z,s)|\;,\\
    G_3^{\sharp}(y,s) &= \sup\limits_{(z,r) \in (y,s)\cdot
    U}|G_3(y,s)-G_3(z,r)|\;,\\
    G_4^{\sharp}(y) &= \sup\limits_{z\in y+U_0}
    |G_4(y)-G_4(z)|\,,
 \end{split}
\end{equation}
%
%in (\ref{Kb}). 
where $(y,s)\cdot U = \{(y+sx,st)~:~(x,t) \in U\}$.
Clearly, the functions $G_i^{\sharp}$ depend on
$\alpha,\beta$ and obey a similar behavior as the functions
$G_i$ in (\ref{f17}) and moreover, the functions appearing in
(\ref{bddcond}) possess this behavior for the same reason.
The integrals in (\ref{f17}) are uniformly bounded in $\alpha\leq
\alpha_0$ and $\beta \leq \beta_0$. Using
Propositions \ref{bkernel}, \ref{bddalgebra} we obtain (i) for $K=\osc$.\\ 
For the kernel $\osc^{\ast}$
we have to replace $G_i^{\sharp}(y,s)$ by
$G_i^{\sharp\sharp}(y,s)$ defined via
\begin{equation}\nonumber
 \begin{split}
    G_1^{\sharp\sharp}(y,s) &= \sup\limits_{(z,r)\in U^{-1}(y,s)}
    |G_1(y,s)-G_1(z,r)|\,,\\
    G_2^{\sharp\sharp}(y,s) &= \sup\limits_{(z,r) \in U^{-1}(y,s)}
    |G_2(y,s)-G_2(z,r)|\,,\\
    G_3^{\sharp\sharp}(y,s) &= \sup\limits_{(z,s) \in [U_0 \times \{1\}](y,s)}|G_3(y,s)-G_3(z,s)|\,,\\
    G_4^{\sharp\sharp}(y) &= \sup\limits_{z\in y+U_0}
    |G_4(y)-G_4(z)|\,,
 \end{split}
\end{equation}
where $U^{-1} = \{(-x/t,1/t)~:~(x,t) \in U\}$ and $U\cdot (y,s) = \{(x+ty,st)~:~(x,t) \in U\}$.
Analogous arguments give (i) for $\mbox{osc}^{\ast}$. For the proof of (ii) we use the continuity of the functions $G_i$ and
argue analogously as in \cite[Lem.\ 4.6(ii)]{Gr91}.
\eproof

Let us state the first discretization result. 

\begin{Theorem}\label{discpeetre} Let $1\leq q \leq \infty$, $a>0$, $w\in \mathcal{W}^{\alpha_3}_{\alpha_1,\alpha_2}$, $B(\R)$
satisfying (B1) and (B2), and $\cf = \{\varphi_\x\}_{\x\in X}$ be an admissible continuous  wavelet frame. There exist $\alpha_0>0$ 
and $\beta_0>1$,
such that for all $0<\alpha\leq \alpha_0$ and $1< \beta\leq \beta_0$ there is a discrete wavelet frame 
$\cf_d = \{\varphi_{\x_{j,k}}\}_{j\in \N_0, k \in \Z}$ with $\x_{j,k} = (\alpha k \beta^{-j}, \beta^{-j})$
and a corresponding dual frame $\mathcal{E}_d = \{e_{j,k}\}_{j\in \N_0, k\in \Z}$ 
such that 
\begin{description}
    \item(a)
    \begin{equation*}
      \begin{split}  
        \|f|\CoP\| &\asymp \|\{\langle f,\varphi_{x_{j,k}}\rangle\}_{j\in \N_0, k\in \Z}|(P^w_{B,q,a})^{\sharp}\|
                   \asymp \|\{\langle f,e_{j,k} \rangle\}_{j\in \N_0, k\in \Z}|(P^w_{B,q,a})^{\sharp}\|\,.
      \end{split}  
    \end{equation*}
    \item(b) If $f\in \CoP$ then the series
    \begin{equation*}
      \begin{split}  
        f &= \sum\limits_{j\in \N_0}\sum\limits_{k\in \Z}\langle f,e_{j,k}\rangle \varphi_{x_{j,k}}
        =\sum\limits_{j\in \N_0}\sum\limits_{k\in \Z}\langle f,\varphi_{x_{j,k}}\rangle e_{j,k}
      \end{split}  
    \end{equation*}
    converge unconditionally in the norm of $\CoP$ if the finite sequences are dense in 
    $(P^w_{B,q,a})^{\sharp}$ and with
    weak$^\ast$--convergence induced by $(\mathcal{H}^1_v)^{\sim}$ otherwise.
\end{description}
\end{Theorem}
\bproof The assertion is a consequence of our abstract Theorem \ref{discr2}. Due to the choice of $v = v_{w,B,q}$, see \eqref{rem1}, we know by Lemma \ref{R(Y)} that $R_{\cf}$ belongs to $\mathcal{B}_{Y,m}$ for $Y=P^w_{B,q,a}$, and that $R_{\cf}(P^w_{B,q,a}) \hookrightarrow L^{1/v}_{\infty}(X)$. Hence, we have that $\cf$ satisfies $(F_{v,Y})$ for $Y = P^{w}_{B,q,a}(X)$. As a consequence of Proposition \ref{bddkernels}, the kernels $\mbox{osc}_{\alpha,\beta}$ and $\mbox{osc}^{\ast}_{\alpha,\beta}$ are 
bounded operators from $Y$ to $Y$, and the norms of $\mbox{osc}_{\alpha,\beta}$ and $\mbox{osc}^{\ast}_{\alpha,\beta}$ %even 
tend to zero when $\alpha \to 0$ and $\beta \to 1$. 
%This is mainly a consequence of Proposition \ref{bkernel} and Lemma \ref{help1}. 
Choosing $v \geq v_{w,B,q}$ and the weight $m$ accordingly, 
we obtain by analogous arguments that the norms of $\mbox{osc}_{\alpha,\beta}$ and $\mbox{osc}^{\ast}_{\alpha,\beta}$ in $\mathcal{A}_m$ tend to zero. 
Therefore, we have $\mathcal{F} \in D[\delta,m,Y]$ for every $\delta>0$\,. In particular, $\mathcal{F}$
satisfies $D[1,1,Y]$\,. 
\eproof

\subsection{Wavelet bases}

In the sequel we are interested in the discretization of coorbits with respect to Peetre type spaces
via $d$-variate wavelet bases of the following type. 
%Corresponding to the discrete situation in 
According to Lemma \ref{dwavelet}
we start with a scaling function
$\psi^0$ and wavelet $\psi^1$ belonging to $L_2(\re)$. 
Let further $E = \{0,1\}^d$. For $c\in E$ we define the function $\psi^c:\R \to \re$ by the tensor product
$\psi^c = \bigotimes_{i=1}^d \psi^{c_i}$, i.e., $\psi^c(x) =  \prod_{i=1}^d \psi^{c_i}(x_i)$. The frame $\Psi^c$ on $X$ is given by $\Psi^c = \{\psi^c_\z\}_{\z\in X}$, where for $c\neq 0$
$$
    \psi^c_{(x,t)} = \left\{\begin{array}{rcl}
                        T_x\mathcal{D}^{L_2}_t \psi^c&:&0<t<1\,,\\
                        T_x \psi^c&:& t=\infty\,,
                     \end{array}\right.
$$ and 
$$
    \psi^0_{(x,t)} = \left\{\begin{array}{rcl}
                        0&:&0<t<1\,,\\
                        T_x \psi^0&:& t=\infty\,.
                     \end{array}\right.
$$
This construction leads to a family of continuous systems $\Psi^c,\,c\in E$. Our aim is to apply Theorem
\ref{wbases2} in order to achieve wavelet basis characterizations of the Besov-Lizorkin-Triebel type spaces $\CoP$ and $\CoL$.
In order to apply the abstract result in Theorem \ref{wbases2} we have to consider the Gramian cross kernels $K_c$ and $K^{\ast}_c$ related to the covering $\mathcal{U}^{\alpha,\beta}$ defined by
$$
  K_c(\x,\y) = \sup\limits_{\z\in Q_\x} |G(\cf,\cg_c)(\z,\y)|, \quad \x,\y \in X\,, c \in E\,,
$$
and $K^{\ast}_c(\x,\y) = K_c(\y,\x)$, see \eqref{eq-28} and \eqref{eq-26}.

\begin{lemma} \label{wavres} Let $1\leq q \leq \infty$, $a>0$, 
$w\in \mathcal{W}^{\alpha_3}_{\alpha_1,\alpha_2}$, $B(\R)$
satisfying (B1) and (B2). Let further $\cf$ be an admissible continuous frame, $\cg_c$ be the frames from above, and $K_c, K^{\ast}_c$, $c \in E$, the corresponding Gramian cross kernels. The weight $v_{w,B,q}$ is given by \eqref{rem1} and 
$m_{w,B,q}$ denotes its associated weight. Assume further that the functions $\psi^0, \psi^1$ satisfy $(D)$ and $(S_K)$, and that $\psi^1$ also satisfies $(M_{L-1})$. 
\begin{itemize}
\item[(i)] Under the assumption
\begin{equation}
\label{eq-27}
    K,L>\max\Big\{\frac{d}{2}+|\alpha_3-\delta_2|+\Big|\frac{d}{q}-\alpha_1-\gamma_2\Big|, 
    -\frac{d}{2}+2|\alpha_3-\delta_2|+\Big|\frac{d}{q}-\alpha_1-\gamma_2\Big|\Big\}
\end{equation}
we have $K_c,K^{\ast}_c \in \mathcal{A}_{m_{w,B,q}}$ for all $c\in E$\,.
\item[(ii)] If
\begin{equation}\label{eq-70}
  K,L >  \max\Big\{\frac{d}{q}-\frac{d}{2}+\alpha_2+a, \alpha_1+2a+\frac{d}{2}-\frac{d}{q}\Big\}
\end{equation}
then the kernels $K_c, K^{\ast}_c$ define bounded operators from $P^w_{B,q,a}(X)$ to $P^w_{B,q,a}(X)$.
\end{itemize}
\end{lemma}
\bproof We start with $c \in E, c \neq 0$. 
The following is analogous to the treatment of $\mbox{osc}^{\ast}$ in Proposition \ref{bddkernels}. 
As before we use the sets $U =
[-\alpha,\alpha]^d\times [\beta^{-1},\beta]$ and $U_0 =
[-\alpha,\alpha]^d$. A straightforward computation (analogously to the proof
of Proposition \ref{bddkernels}) gives the bounds \eqref{Kb} for the kernel $K_c$ with 
\begin{equation}\nonumber
 \begin{split}
    G^c_1(y,s) &= \sup\limits_{(z,r) \in U^{-1}(y,s)}|(W_{\Phi}\psi^c)(z,r)|\quad,\quad 0<s<\infty, y\in \R\,,\\
    G^c_2(y,s) &= \sup\limits_{(z,r)\in U^{-1}(y,s)}|(W_{\Phi_0}\psi^c)(z,r)|\quad,\quad 1<s<\infty, y\in \R\,,\\
    G^c_3(y,s) &= \sup\limits_{(z,s)\in [U_0\times \{1\}](y,s)} |(W_{\Phi}\psi^c)(z,s)|\quad,\quad 0<s<1, y\in \R\,,\\
    G^c_4(y)   &= \sup\limits_{z\in U_0+y} |(W_{\Phi_0}\psi^c)(z,1)|\quad,\quad y \in \R\,.
 \end{split}
\end{equation}
See the proof of Proposition \ref{bddkernels} for the used notation. Since $\psi^c$ satisfies $(D)$,$(S_K)$, $(M_{L-1})$, and $\Phi$ satisfies $(M_J)$ for all $J\in \N$, we obtain with the help of Lemma \ref{help1} the following estimates, valid
 for all $N\in \N$,
\begin{equation}\label{eq-30}
   \begin{split}
      |G^c_1(y,s)| &\leq C_N \frac{s^{K+d/2}}{(1+s)^{2K+d}}\Big(1+\frac{|y|}{s+1}\Big)^{-N}\,,\\
      |G^c_2(y,s)| &\leq C_N s^{-(L+d/2)}\Big(1+\frac{|y|}{s}\Big)^{-N}\,,\\
      |G^c_3(y,s)| &\leq C_N \frac{s^{K+d/2}}{(1+|y|)^{N}}\,,\\
      |G^c_4(y)| &\leq \frac{C_N}{(1+|x|)^N}\,.
   \end{split}
\end{equation}
Now we consider the kernels $K^{\ast}_c$, $c \in E$. In this case we obtain \eqref{Kb} with
\begin{equation}\nonumber
 \begin{split}
    G^{\ast,c}_1(y,s) &= \sup\limits_{(z,r) \in (y,s)\cdot U}|(W_{\psi^c}\Phi)(z,r)|\quad,\quad 0<s<\infty, y\in \R\,,\\
    G^{\ast,c}_2(y,s) &= \sup\limits_{(z,s)\in (y,s)\cdot[U_0\times \{1\}]}|(W_{\psi^c}\Phi)(z,s)|\quad,\quad 1<s<\infty, y\in \R\,,\\
    G^{\ast,c}_3(y,s) &= \sup\limits_{(z,r)\in (y,s)\cdot U} |(W_{\psi^c}\Phi_0)(z,s)|\quad,\quad 0<s<1, y\in \R\,,\\
    G^{\ast,c}_4(y)   &= \sup\limits_{z\in y+U_0} |(W_{\psi^c}\Phi_0)(z,1)|\quad,\quad y \in \R\,.
 \end{split}
\end{equation}
See again the proof of Proposition \ref{bddkernels} for the used notation.
The corresponding estimates are similar to \eqref{eq-30}, we just have to swap the role of $K$ and $L$. Hence, Proposition \ref{bddalgebra} implies that $K_c, K_c^{\ast}$ belong to $\mathcal{A}_{m_{w,B,q}}$ if \eqref{eq-27} is satisfied.        
Similar, Proposition \ref{bkernel} implies that the operators $K_c,K_c^{\ast}$ map $P^w_{B,q,a}$ boundedly into $P^w_{B,q,a}$ if 
\eqref{eq-70} is satisfied. In case $c = 0$ we have $G^0_1 = G^0_2 = 0$,
and $G_3(y,s)$, $G_4(y)$. The same conditions on $K$ and $L$ lead to the boundedness of the operators $K_0, K^{\ast}_0$. \eproof

Now we are ready for the discretization of $\CoP$ and $\CoL$ in terms of orthonormal wavelet bases. We only state the results for
$\Co P^{w}_{B,q,a}$. For $\Co L^{w}_{B,q,a}$ it is literally the same. 

\begin{Theorem}\label{wbasesP} Let $1\leq q \leq \infty$, $a>0$, $w\in \mathcal{W}^{\alpha_3}_{\alpha_1,\alpha_2}$, $B(\R)$
satisfying (B1) and (B2), and $\cf$ be an admissible continuous wavelet frame. Assume that 
$\psi^0, \psi^1 \in L_2(\re)$ generate an orthonormal wavelet basis of $L^2(\R)$
in the sense of Lemma \ref{dwavelet}
where $\psi^0$ satisfies $(D)$, $(S_K)$, and $\psi^1$ satisfies $(D)$, $(S_K)$, $(M_{L-1})$ such that 
\begin{equation}
  \begin{split}
     K,L>&\max\Big\{\frac{d}{2}+|\alpha_3 - \delta_2|+\Big|\frac{d}{q}-\alpha_1-\gamma_2\Big|,
     -\frac{d}{2}+2|\alpha_3-\delta_2|+\Big|\frac{d}{q}-\alpha_1-\gamma_2\Big|,\\
     &~~~~~~~~\frac{d}{q}-\frac{d}{2}+\alpha_2+a,\alpha_1+2a+\frac{d}{2}-\frac{d}{q}\Big\}\,.
  \end{split}
\end{equation}  
Then every $f\in \CoP$ has the decomposition 
\begin{equation}\label{eq-32}
   \begin{split}
        f  =&  \sum\limits_{c\in E}\sum\limits_{k\in \Z}\lambda^c_{0,k} \psi^c(\cdot-k)+\sum\limits_{c\in E\setminus \{0\}}\sum\limits_{j\in \N}\sum\limits_{k\in \Z}  \lambda^c_{j,k}2^{\frac{jd}{2}}\psi^c(2^j\cdot-k),
   \end{split}
\end{equation}
where the sequences $\lambda^c = \{\lambda^c_{j,k}\}_{j\in \N_0,k\in\Z}$ defined by 
$$
   \lambda^c_{j,k} = \langle f,2^{\frac{jd}{2}}\psi^c(2^j\cdot-k)\rangle\quad,\quad j\in \N_0,k\in \Z\,,
$$belong to the sequence space $(P^w_{B,q,a})^{\sharp} = (P^w_{B,q,a})^{\sharp}(\mathcal{U})$
for every $c\in E$, where $\mathcal{U} = \mathcal{U}^{1,2}$ is the covering introduced in Section
\ref{seq} with $\alpha=1,\beta=2$.

Conversely, an element $f\in (\mathcal{H}^1_{v_{w,B,q}})^{\sim}$ belongs to $\CoP$ if all sequences $\lambda^c(f)$ belong to $(P^w_{B,q,a})^{\sharp}({\mathcal U}^{1,2})$.
% with $\alpha = 1$ and $\beta = 2$.
The convergence in \eqref{eq-32} is in the norm of $\CoP$ if the finite sequences are dense in
$(P^w_{B,q,a})^{\sharp}$. In general, we have weak$^{\ast}$--convergence. 
\end{Theorem}
\bproof We apply the abstract Theorem \ref{wbases2}. First of all, Lemma \ref{R(Y)} shows that $\cf$ satisfies $(F_{v,Y})$ for $Y = P^w_{B,q,a}$ and $v = v_{w,B,q}$ given by \eqref{rem1}. The embedding $(P^w_{B,q,a})^{\sharp} \hookrightarrow (L_{\infty}^{1/v})^{\sharp}$ is ensured by the abstract Lemma \ref{Linftyemb} in combination with Lemma \ref{chainembappl}. The required boundedness of the Gramian kernels is showed in Lemma \ref{wavres}. \eproof

\begin{remark} The space $\Co L^{w}_{B,q}$ can be discretized in the same way.
According to Remark \ref{Tx} the corresponding conditions in Proposition \ref{bkernel}
are much weaker. The parameter $a$ is not needed here. We leave the
details to the reader.  
\end{remark}

\section{Examples -- generalized Besov-Lizorkin-Triebel spaces}
\label{examples}
The main class of examples are represented by 
the scales of Besov-Lizorkin-Triebel-Morrey spaces and
weighted Besov-Lizorkin-Triebel spaces, where 
we consider $B(\R)$ to be a Morrey space $M_{u,p}(\R)$, see Definition \ref{def:morrey} below, or a weighted Lebesgue
space $L_p(\R,v)$, see \eqref{wLp}. In the sequel we consider only weight functions $v$ such that (B2) is satisfied for $L_p(\R,v)$.

\subsection{Generalized $\mathbf{2}$-microlocal Besov-Lizorkin-Triebel spaces with Muckenhoupt weights}
A large class of examples is given by the scales of inhomogeneous 
$2$-microlocal Besov and Lizorkin-Triebel spaces
$B^w_{p,q}(\R,v)$ and $F^w_{p,q}(\R,v)$ on $\R$ with Muckenhoupt weights. 
These spaces represent a symbiosis of the spaces 
studied by Kempka \cite{Ke09,Ke10,Ke11}, Bui \cite{Bu82,Bu84}, Bui et al.\ 
\cite{BuPaTa96, BuPaTa97}, and Haroske et al.\  %, Piotrowska \cite{HaPi08}. 
\cite{HaPi08}. 
The scales $B^w_{p,q}(\R,v)$ and $F^w_{p,q}(\R,v)$
contain the classical inhomogeneous Besov-Lizorkin-Triebel spaces. 
For their definition, basic properties, and
results on atomic decompositions we mainly refer to Triebel's monographs \cite{Tr83, Tr92, Tr06,Tr08}. % and the recent one \cite{Tr08}.

Let us briefly recall the definition and some basic facts on Muckenhoupt weights. 
%We concentrate on weights belonging to some Muckenhoupt class $\mathcal{A}_p$. 
A locally integrable
function $v:\R \to \re_+$
belongs to $\mathcal{A}_p$, $1<p<\infty$, if the famous Muckenhoupt condition
$$
  \Big(\frac{1}{|B(y,r)|}\int_{B(y,r)}v(x)\,dx\Big)^{1/p}\cdot \Big(\frac{1}{|B(y,r)|}\int_{B} v(x)^{-p'/p}\,dx\Big)^{1/p'}\leq A,\quad \mbox{ for all }y\in \R, r >0\,,
$$
holds, where $1/p+1/p' = 1$ and $A$ is some constant independent of $y$ and $r$. 
The $\mathcal{A}_p$-condition implies the condition
\begin{equation}\label{eq45}
  \int_{\R}(t+|x-y|)^{-dp}v(y)\,dy \leq ct^{-dp}\int_{B(x,t)} v(y)\,dy\quad,\quad x\in \R, t>0\,,
\end{equation}
where c is independent of $x$ and $t$. See for instance \cite{Bu82} and the references given there. We further put 
$$
    \mathcal{A}_{\infty} := \bigcup\limits_{p>1} \mathcal{A}_p\,.
$$

\begin{lemma}\label{muckenh} Let $v \in {\mathcal A}_p$ for some $1<p<\infty$. 
Then the space $L_q(\R,v)$, $1\leq q\leq \infty$, satisfies property (B2) with $\gamma_1 = 0$,$\gamma_2 = dp/q$ and 
$\delta_1 = dp/q$ and $\delta_2 =-dp/q$.
\end{lemma}

\bproof Let $x \in \R$ and $0<t<1$ arbitrary. Let further $Q = [-1,1]^d$. We estimate using \eqref{eq45}
  \begin{equation}\nonumber
     \begin{split}
       t^{dp/q}t^{-dp/q}\Big(\int_{x+tQ} v(y)\,dy\Big)^{1/q} &\gtrsim t^{dp/q} \Big(\int_{\R} (1+|x-y|)^{-dp}v(y)\,dy\Big)^{1/q}\\ 
       &\gtrsim \frac{t^{dp/q}}{(1+|x|)^{dp/q}}\Big(\int_{\R} (1+|y|)^{-dp}v(y)\,dy\Big)^{1/q}\\
       &\gtrsim t^{dp/q}(1+|x|)^{-dp/q}\,,
     \end{split}  
  \end{equation}
  since $v$ is supposed to be locally integrable. 
  The estimate from above proceeds as follows, 
  \begin{equation}\nonumber
     \begin{split}
       \Big(\int_{x+tQ}v(y)\,dy\Big)^{1/q} &\lesssim \Big(\int_{B(0,|x|+\sqrt{d}t)}v(y)\,dy\Big)^{1/q}\\
       &\lesssim \Big(\int_{B(0,|x|+\sqrt{d}t)}(1+|y|)^{-dp}(1+|y|)^{dp}v(y)\,dy\Big)^{1/q}\\
       &\lesssim (1+|x|)^{dp/q}\Big(\int_{B(0,1)}v(y)\,dy\Big)^{1/q}\lesssim (1+|x|)^{dp/q}\,,
     \end{split}  
  \end{equation}
  where we used \eqref{eq45} in the last step. \eproof
    
The crucial tool in the theory of Muckenhoupt weights is the vector-valued Fefferman--Stein maximal inequality, see 
for instance \cite[Lem.\ 1.1]{Bu82} or \cite[Thm.\ 2.11]{HaPi08} and the references given there. 

\begin{lemma}\label{muckenmax} Let $1<p<\infty$, $1<q\leq \infty$, $v\in \mathcal{A}_p$, and $\{f_j\}_j$ be a sequence in $L_p(\R,v)$. Then 
$$
    \Big\|\Big(\sum\limits_{j}|M f_j|^q\Big)^{1/q}|L_p(\R,v)\Big\| \lesssim \Big\|\Big(\sum\limits_{j}|f_j|^q\Big)^{1/q}|L_p(\R,v)\Big\|\,.
$$
 
\end{lemma}
The definition of the spaces $B^w_{p,q}(\R,v)$ and $F^w_{p,q}(\R,v)$ relies on a dyadic decomposition of unity, see also \cite[2.3.1]{Tr83}.

\begin{definition}\label{decunity} Let $\Phi(\R)$ be the collection of all systems $\{\varphi_j(x)\}_{j\in
\n} \subset \mathcal{S}(\R)$ with the following properties:
\begin{description}
    \item(i) $\varphi_j(x) = \varphi(2^{-j}x)\quad,\quad j\in
    \N$\,,
    \item(ii) $\supp \varphi_0 \subset \{x\in \R~:~|x|\leq
    2\}\,,\quad
    \supp \varphi \subset \{x\in \R~:~1/2 \leq |x| \leq 2\}$\,, % and
    \item(iii) $\sum\limits_{j=0}^{\infty} \varphi_j(x) = 1$ for every $x\in \re$\,.
\end{description}
\end{definition}
\noindent
In order to define the spaces $F^{w}_{p,q}(\R,v)$ and $B^w_{p,q}(\R,v)$ 
%Recall Definition \ref{def:weight}
%and the weights $w \in \mathcal {W}^{\alpha_3}_{\alpha_1,\alpha_2}$. 
we define for a weight $w \in \mathcal{W}^{\alpha_3}_{\alpha_1,\alpha_2}$
a semi-discrete counterpart $\{w_j\}_{j \in \N_0}$, corresponding to an
%consider a semi-discrete counterpart which corresponds to the 
%definition of 
an admissible weight sequence in the sense of \cite{Ke09,Ke10,Ke11}. 
%Indeed, for $w \in \mathcal{W}^{\alpha_3}_{\alpha_1,\alpha_2}$ 
We put
\begin{equation}\label{def:wj}
  \begin{split}
    w_j(x) = \left\{\begin{array}{rcl}
                       w(x,2^{-j})&:&j\in \N\,,\\
                       w(x,\infty)&:&j=0\,.
                    \end{array}\right.
  \end{split}
\end{equation}

\begin{definition}\label{inhom} Let $v\in \mathcal{A}_{\infty}$ and 
$\{\varphi_j(x)\}_{j= 0}^{\infty} \in \Phi(\R)$.
Let further $w \in \mathcal{W}^{\alpha_3}_{\alpha_1,\alpha_2}$ with associated weight sequence
$\{w_j\}_{j\in \n}$ defined in \eqref{def:wj} 
and $0<q\leq \infty$. Put %$\Phi_j = \cf^{-1}\varphi_j$.
$\widehat{\Phi}_j = \varphi_j$.
\begin{description}
 \item(i) For $0<p\leq \infty$ we define (modification if $q=\infty$) 
 \begin{equation}\nonumber
  \begin{split}  
    B^w_{p,q}(\R,v) &= \Big\{f\in \mathcal{S}'(\R):\\
    &\|f|B^w_{p,q}(\R,v)\| = \Big(\sum\limits_{j=0}^{\infty}
    \|w_j(x)(\Phi_j \ast f)(x)|L_p(\R,v)\|^q\Big)^{1/q}<\infty\Big\}\,.
  \end{split}  
 \end{equation}
 \item(ii) For $0<p<\infty$ we define (modification if $q=\infty$) 
 \begin{equation}\nonumber
  \begin{split}  
    F^w_{p,q}(\R,v) &= \Big\{f\in \mathcal{S}'(\R):\\
    &\|f|F^w_{p,q}(\R,v)\| = \Big\|\Big(\sum\limits_{j=0}^{\infty}
    |w_j(x)(\Phi_j \ast f)(x)|^q\Big)^{1/q}|L_p(\R,v)\Big\|<\infty\Big\}\,.
  \end{split}
 \end{equation} 
\end{description}
\end{definition}

\begin{remark}  Let us discuss some special cases of the above defined scales. In the particular case $v\equiv 1$ and
$w(x,t) = t^{-s}$ we obtain the classical Besov-Lizorkin-Triebel spaces $B^s_{p,q}(\R)$ and $F^s_{p,q}(\R)$, see Triebel's
monographs \cite{Tr83,Tr92,Tr06} for details and historical remarks.
The choice $v\equiv 1$ leads to the generalized $2$-microlocal Besov-Lizorkin-Triebel spaces $B^{w}_{p,q}(\R)$ and $F^{w}_{p,q}(\R)$ studied systematically 
by Kempka \cite{Ke09,Ke10,Ke11}. The weight $w(x,t) = t^{-s}$ yields the 
Besov-Lizorkin-Triebel spaces with Muckenhoupt weights $B^s_{p,q}(\R,v)$ and $F^s_{p,q}(\R,v)$ already treated in 
Bui \cite{Bu82} and Haroske et al.\ %, Piotrowska
\cite{HaPi08} to mention just a few. 
 
\end{remark}

Unfortunately, this definition is not suitable %to %see that these spaces can be
to identify these spaces as certain coorbits. 
The connection to our spaces $\Co P^w_{B,q,a}(X)$ and $\Co L^w_{B,q,a}(X)$ is established by the theorem below. 
First, we replace the system $\{\varphi_j\}_{j}$ by a more general one and secondly, 
we prove a so-called continuous characterization, where we replace the discrete dilation parameter $j\in \n$ by $t>0$ and
the sums by integrals over $t$. Characterizations of this type have some history and are usually referred to as characterizations via 
local means. For further references and some historical facts we mainly refer to \cite{Tr88, Tr92, BuPaTa96, BuPaTa97, Ry99a} and
in particular to the recent contribution \cite{T10}, which provides a complete and self-contained reference. Essential for what follows are functions $\Phi_0, \Phi \in \mathcal{S}(\R)$ satisfying
the so-called Tauberian conditions 
\begin{equation}\label{condphi1}
\begin{split}
  |\widehat{\Phi}_0(x)| > 0 \quad &\mbox{ on }\quad \{|x| < 2\varepsilon\}\\
  |\widehat{\Phi}(x)| > 0 \quad&\mbox{ on }\quad \{\varepsilon/2<|x|< 2\varepsilon\}\,,
\end{split}
\end{equation}
for some $\varepsilon>0$, and
\begin{equation}\label{condphi2}
   D^{\bar{\alpha}} \widehat{\Phi}(0) = 0\quad\mbox{for all}\quad
   |\bar{\alpha}|_1 \leq R
\end{equation}
for some $R+1 \in \n $. If $R+1 = 0$ then the condition \eqref{condphi2} is void. We call the functions $\Phi_0$ and $\Phi$ kernels for local
means and use the notations $\Phi_k = 2^{kd}\Phi(2^k\cdot)$, $k\in \N$, as well as
$\Psi_t = \mathcal{D}_t \Psi = t^{-d}\Phi(\cdot/t)$, and the well-known classical Peetre maximal function
$$
    (\Psi^{\ast}_t f)_{a}(x) = \sup\limits_{y \in \R}\frac{|(\Psi_t \ast f)(x+y)|}
    {(1+|y|/t)^{a}}\quad,\quad x\in \R\,,t>0\,,
$$
originally introduced by Peetre in \cite{Pe75}.
The second ingredient is a Muckenhoupt weight $v \in \mathcal{A}_{\infty}$. The critical index $p_0$ is defined by % the critical index $p_0$ by
\begin{equation}\label{p0}
    p_0 := \inf\{p:v \in \mathcal{A}_p\}\,.
\end{equation}

\begin{proposition}\label{contchar} Let $w \in \mathcal{W}^{\alpha_3}_{\alpha_1,\alpha_2}$ and  
$v$ belong to the class $\mathcal{A}_{\infty}$, where $p_0$ is given by \eqref{p0}. Choose functions $\Phi_0,\Phi \in \mathcal{S}(\R)$
satisfying \eqref{condphi1} and \eqref{condphi2} with $R+1>\alpha_2$.
\begin{description}
 \item(i) If $0<q\leq \infty$, $0<p\leq\infty$, and $a>\frac{dp_0}{p}+\alpha_3$ then, for both $i=1,2$, 
 $$
    B^w_{p,q}(\R,v) = \{f\in \mathcal{S}'(\R)~:~\|f|B^w_{p,q}(\R,v)\|_i < \infty\}\quad,\quad i=1,2,
 $$
 where 
 \begin{eqnarray}\nonumber
     \|f|B^w_{p,q}(\R,v)\|_1 &=& \|w(x,\infty)(\Phi_0 \ast f)(x)|L_p(\R,v)\|\\ 
     &&+ \Big(\int_{0}^1
     \|w(x,t)(\Phi_t \ast f)(x)|L_p(\R,v)\|^q\frac{dt}{t}\Big)^{1/q}\,, \label{c1}\nonumber\\
     \|f|B^w_{p,q}(\R,v)\|_2 &=& \|w(x,\infty)(\Phi_0^{\ast}f)_a(x)|L_p(\R,v)\|\\
     &&+ \Big(\int_{0}^1
     \|w(x,t)(\Phi_t^{\ast}f)_a(x)|L_p(\R,v)\|^q\frac{dt}{t}\Big)^{1/q}\,.\label{c2}\nonumber
 \end{eqnarray}
Moreover, $\|\cdot|B^w_{p,q}(\R,v)\|_i$, $i=1,2$, are equivalent quasi-norms in $B^w_{p,q}(\R,v)$\,.
\item(ii) If $0<q\leq \infty$, $0<p < \infty$, and $a>d\max\{p_0/p,1/q\}+\alpha_3$ then
 $$
    F^w_{p,q}(\R,v) = \{f\in \mathcal{S}'(\R)~:~\|f|F^w_{p,q}(\R,v)\|_i < \infty\}\quad,\quad i=1,2,
 $$
 where 
 \begin{eqnarray}
     \|f|F^w_{p,q}(\R,v)\|_1 &=& \|w(x,\infty)(\Phi_0 \ast f)(x)|L_p(\R,v)\|\label{eq40}\\ 
     &&+ \Big\|\Big(\int_{0}^1
     |w(x,t)(\Phi_t \ast f)(x)|^q\frac{dt}{t}\Big)^{1/q}|L_p(\R,v)\Big\|\,, \nonumber\\
     \|f|F^w_{p,q}(\R,v)\|_2 &=& \|w(x,\infty)(\Phi_0^{\ast}f)_a(x)|L_p(\R,v)\|\label{eq41}\\
     &&+ \Big\|\Big(\int_{0}^1
     |w(x,t)(\Phi_t^{\ast}f)_a(x)|^q\frac{dt}{t}\Big)^{1/q}|L_p(\R,v)\Big\|\,.\nonumber
 \end{eqnarray}
Moreover, $\|\cdot|F^w_{p,q}(\R,v)\|_i$, $i=1,2$, are equivalent quasi-norms in $F^w_{p,q}(\R,v)$\,.
\end{description}
 
\end{proposition}
\bproof We only prove (ii) since the proof of (i) is analogous. The arguments are more or less the same as in 
the proof of \cite[Thm.\ 2.6]{T10}. 
Let us provide the necessary modifications. 

\noindent
{\em Step 1:} At the beginning of Substep 1.3 in the proof of \cite[Thm.\ 2.6]{T10} we proved a crucial inequality 
stating that for $r>0$ and $a<N$
\begin{equation}\label{2.30}
    (\Phi^{\ast}_{2^{-\ell}t}f)_{a}(x)^r\leq
    c\sum\limits_{k\in \n} 2^{-kNr}2^{(k+\ell)d}
    \int_{\R}\frac{
    |((\Phi_{k+\ell})_t\ast
    f)(y)|^r}{(1+2^{\ell}|x-y|)^{a
    r}}\,dy\,,
\end{equation}
where $c$ is independent of $f,x,t$ and $\ell$ 
but may depend on $N$ and $a$. In case $\ell = 0$ we have to replace 
$(\Phi^{\ast}_{2^{-\ell}t}f)_{a}(x)$ by $(\Phi^{\ast}_0f)_a(x)$ on the left-hand side and
$(\Phi_{k+\ell})_t$ by $\Phi_{k+\ell}=\Phi_0$ for $k=0$ on the right-hand side. We modify \eqref{2.30} by multiplying
with $|w(x,2^{-\ell}t)|^r$ on both sides ($|w(x,\infty)|^r$ in case $\ell = 0$). 
By using $w(x,2^{-\ell}t) \lesssim
2^{k\alpha_1}(1+2^{\ell}|x-y|)^{\alpha_3}w(y,2^{-(k+\ell)}t)$, which follows from $(W1)$ and $(W2)$, this gives the following modified relation
\begin{equation}\label{2.31}
    \begin{split}
      &|(\Phi^{\ast}_{2^{-\ell}t}f)_{a}(x)w(x,2^{-\ell}t)|^r\\
      &~~~\leq c\sum\limits_{k\in \n} 2^{-k(N-\alpha_1)r}2^{(k+\ell)d}
      \int_{\R}\frac{
      |((\Phi_{k+\ell})_t\ast
      f)(y)w(y,2^{-(k+\ell)}t)|^r}{(1+2^{\ell}|x-y|)^{(a-\alpha_3)
      r}}\,dy\,.
    \end{split}
\end{equation}
Now we choose $r>0$ in a way such that (a) $r(a-\alpha_3)>d$, (b)  $p/r,q/r>1$, and (c) $p_0<p/r$\,.

%\begin{description}
% \item(a) $r(a-\alpha_3)>d$,
% \item(b) $p/r,q/r>1$,
% \item(c) and $p_0<p/r$\,.
%\end{description}

Let us shortly comment on these conditions. Condition (a) is needed in order to replace the convolution integral on the
right-hand side of \eqref{2.31} by the Hardy--Littlewood maximal function $M[|w(\cdot,2^{-(k+\ell)}t)((\Phi_{k+\ell})_t\ast f)|^r](x)$ via Lemma \ref{majprop}. Conditions (b) and (c) are necessary in order to guarantee the Fefferman-Stein
maximal inequality, see Lemma \ref{muckenmax}, in the space $L_{p/r}(\ell_{q/r},\R,v)$, where we
use that $v \in \mathcal{A}_{p/r}$ as a consequence of (c) and \eqref{p0}.\\
Since $p_0\geq 1$ the conditions (a),(b),(c) are satisfied if 
$$
    \frac{d}{a-\alpha_3} < r < \min\Big\{\frac{p}{p_0},q\Big\} 
$$
which is possible if we assume $a>\alpha_3 + d\max\{p_0/p,1/q\}$\,. Now we can proceed analogously as done in Substep 1.3 of the proof
of Theorem 2.6 in \cite{T10} and obtain the equivalence of $\|f|F^w_{p,q}(\R,v)\|_1$ and $\|f|F^w_{p,q}(\R,v)\|_2$ on $\mathcal{S}'(\R)$. With the same type
of argument but some minor modifications we show that their discrete counterparts (in the sense of Definition \ref{inhom}) and \eqref{eq40}, \eqref{eq41} 
are equivalent as well. 

\noindent
{\em Step 2:} It remains to show that we can change from the system from $(\Phi_0,\Phi)$
to a system $(\Psi_0,\Psi)$ satisfying \eqref{condphi1}, \eqref{condphi2}. We argue as in Step 2 of the proof of \cite[Thm.\ 2.6]{T10}. 
There we obtain the crucial inequality 
$$
    (\Psi^{\ast}_{\ell}f)_a(x) \lesssim \sum\limits_{k=0}^{\infty} (\Phi_k^{\ast}f)_a(x)
        \left\{\begin{array}{lcl}
            2^{(\ell-k)(L+1-a)}&:&k>\ell\\
            2^{(k-\ell)(R+1)}&:&\ell \geq k
        \end{array}\right.,
$$
where $L$ can be chosen arbitrarily large. Multiplying both sides with $w_{\ell}(x)$ and using 
$$
   w_{\ell}(x) \lesssim w_k(x)\left\{\begin{array}{lcl}
            2^{(k-\ell)\alpha_1}&:&k>\ell,\\
            2^{(\ell-k)\alpha_2}&:&\ell \geq k,
        \end{array}\right.
$$
we obtain 
$$
    w_{\ell}(x)(\Psi^{\ast}_{\ell}f)_a(x) \lesssim \sum\limits_{k=0}^{\infty} w_k(x)(\Phi_k^{\ast}f)_a(x)
        \left\{\begin{array}{lcl}
            2^{(\ell-k)(L+1-a-\alpha_1)}&:&k>\ell,\\
            2^{(k-\ell)(R+1-\alpha_2)}&:&\ell \geq k,
        \end{array}\right..
$$
With our assumption $R+1>\alpha_2$ we obtain finally 
$$
    w_{\ell}(x)(\Psi^{\ast}_{\ell}f)_a(x) \lesssim \sum\limits_{k=0}^{\infty} 2^{-|k-\ell|\delta}w_k(x)(\Phi_k^{\ast}f)_a(x)\,,
$$
where $\delta = \min\{1,R+1-\alpha_2\}$\,. Now we use a straightforward generalization of the convolution 
Lemma 2 in \cite{Ry99a} and obtain immediately the desired result
$$
    \|w_{\ell}(\Psi^{\ast}_{\ell}f)_a|L_p(\ell_q,\R,v)\| \lesssim \|w_k(\Phi_k^{\ast}f)_a|L_p(\ell_q,\R,v)\|\,.
$$
This together with Step 1 and Definition \ref{inhom} concludes the proof. \eproof

\begin{Theorem}\label{mainex} Let $w\in \mathcal{W}^{\alpha_3}_{\alpha_1,\alpha_2}$, 
$v \in \mathcal{A}_{\infty}$, and $p_0$ given by \eqref{p0}. 
We choose an admissible continuous wavelet frame $\cf$ according to Definition \ref{admfr}. Let further $1\leq p\leq \infty$ ($p<\infty$ in $F$-case) and 
$1\leq q\leq \infty$. Putting 
\begin{equation}\label{eq60}
      \tilde{w}(x,t) := \left\{\begin{array}{rcl}
                                  t^{d/q-d/2}w(x,t)&:& 0<t\leq 1\,,\\
                                  w(x,\infty)&:& t=\infty\,,
                               \end{array}\right.
\end{equation}
and $B(\R) = L_p(\R,v)$ we have the following
identities in the sense of equivalent norms
$$
      B^{w}_{p,q}(\R,v) = \Co (L^{\tilde{w}}_{B,q,a},\cf)
$$
if $a>\frac{dp_0}{p}+\alpha_3$ and
$$
      F^{w}_{p,q}(\R,v) = \Co (P^{\tilde{w}}_{B,q,a},\cf)
$$
if $a>d\max\{p_0/p,1/q\}+\alpha_3$\,.
\end{Theorem}
\bproof For $t\in (0,1)$ and $x\in \R$ it holds
$$
   (V_{\cf} f)(x,t) = \left(\mathcal{D}^{L_2}_t \Phi(-\cdot) \ast \bar{f}\right)(x) = t^{d/2} \left(\mathcal{D}_t \Phi(-\cdot) \ast \bar{f}\right)(x)\,
$$
with an obvious modification in case $t=\infty$. Hence the identities 
are consequences of Definition \ref{defFS}, Theorem \ref{props}(c), and Proposition 
\ref{contchar}.\eproof\\
\newline
Now we are prepared for the discretization result, which
we state only for $F^w_{p,q}(\R,v)$. The conditions for $B^w_{p,q}(\R,v)$ are the same. 
We use the covering ${\mathcal U} = {\mathcal U}^{1,2}$ in Section \ref{seq} for $\alpha=1,\beta=2$,
and the associated sequence spaces
%Recall the definition of the sequence spaces in Paragraph \ref{seq}. We use $\beta = 2$ and $\alpha = 1$ for 
%$\mathcal{U}^{\alpha,\beta}$ and put 
$f^w_{p,q}(v) := p^{\tilde{w}}_{B,q}$ and $b^w_{p,q}(v) = \ell^{\tilde{w}}_{B,q}$ where $B(\R) = L_p(\R,v)$. We get
\begin{equation}\label{fw}
  \begin{split}  
    \|\{\lambda_{j,k}\}_{j,k}|f^{w}_{p,q}(v)\|
    &= \Big\|
    \sum\limits_{k\in \Z}w(k,\infty)|
    \lambda_{0,k}|\chi_{0,k}|L_p(\R,v)\Big\|\\
    &~~~+\Big\|\Big(\sum\limits_{j\in \n}\Big[2^{jd/2}\sum\limits_{k\in \Z}w(k 2^{-j},2^{-j})|\lambda_{j,k}|\chi_{j,k}\Big]^q\Big)^{1/q}|L_p(\R,v)\Big\|\,,
  \end{split}  
\end{equation}
as well as
\begin{equation}\label{bw}
  \begin{split}  
    \|\{\lambda_{j,k}\}_{j,k}|b^{w}_{p,q}(v)\|
    &=  \Big\|
    \sum\limits_{k\in \Z}w(k,\infty)|
    \lambda_{0,k}|\chi_{0,k}|L_p(\R,v)\Big\|\\
    &~~~+\Big(\sum\limits_{j\in \n}\Big\|2^{jd/2}
    \sum\limits_{k\in \Z}w(k 2^{-j},2^{-j})|\lambda_{j,k}|\chi_{j,k}|L_p(\R,v)\Big\|^q\Big)^{1/q}\,.
  \end{split}  
\end{equation}
Note, that Corollary \ref{corsequ} is applicable with 
$$
    \frac{d}{a-\alpha_3} < r < \min\Big\{\frac{p}{p_0},q\Big\} 
$$
as a consequence of $a>d\max\{p_0/p,1/q\}+\alpha_3$\,.

\begin{Theorem}\label{wbasesmucken} Let $w\in \mathcal{W}^{\alpha_3}_{\alpha_1,\alpha_2}$, $v \in \mathcal{A}_{\infty}$, and $p_0$ defined by \eqref{p0}.
Let further $1\leq p\leq \infty$ ($p<\infty$ in $F$-case) and $1\leq q\leq \infty$. Assume that
$\psi^0, \psi^1 \in L_2(\re)$ generate a wavelet basis in the sense of Lemma \ref{dwavelet}, 
where $\psi^0$ satisfies $(D)$, $(S_K)$, and $\psi^1$ 
satisfies $(D)$, $(S_K)$, $(M_{L-1})$ such that 
\begin{equation}\label{KLmucken}
  \begin{split}
     K,L>&\max\Big\{\Big|\max\Big\{\alpha_1,\frac{d}{2}-\frac{d}{q}\Big\}\Big|+d\max\Big\{1,\frac{2p_0}{p}\Big\}+\alpha_3,\\
     &~~~~~~~~\Big|\max\Big\{\alpha_1,\frac{d}{2}-\frac{d}{q}\Big\}\Big|+d\max\Big\{\frac{p_0}{p},\frac{3p_0}{p}-1\Big\}+2\alpha_3,\\
     &~~~~~~~~\max\Big\{\alpha_2,\frac{d}{q}-\frac{d}{2}\Big\}+d\max\Big\{\frac{p_0}{p},\frac{1}{q}\Big\}+\alpha_3,\\
     &~~~~~~~~\max\Big\{\alpha_1,\frac{d}{2}-\frac{d}{q}\Big\}+2d\max\Big\{\frac{p_0}{p},\frac{1}{q}\Big\}+2\alpha_3\Big\}\,.
  \end{split}
\end{equation}  
Then every $f\in F^w_{p,q}(\R,v)$ has the decomposition 
\begin{equation}\label{eq36}
   \begin{split}
        f  =&  \sum\limits_{c\in E}\sum\limits_{k\in \Z}\lambda^c_{0,k} \psi^c(\cdot-k)+\sum\limits_{c\in E\setminus \{0\}}\sum\limits_{j\in \N}\sum\limits_{k\in \Z}  \lambda^c_{j,k}2^{\frac{jd}{2}}\psi^c(2^j\cdot-k),
   \end{split}
\end{equation}
where the sequences $\lambda^c = \{\lambda^c_{j,k}\}_{j\in \N_0,k\in\Z}$ defined by 
$$
   \lambda^c_{j,k} = \langle f,2^{\frac{jd}{2}}\psi^c(2^j\cdot-k)\rangle\quad,\quad j\in \N_0,k\in \Z\,,
$$ 
belong to the sequence space $f^w_{p,q}(v)$
for every $c\in E$.

Conversely, an element $f\in (\mathcal{H}_{v_{\tilde{w},B,q}}^1)^{\sim}$ belongs to $F^w_{p,q}(\R,v)$ if all sequences $\lambda^c(f)$ belong to $f^w_{p,q}(v)$.
The convergence in \eqref{eq36} is in the norm of $F^w_{p,q}(\R,v)$ 
if the finite sequences are dense in 
$f^w_{p,q}(v)$. In general, we have 
weak$^{\ast}$-convergence induced by $(\mathcal{H}_{v_{\tilde{w},B,q}}^1)^{\sim}$.
\end{Theorem}
\bproof The statement is a consequence of Theorem \ref{wbasesP}, Lemma \ref{muckenh}, Theorem \ref{mainex}, and Corollary \ref{corsequ}. Indeed, the parameters $\tilde{\alpha}_1,\tilde{\alpha}_2,$ and $\tilde{\alpha}_3$ according to $\tilde{w}$ 
are given by $\tilde{\alpha}_1 = (\alpha_1+d/q-d/2)_+$, $\tilde{\alpha}_2 = (\alpha_2+d/2-d/q)_+$ and $\tilde{\alpha}_3 = \alpha_3$.
\eproof

\begin{remark} The conditions in the $B$-case are slightly weaker. Since we have then $a>dp_0/p+\alpha_3$, see Theorem \ref{mainex}, we can replace the term $d\max\{p_0/p,1/q\}$ by $dp_0/p$ in 
$\eqref{KLmucken}$.
\end{remark}

Without the weight $v$, i.e., $v\equiv 1$, we obtain wavelet characterizations for the generalized $2$-microlocal spaces
studied by Kempka in \cite{Ke09,Ke10,Ke11}.
\begin{Theorem}\label{kempka1} Let $w\in \mathcal{W}^{\alpha_3}_{\alpha_1,\alpha_2}$ and $1\leq p\leq \infty$ ($p<\infty$ in $F$-case), $1\leq q\leq \infty$. Let further $\psi^0, \psi^1 \in L_2(\re)$ generate a wavelet basis in the sense of Lemma \ref{dwavelet} where $\psi^0$ satisfies $(D)$, $(S_K)$, and $\psi^1$ satisfies $(D)$, $(S_K)$, $(M_{L-1})$ such that 
   \begin{equation}\label{KLkempka}
     \begin{split} 
        K,L>&\max\Big\{\max\Big\{\alpha_2,\frac{d}{q}-\frac{d}{2}\Big\}+d\max\Big\{\frac{1}{p},\frac{1}{q}\Big\}+\alpha_3,\\
        &~~~~~~~~\Big|\max\Big\{\alpha_1,\frac{d}{2}-\frac{d}{q}\Big\}\Big|+
        d\max\Big\{\frac{1}{p},1-\frac{1}{p}\Big\}+\alpha_3,\\
        &~~~~~~~~\max\Big\{\alpha_1,\frac{d}{2}-\frac{d}{q}\Big\}+2d\max\Big\{\frac{1}{p},\frac{1}{q}\Big\}+2\alpha_3,\\
        &~~~~~~~~\Big|\max\Big\{\alpha_1,\frac{d}{2}-\frac{d}{q}\Big\}\Big|+2\alpha_3\Big\}\,.
    \end{split}
   \end{equation}
 Then the generalized $2$-microlocal spaces $F^w_{p,q}(\R)$ can be discretized
 in the sense of Theorem \ref{wbasesmucken}.
 \end{Theorem}
 \bproof We apply Theorem \ref{wbasesP} with $B = L_p(\R)$ in connection with Theorem \ref{mainex} where in this case
 $a>d\max\{1/p,1/q\}+\alpha_3$. We can use Example \ref{Lp}
 instead of Lemma \ref{muckenh}, where $\gamma_2 = d/p$ and $\delta_2 = 0$.\eproof

\begin{remark}\label{remk} The conditions in the $B$-case are slightly weaker. Since we have then $a>d/p+\alpha_3$, see Theorem \ref{mainex} with $p_0=1$, we can replace the term $d\max\{1/p,1/q\}$ by $d/p$ in 
$\eqref{KLkempka}$.
\end{remark}

Finally, we obtain characterizations for the classical Besov-Lizorkin-Triebel spaces by putting $w(x,t) = t^{-s}$.

 \begin{Theorem}\label{classbtl} Let $1\leq p,q\leq \infty$ ($p<\infty$ in the $F$-case) and $s\in \re$. Assume that 
  $\psi^0, \psi^1 \in L_2(\re)$ generate a wavelet basis in the sense of Lemma \ref{dwavelet} and let $\psi^0$ satisfy $(D)$, $(S_K)$, and $\psi^1$ satisfy $(D)$, $(S_K)$, $(M_{L-1})$.
\begin{itemize}
\item[(i)] Assuming that 
   \begin{equation}\nonumber
      \begin{split}
        K,L>&\max\Big\{\max\Big\{s,\frac{d}{q}-\frac{d}{2}\Big\}+d\max\Big\{\frac{1}{p},\frac{1}{q}\Big\},\\
        &~~~~~~~~\Big|\min\Big\{s,\frac{d}{q}-\frac{d}{2}\Big\}\Big|+d\max\Big\{\frac{1}{p},1-\frac{1}{p}\Big\},\\
        &~~~~~~~-\min\Big\{s,\frac{d}{q}-\frac{d}{2}\Big\}+2d\max\Big\{\frac{1}{p},\frac{1}{q}\Big\}\Big\}\,.
      \end{split}
   \end{equation}
then the classical inhomogeneous Lizorkin-Triebel spaces $F^s_{p,q}(\R)$ can be discretized in the sense of 
 Theorem \ref{wbasesmucken}.
 \item[(ii)] In case 
    \begin{equation}\nonumber  
      \begin{split}  
        K,L>&\max\Big\{\max\Big\{s,\frac{d}{q}-\frac{d}{2}\Big\}+\frac{d}{p},
        \Big|\min\Big\{s,\frac{d}{q}-\frac{d}{2}\Big\}\Big|+d\max\Big\{\frac{1}{p},1-\frac{1}{p}\Big\},\\
        &~~~~~~~~-\min\Big\{s,\frac{d}{q}-\frac{d}{2}\Big\}+\frac{2d}{p}\Big\}\,.
       \end{split} 
   \end{equation}
   the classical Besov spaces $B^s_{p,q}(\R)$ can be discretized in the sense of Theorem \ref{wbasesmucken}. 
   \end{itemize}
\end{Theorem}
\bproof We apply Theorem \ref{wbasesP} with $B = L_p(\R)$ and $w(x,t) = t^{-s}$ in connection with Theorem \ref{mainex}. This gives $\tilde{\alpha}_2 = (s+d/2-d/q)_+$, $\tilde{\alpha}_1 = -(s+d/2-d/q)_-$, and $\tilde{\alpha}_3 = 0$. In the B-case we have therefore $a>d/p$, while $a>d/\min\{p,q\}$
in the $F$-case\,.
\eproof

\begin{remark}\label{splres} Theorem \ref{wbasesmucken}, Theorem \ref{kempka1}, and Theorem \ref{classbtl} provide in particular characterizations in terms of orthonormal spline wavelets, see Appendix \ref{splines}. Indeed, we have that 
$\psi^1 = \psi_m$ satisfies $(M_{L-1})$ for $L=m$ and $\psi^0 = \varphi_m, \psi^1 = \psi_m$ satisfy $(D)$ and $(S_{K})$ for $K<m-1$. 
\end{remark}

\begin{remark} Since all our results rely on the abstract Theorem \ref{wbases2} we are able to use even biorthogonal wavelets \cite{CoDaFe92}. The conditions on the smoothness and the vanishing moments of the wavelet and dual wavelet are similar. See also \cite{Ky03} for earlier results in this direction. 
\end{remark}

\subsection{Besov-Lizorkin-Triebel spaces with doubling weights}
\label{doubling}
We intend to extend the definition of weighted Besov-Lizorkin-Triebel spaces also to 
general doubling weights and give corresponding atomic and wavelet decompositions. It is well-known that, for 
general doubling weights, Lemma \ref{muckenmax} does not apply. Therefore, the first challenge
is to define certain spaces in a reasonable way, i.e., to get at least the independence of the definition of the used dyadic decomposition of unity. 

Bownik \cite{Bo07} approaches such definition by adapting the classical $\varphi$-transform 
due to Frazier and Jawerth \cite{FrJa90} to the weighted anisotropic situation. A replacement of Lemma \ref{muckenmax} is used to this end, where the classical Hardy--Littlewood maximal function is defined with respect to the doubling measure. 

Our approach is entirely different and relies on the fact that the spaces defined below can be interpreted as certain coorbits which allows to exploit our Theorems \ref{discpeetre} and \ref{wbasesP}.

A weight $v:\R \to \re_+$ is called doubling if,
$$
  \int_{B(x,2r)} v(y)\,dy \leq C\int_{B(x,r)}v(y)\,dy\quad,\quad x\in \R, r>0\,,
$$
for some positive constant $C> 1$ independent of $r$ and $x$. Note that Muckenhoupt weights in $\mathcal{A}_\infty$ are doubling, but there exist doubling
weights which are not contained in ${\mathcal A}_\infty$. For a construction of such a weight we refer
to \cite{FeMu74}. However, doubling weights are suitable in our context. We start by proving 
that the weighted Lebesgue space 
$L_p(\R,v)$ satisfies property (B2) (note that (B1) is immediate).

\begin{lemma}\label{doublemma} Let $v:\R \to (0,\infty)$ be a doubling weight with doubling constant $C\geq 1$. Then $L_p(\R,v)$ satisfies property (B2) with 
\begin{equation}\label{eq50}
  \gamma_1 = \gamma_2 = 0~,~\delta_2 = \frac{-\log_2 C}{p}~,~\delta_1 = \frac{\log_2 C}{p}\,.
\end{equation}
\end{lemma}
\bproof The idea is
   that, as a consequence of the doubling condition, $v$ cannot decay and grow too fast. On the one hand, we have
   \begin{equation}\label{f1}
      \int_{x+tQ}v(y)\,dy \leq \int_{B(0,|x|+\sqrt{d})}v(y)\,dy \lesssim (1+|x|)^{\gamma}\quad,\quad x\in \R\,,
   \end{equation}
   where $\gamma = \log_2 C$.
   On the other hand
   $$
      \int_{x+tQ}v(y)\,dy \geq  \frac{1}{C^n}\int_{B(x,|x|+\sqrt{d})}v(y)\,dy \geq \frac{1}{C^n} \int_{B(0,1)}v(y)\,dy\,,
   $$
   where $n = \lceil \log_2(c(1+|x|)) \rceil$. Hence, we get
   \begin{equation}\label{f2}
      \int_{x+Q}v(y)\,dy \gtrsim (1+|x|)^{-\gamma}\,.
   \end{equation}
   Finally, \eqref{f1} and \eqref{f2} imply that $(B2)$ is satisfied with the parameters in \eqref{eq50}\,.\eproof
   
In order to avoid Fefferman--Stein maximal inequalities (which we indeed do not have here) we modify Definition \ref{inhom} for a general doubling measure $v$ as follows.
\begin{definition}\label{inhomdoubling} Let $v$ be a doubling weight and 
$\{\varphi_j(x)\}_{j=0}^{\infty} \in \Phi(\R)$, $0<q\leq \infty$, and $a>0$. 
Let $w \in \mathcal{W}^{\alpha_3}_{\alpha_1,\alpha_2}$ and define the weight sequence $\{w_j\}_{j\in \n}$ as in \eqref{def:wj}. Put %$\Phi_j = \cf^{-1}\varphi_j$.
$\widehat{\Phi}_j = \varphi_j$.
\begin{description}
 \item(i) For $0<p\leq \infty$ we define (modification if $q=\infty$) 
 \begin{equation}\nonumber
  \begin{split}  
    B^w_{p,q,a}(\R,v) &= \Big\{f\in \mathcal{S}'(\R):\\
    &\hspace{-0.5cm}\|f|B^w_{p,q,a}(\R,v)\| = \Big(\sum\limits_{j=0}^{\infty}
    \Big\|w_j(x)\sup\limits_{z\in \R} \frac{|(\Phi_j \ast f)(x+z)|}
    {(1+2^{j}|z|)^a}|L_p(\R,v)\Big\|^q\Big)^{1/q}<\infty\Big\}\,.
  \end{split}  
 \end{equation}
 \item(ii) For $0<p<\infty$ we define (modification if $q=\infty$)
 \begin{equation}\nonumber
  \begin{split}  
    F^w_{p,q,a}(\R,v) &= \Big\{f\in \mathcal{S}'(\R):\\
    &\hspace{-0.5cm}\|f|F^w_{p,q,a}(\R,v)\| = \Big\|\Big(\sum\limits_{j=0}^{\infty}\Big|
    w_j(x)\sup\limits_{z\in \R} \frac{|(\Phi_j \ast f)(x+z)|}
    {(1+2^{j}|z|)^a}\Big|^q\Big)^{1/q}|L_p(\R,v)\Big\|<\infty\Big\}\,.
  \end{split}
 \end{equation} 
\end{description}
\end{definition}
Here we have a counterpart of Proposition \ref{contchar} stating 
that $\|\cdot|F^w_{p,q,a}(\R,v)\|_2$ and\linebreak $\|\cdot|B^w_{p,q,a}{\R,v}\|_2$ are equivalent
characterizations (for all $a>0$) for the $F$ and $B$-spaces, respectively. To show this, we switch in a first step 
from one system $\Phi$ to another system $\Psi$ in the discrete characterization given in Definition \ref{inhomdoubling}. Indeed, we argue analogously as in Step 2 of the proof of Proposition \ref{contchar}, see also \cite[Thm.\ 2.6]{T10}. Note, that we did not use a Fefferman--Stein maximal inequality there. With a similar argument we switch in a second step from the discrete characterization to the continuous characterization (Prop. \ref{contchar}) using the same sytem $\Phi$. 

Consequently, we identify Besov-Lizorkin-Triebel spaces with doubling weights as coorbits.
%obtain the following result
\begin{Theorem}\label{coo:doubling} Let $w \in \mathcal{W}^{\alpha_3}_{\alpha_1,\alpha_2}$ and $v$ be a general doubling weight with doubling constant $C>1$. Choose $\cf$ to be an admissible continuous wavelet frame according to Definition \ref{admfr}. 
Let further $1\leq p,q\leq \infty$ and $a>0$.
Putting $\tilde{w}(x,t)$ as in \eqref{eq60} and $B(\R) = L_p(\R,v)$ we have the following
identities in the sense of equivalent norms
\begin{equation}\nonumber
  \begin{split}    
      B^{w}_{p,q,a}(\R,v) &= \Co (L^{\tilde{w}}_{B,q,a},\cf)\,,\\
      F^{w}_{p,q,a}(\R,v) &= \Co (P^{\tilde{w}}_{B,q,a},\cf).
  \end{split}    
\end{equation}
\end{Theorem}
Based on Theorem \ref{wbasesP} we immediately arrive at one of our main discretization results. We state it only for the $F$-spaces. The conditions for the $B$-spaces are the same. 

\begin{Theorem}\label{wbasesdoubling} Let $w \in \mathcal{W}^{\alpha_3}_{\alpha_1,\alpha_2}$ and let $v$ be a doubling weight with
doubling constant $C>1$. Assume $1\leq p,q\leq \infty$, $a>0$, and let 
$\psi^0, \psi^1 \in L_2(\re)$ generate an orthonormal wavelet basis in the sense of Lemma \ref{dwavelet} where
$\psi^0$ satisfies $(D)$, $(S_K)$, and $\psi^1$ satisfies $(D)$, $(S_K)$, $(M_{L-1})$ such that 
 \begin{equation}\nonumber
  \begin{split}
     K,L>&\max\Big\{\Big|\max\Big\{\alpha_1,\frac{d}{2}-\frac{d}{q}\Big\}\Big|+d+\frac{\log_2 C}{p}+\alpha_3,\\
     &~~~~~~~~\Big|\max\Big\{\alpha_1,\frac{d}{2}-\frac{d}{q}\Big\}\Big|+2\Big(\frac{\log_2 C}{p}+\alpha_3\Big),\\
     &~~~~~~~~\max\Big\{\alpha_2,\frac{d}{q}-\frac{d}{2}\Big\}+a,\max\Big\{\alpha_1,\frac{d}{2}-\frac{d}{q}\Big\}+2a
     \Big\}\,.
  \end{split}
\end{equation}  
Then every $f\in F^w_{p,q,a}(\R,v)$ has the decomposition \eqref{eq36},
where the sequences\linebreak $\lambda^c = \{\lambda^c_{j,k}\}_{j\in \N_0,k\in\Z}$, $c \in E$,
are contained in $(P^{\tilde{w}}_{L_p(\R,v),q,a})^\sharp({\mathcal U}^{1,2})$. The latter is equivalent to \eqref{sequ_1} with $B(\R)$ replaced by $L_p(\R,v)$, $\beta^{dj/q}$ by $\beta^{dj/2}$, and with $\beta = 2$ and $\alpha = 1$.

Conversely, an element $f\in (\mathcal{H}^1_{v_{\tilde{w},B,q}})^{\sim}$ belongs to $F^w_{p,q,a}(\R,v)$ if all sequences $\lambda^c(f)$, $c \in E$, belong to $(P^{\tilde{w}}_{L_p(\R,v),q,a})^\sharp({\mathcal U}^{1,2})$. The convergence in \eqref{eq36} is in the norm of $F^w_{p,q,a}(\R,v)$ if the finite sequences are dense in 
$(P^{\tilde{w}}_{L_p(\R,v),q,a})^\sharp({\mathcal U}^{1,2})$. In general, we have weak$^{\ast}$-convergence. 
\end{Theorem}

\subsection{Generalized 2-microlocal Besov-Lizorkin-Triebel-Morrey spaces} 
\label{Morrey}
Several applications in PDEs require the investigation of smoothness spaces constructed on Morrey spaces \cite{Mo38}.
The spaces of Besov-Lizorkin-Triebel-Morrey type are currently a very active research area. We refer to 
Sawano \cite{Sa08}, Sawano et al.\ \cite{SaTa07, SaYaYu10}, Tang--Xu \cite{TaXu05} as well as to the recent monograph by Yuan et al.\ \cite{YuSiYa10} and the references given there. 
Our intention in the current paragraph is to provide wavelet decomposition theorems as consequences of the fact 
that the mentioned spaces can be interpreted as coorbits. Note that \cite{Sa08,SaTa07} and \cite{YuSiYa10} have already dealt 
with atomic and wavelet decompositions of these spaces. Our results have to be compared with the ones in there, see the list at the end of Subsection \ref{overBe} above. 

We start with the definition of the Morrey space $M_{q,p}(\R)$ on $\R$.
\begin{definition}\label{def:morrey} Let $0<p\leq q\leq \infty$. Then the Morrey space $M_{q,p}(\R)$ is defined as the
collection of all measurable and locally Lebesgue-integrable functions $f$ with finite (quasi-)norm
\begin{equation}\label{eq51}
    \|f|M_{q,p}(\R)\| = \sup\limits_{R>0,\,x\in \R} R^{d(1/q-1/p)}\Big(\int_{B_R(x)}|f(y)|^p\,dy\Big)^{1/p}
\end{equation}
if $p<\infty$, where $B_R(x)$ denotes the Euclidean ball with radius $R>0$ and center $x\in \R$. In the case $p=\infty$
we put $M_{\infty,\infty}(\R) := L_{\infty}(\R)$.
\end{definition}

These spaces -- studied first by Morrey \cite{Mo38} --  generalize the ordinary Lebesgue spaces. Indeed, we have 
$M_{p,p}(\R) = L_p(\R)$, $0<p\leq \infty$. %They were first studied by Morrey \cite{Mo38}. 
In the case $q<p$ the quantity \eqref{eq51} is infinite as soon as $f \neq 0$, so that $M_{q,p}(\R) = \{0\}$. If $p\geq 1$ 
then $M_{q,p}(\R)$ is a Banach space, otherwise a quasi-Banach space. The following Lemma ensures that $M_{q,p}(\R)$ satisfies (B2). 

\begin{lemma}\label{morlemma} Let $0<p\leq q \leq \infty$. Then the space $M_{q,p}(\R)$ satisfies (B2) with $\gamma_1 = \gamma_2 = d/q$
and $\delta_1 = \delta_2 = 0$.
\end{lemma}
\bproof Fix $x\in \R$ and $0<t<1$. We consider the norm of the characteristic function $\chi^{\alpha}_{(x,t)}$ of the cube 
$Q^{\alpha}_{(x,t)} = x+t[-\alpha,\alpha]^d$ in $M_{q,p}(\R)$\,. By \eqref{eq51} we obtain immediately
  $$
        \|\chi^{\alpha}_{(x,t)}|M_{q,p}(\R)\| \gtrsim t^{d(1/q-1/p)}t^{d/p} = t^{d/q}\,.
  $$
For the reverse estimate we use the well-known fact that $L_q(\R) \hookrightarrow M_{q,p}(\R)$, see \cite{Mo38}. 
Therefore, we have
$$
    \|\chi^{\alpha}_{(x,t)}|M_{q,p}(\R)\| \lesssim \|\chi^{\alpha}_{(x,t)}|L_q(\R)\| \lesssim t^{d/q}\,
$$ 
which concludes the proof. \eproof

We define the 2-microlocal Besov-Morrey spaces $B^{w,u}_{p,q}(\R)$ and 
Lizorkin-Triebel-Morrey spaces $F^{w,u}_{p,q}(\R)$ by replacing 
$L_p(\R,v)$ by $M_{u,p}(\R)$, $u>p$, in Definition \ref{inhom}. Here $w \in \mathcal{W}^{\alpha_3}_{\alpha_1,\alpha_2}$ is a weight 
function and $0<p,q\leq \infty$, $0<p<u\leq \infty$, where $p<\infty$ 
in the $F$-case. This is a straightforward generalization of the 
spaces appearing in \cite{TaXu05,Sa08, SaTa07, SaYaYu10, YuSiYa10}. With exactly the same proof 
techniques we obtain a counterpart of Proposition \ref{contchar} under the conditions $a>d/p+\alpha_3$ in the $B$-case 
and $a>d\max\{1/p,1/q\}+\alpha_3$ in the $F$-case. One uses a vector-valued 
Fefferman--Stein type maximal inequality for the space $M_{u,p}(\ell_{q},\R)$, where $1<p\leq u<\infty$ and $1 < q \leq \infty$, see
\cite{TaXu05} and \cite{ChFr87} for the case $q = \infty$.

As a consequence, the Besov-Lizorkin-Triebel-Morrey spaces can be identified as coorbits, i.e., the following counterpart to Theorem \ref{mainex} holds.

\begin{Theorem}\label{coo:morrey} Under the assumptions of Theorem \ref{mainex}
we have in the sense of equivalent norms
$$
      B^{w,u}_{p,q}(\R) = \Co (L^{\tilde{w}}_{M_{u,p},q,a},\cf)
$$
if $1\leq p\leq u\leq \infty$, $1\leq q\leq \infty$, and $a>\frac{d}{p}+\alpha_3$
as well as
$$
      F^{w,u}_{p,q}(\R) = \Co (P^{\tilde{w}}_{M_{u,p},q,a},\cf)
$$
for $1\leq p\leq u < \infty$, $1\leq q\leq \infty$, and $a>d\max\{1/p,1/q\}+\alpha_3$.
\end{Theorem}
Since Corollary \ref{corsequ} is applies for the space $M_{u,p}(\ell_q,\R)$ with
$$
    \frac{d}{a-\alpha_3} < r < \min\{p,q\} 
$$
we may (and do) define the sequence spaces $f^{w,u}_{p,q}$ and $b^{w,u}_{p,q}$ just by replacing 
$L_p(\R,v)$ by $M_{u,p}(\R)$ in \eqref{fw} and \eqref{bw}.

The main result of this subsection is the following

\begin{Theorem}\label{morreydec} Let $w\in \mathcal{W}^{\alpha_3}_{\alpha_1,\alpha_2}$ 
and $1\leq p\leq u< \infty$, $1\leq q\leq \infty$. Assume that $\psi^0, \psi^1 \in L_2(\re)$ generate an orthonormal wavelet basis in the sense of Lemma \ref{dwavelet} and let $\psi^0$ satisfy $(D)$, $(S_K)$, and $\psi^1$ satisfy $(D)$, $(S_K)$, $(M_{L-1})$ such that 
\begin{equation}\nonumber
  \begin{split}  
        K,L>&\max\Big\{\max\Big\{\alpha_2,\frac{d}{q}-\frac{d}{2}\Big\}+d\max\Big\{\frac{1}{p},\frac{1}{q}\Big\}+\alpha_3,\\
        &~~~~~~~~\Big|\max\Big\{\alpha_1,\frac{d}{2}-\frac{d}{q}\Big\}\Big|+
        d\max\Big\{\frac{1}{u},1-\frac{1}{u}\Big\}+\alpha_3,\\
        &~~~~~~~~\max\Big\{\alpha_1,\frac{d}{2}-\frac{d}{q}\Big\}+2d\max\Big\{\frac{1}{p},\frac{1}{q}\Big\}+2\alpha_3,\\
        &~~~~~~~~\Big|\max\Big\{\alpha_1,\frac{d}{2}-\frac{d}{q}\Big\}\Big|+2\alpha_3\Big\}\,.
    \end{split}
\end{equation}
Then every $f\in F^{w,u}_{p,q}(\R)$ has the decomposition \eqref{eq36}
where the sequences $\lambda^c = \{\lambda^c_{j,k}\}_{j\in \N_0,k\in\Z}$ 
belong to the sequence space $f^{w,u}_{p,q}$
for every $c\in E$.

Conversely, 
an element $f\in (\mathcal{H}_{v_{\tilde{w},B,q}}^1)^{\sim}$ belongs to $F^{w,u}_{p,q}(\R)$ if all sequences $\lambda^c(f)$ belong to $f^{w,u}_{p,q}$.
The convergence in \eqref{eq36} is considered in $F^{w,u}_{p,q}(\R)$ if the finite sequences are dense in 
$f^{w,u}_{p,q}(\R)$. In general, we have 
weak$^{\ast}$-convergence. 
\end{Theorem}

\begin{remark} The modifications for the $B$-spaces are according to Remark \ref{remk}.
\end{remark}

\subsection{Spaces of dominating mixed smoothness}
\label{dommix}
Recently, there has been an increasing interest in function spaces of dominating mixed smoothness, see 
\cite{Ba03,ScTr87,SiUl09, SiUl10,TDiff06, T08_2,  Vy06}.
Their structure is suitable for treating high-dimensional approximation and integration problems efficiently
and overcome the so-called curse of dimensionality to some extent. These spaces can as well be treated in terms of our generalized coorbit space theory.
We briefly describe this setting. In a certain sense these spaces 
behave like the isotropic ones for $d=1$, and consequently, the proofs operate by iterating the techniques 
from Subsections \ref{PeetSp}, \ref{seq}, and \ref{coorbits}. 

We start with a definition of mixed Peetre spaces on $\bar{X} = X \times \cdots \times X$, where $X = \mathbb{R} \times [(0,1] \cup \{\infty\}]$, as a tensorized version of Definition \ref{defFS}. Our definition is
motivated by equivalent characterizations of dominating mixed spaces, which are obtained by a 
combination of the techniques in \cite{T10} with \cite{T08_2, Vy06}.

\begin{definition}\label{defFSmixed} Let $1\leq p, q\leq \infty$ and $\bar{a}>1$. Let further $\bar{r} \in \R$. We define by 
\begin{equation}\nonumber
  \begin{split}
     P^{\bar{r}}_{p,q,\bar{a}}(\bar{X}) &= \{F:\bar{X} \to \C~:~\|F|P^{w}_{p,q,\bar {a}}\| < \infty\}\,,\\
     L^{\bar{r}}_{p,q,\bar{a}}(\bar{X}) &= \{F:\bar{X} \to \C~:~\|F|L^{\bar{r}}_{p,q,\bar{a}}\| < \infty\}
  \end{split}   
\end{equation}
two scales of Banach function spaces on $\bar{X}$, where we put
\begin{equation}\nonumber
  \begin{split} 
    &\|F|P^{\bar{r}}_{p,q,\bar{a}}\| := \Big\|\sup\limits_{z\in \R}\frac{|F((x_1+z_1,\infty),\hdots,(x_d+z_d,\infty))|}{(1+|z_1|)^{a_1}\cdots (1+|z_d|)^{a_d}}|L_p(\R)\Big\|\\
    &+ \sum\limits_{\substack{A \subset \{1,\hdots,d\}\\A \neq \emptyset}} \Big\|\Big(\int_{0}^1 \cdots
    \int_{0}^1 \Big[\sup\limits_{z\in \R}\frac{|F((x_1+z_1,t_1),\hdots,(x_d+z_d,t_d))|}{\prod\limits_{i\in A}(1+|z_i|/t_i)^{a_i}\prod\limits_{i\notin A}(1+|z_i|)^{a_i}}
   \Big]^q\prod\limits_{i\in A}t_{_i}^{-r_iq}\frac{dt_{i}}{t^{2}_{i}}\Big)^{1/q}|L_p(\R)\Big\|
  \end{split} 
\end{equation}\nonumber
and 
\begin{equation}
  \begin{split} 
    &\|F|L^{\bar{r}}_{p,q,\bar{a}}\|:=\Big\|\sup\limits_{z\in \R}\frac{|F((x_1+z_1,\infty),\hdots,(x_d+z_d,\infty))|}{(1+|z_1|)^{a_1}\cdots (1+|z_d|)^{a_d}}|L_p(\R)\Big\|\\
    &+ \sum\limits_{\substack{A \subset \{1,\hdots,d\}\\A \neq \emptyset}} \Big(\int_{0}^1 \cdots
    \int_{0}^1 \Big\|\sup\limits_{z\in \R}\frac{|F((x_1+z_1,t_1),\hdots,(x_d+z_d,t_d))|}{\prod\limits_{i\in A}(1+|z_i|/t_i)^{a_i}\prod\limits_{i\notin A}(1+|z_i|)^{a_i}}
   |L_p(\R)\Big\|^q\prod\limits_{i\in A}t_i^{-r_{i}q}\frac{dt_{i}}{t^2_{i}}\Big)^{1/q}\,.
  \end{split} 
\end{equation}
For fixed $A \subset \{1,\hdots,d\}$ we put $t_i = \infty$ if $i\notin A$. In case $q=\infty$ the integrals over $t_i, i\in A$, 
have to be replaced by 
a supremum over $t_i$.
\end{definition}

\subsubsection*{Associated sequence spaces}
We cover the space $\bar{X}$ by the Cartesian product of the family from Subsection \ref{seq}.
For fixed $\alpha>0$ and $\beta>1$ we consider
the family $\overline{\mathcal{U}}^{\alpha,\beta} = \{\bar{U}_{\bar{j},\bar{k}}\}_{j\in \n^d, k\in \zz^d}$ of subsets
$$
   \bar{U}_{\bar{j},\bar{k}} = U_{j_1,k_1} \times \cdots \times U_{j_d,k_d} \,.
$$
Clearly, we have $\bar{X} \subset \bigcup\limits_{\bar{j}\in \n^d, \bar{k}\in \Z} \bar{U}_{\bar{j},\bar{k}}$.
We use the notation
$$
    \chi_{\bar{j},\bar{k}}(x) = (\chi_{j_1,k_1}\otimes \cdots \otimes \chi_{j_d,k_d})(x) = \prod_{i=1}^d \chi_{j_i,k_i}(x_i) \quad,\quad x\in \R\,.
$$
Iterating dimensionwise the arguments leading to \eqref{eq-0} 
gives the following description for the sequence spaces $(P^{\bar{r}}_{p,q,\bar{a}})^{\sharp}$
and $(L^{\bar{r}}_{p,q,\bar{a}})^{\sharp}$\,.

\begin{Theorem} If $1\leq p,q \leq \infty$, $\bar{a}>1/\min\{p,q\}$, and $\bar{r}\in \R$ then 
\begin{equation}\label{mixedsequ1}
   \|\{\lambda_{\bar{j},\bar{k}}\}_{\bar{j},\bar{k}}|(P^{\bar{r}}_{p,q,\bar{a}})^{\sharp}\|
    \asymp \Big\|\Big(\sum\limits_{\bar{j}\in \n^d}\Big[\beta^{|j|_1/q}\sum\limits_{\bar{k}\in \Z}
    \beta^{\bar{j}\bar{r}}|\lambda_{\bar{j},\bar{k}}|\chi_{\bar{j},\bar{k}}(x)\Big]^q\Big)^{1/q}|L_p(\R)\Big\|
 \end{equation}
and
\begin{equation}\label{mixedsequ2}
   \|\{\lambda_{\bar{j},\bar{k}}\}_{\bar{j},\bar{k}}|(L^{\bar{r}}_{p,q,\bar{a}})^{\sharp}\| \asymp \Big(\sum\limits_{\bar{j}\in \n^d}\Big(\prod\limits_{i=1}^d\beta^{j_i(r_i+1/q-1/p)q}\Big)\Big(\sum\limits_{\bar{k}\in \Z}
   |\lambda_{\bar{j},\bar{k}}|^p\Big)^{q/p}\Big)^{1/q}\,.
\end{equation}
We have $(L^{\bar{r}}_{p,q,\bar{a}})^{\sharp} = (L^{\bar{r}}_{p,q,\bar{a}})^{\flat}$ and 
$(P^{\bar{r}}_{p,q,\bar{a}})^{\sharp} = (P^{\bar{r}}_{p,q,\bar{a}})^{\flat}$, respectively.
\end{Theorem}

\bproof Since we deal here with usual $L_p(\R)$ and $L_p(\ell_q)$-spaces we can use the methods from Corollary \ref{corsequ} 
to obtain \eqref{mixedsequ1} and \eqref{mixedsequ2}. The Hardy--Littlewood maximal operator then acts componentwise. 
For the corresponding maximal inequality  see \cite[Thm.\ 1.11]{Vy06}.
\eproof

\subsubsection*{The coorbits of $L^{\bar{r}}_{p,q,\bar{a}}(\bar{X})$ and $P^{\bar{r}}_{p,q,\bar{a}}(\bar{X})$}
\label{coorbittensor}
We apply the abstract theory in a situation where the index set is given by 
$$\bar{X} = \underbrace{X \times \cdots \times X}_{d-\mbox{times}}$$ 
with 
$X = \re\times [(0,1)\cup \{\infty\}]$. This space is equipped with the product measure 
$\mu_{\bar{X}} = \mu_X \otimes \cdots \otimes \mu_X$, i.e., 
$$
    \int_{\bar{X}} F(\z_1,\cdots,\z_d)\,\mu_{\bar{X}}(d\z) = \int_X\cdots \int_X F(\z_1,\hdots,\z_d) \mu_X(d\z_1)\cdots \mu_X(d\z_d)\,.
$$
We put $\mathcal{H} = L_2(\R)$. We choose an admissible continuous frame $\cf_1 = \{\varphi_\x\}_{\x\in X}$ according
to Definition \ref{admfr}. For $\z = (\z_1,...,\z_d) \in \bar{X}$ we define 
$\bar{\varphi}_{\z} := \varphi_{\z_1}\otimes \cdots \otimes \varphi_{\z_d}$. It is  easy to see
that the system $\bar{\mathcal{F}} = \{\bar{\varphi}_{\z}\}_{\z\in \bar{X}}$ represents a tight continuous frame indexed by $\bar{X}$ in $\mathcal{H}$. The corresponding frame transform is given by
$$
      V_{\bar{\mathcal{F}}}f(\z) = \langle f, \bar{\varphi}_{\z} \rangle\quad,\quad \z\in \bar{X}\,.
$$
For $1\leq p,q \leq \infty$ and $\bar{r}\in \R$ put for $i=1,\hdots,d$
$$
v_{p,q,r_i}(x,t) = \left\{\begin{array}{rcl}
              1&:& t=\infty\\
              \max\{t^{-(1/q-1/p)}t^{-r_i},t^{-(1/p-1/q)}t^{r_i}\}&:& t\in (0,1)\,
       \end{array}\right.
$$
and $v_{p,q,\bar{r}} = v_{p,q,r_i}\otimes \cdots \otimes v_{p,q,r_i}$. Let us define the corresponding coorbit spaces. 

\begin{definition}\label{deftensor} Let $1\leq p,q \leq \infty$, $\bar{r} \in \R$, and $\bar{a}>0$. We define
\begin{equation}\nonumber
    \begin{split}  
    \Co P^{\bar{r}}_{p,q,\bar{a}} &= \Co (P^{\bar{r}}_{p,q,\bar{a}},\bar{\cf}):= \{f\in (\mathcal{H}^1_{v_{p,q,\bar{r}}})^{\sim}~:~V_{\mathcal{\bar{\cf}}}f \in P^{\bar{r}}_{p,q,\bar{r}}(\bar{X})\}\,,\\
    \Co L^{\bar{r}}_{p,q,\bar{a}} &= \Co (L^{\bar{r}}_{p,q,\bar{a}},\bar{\cf}):=\{f \in (\mathcal{H}^1_{v_{p,q,\bar{a}}})^{\sim}~:~
    V_{\mathcal{\bar{\cf}}}f \in L^{\bar{r}}_{p,q,\bar{a}}(\bar{X})\}\,.
  \end{split}  
\end{equation} 
\end{definition}
An iteration of the techniques from Section \ref{Peetre} shows that all the conditions needed for the above definition are valid.

\begin{remark}\label{genmixed} It is also possible to define spaces with dominating mixed smoothness in a more general sense as done 
in the isotropic case, see Definition \ref{defFS}
and the corresponding coorbit spaces. Indeed, it is possible to treat even weighted spaces or 2-microlocal spaces of dominating
mixed smoothness as in the previous subsections.
\end{remark}

\begin{Theorem}\label{prop}  The space $\Co P^{\bar{r}}_{p,q,\bar{a}}$ and $\Co L^{\bar{r}}_{p,q,\bar{a}}$ are Banach spaces and 
 do not depend on the frame $\bar{\cf}$. Furthermore, we have the identity
 $$
    \Co P^{\bar{r}}_{p,q,\bar{a}} = \{f\in \mathcal{S}'(\R)~:~V_{\bar{\mathcal{F}}}f \in \Co P^{\bar{r}}_{p,q,\bar{a}}\}\,,
 $$
 respective for $\Co L^{\bar{r}}_{p,q,\bar{a}}$\,.
\end{Theorem}

\subsubsection*{Relation to classical spaces}
We give the definition of the spaces $S^{\bar{r}}_{p,q}F(\R)$ and $S^{\bar{r}}_{p,q}B(\R)$.  
It is well-known that these spaces can be characterized in a discrete way via so-called local means and Peetre maximal functions
\cite{Vy06,T08_2,HaDiss10,T10}.
%We refer to \cite{Vy06} in connection with
%\cite{T08_2}, see also \cite{HaDiss10} and \cite{T10}. 
Recall the notion of decomposition of unity in Definition \ref{decunity}. We start with 
$d$ systems $\varphi^i\in \Phi(\re)$ for $i=1,\hdots,d$ and 
put
$$
    (\varphi^1 \otimes \cdots \otimes \varphi^d)_{\bar{\ell}}(\xi_1,\hdots,\xi_d) :=
    \varphi^1_{\ell_1}(\xi_1) \cdots  \varphi^d_{\ell_d}(\xi_d)\quad,\quad \xi\in \R,\bar{\ell}\in \n^d. 
$$
\begin{definition}\label{sprd} \rm Let $\bar{r} = (r_1,\hdots,r_d) \in \R$ and $0<q\leq \infty$.
  \begin{description}
      \item(i) Let $0<p\leq\infty$. Then $S^{\bar{r}}_{p,q}B(\dopp{R}^d)$ is the collection of all
      $f\in \mathcal{S}'(\dopp{R}^d)$ such that
      $$
        \|f|S^{\bar{r}}_{p,q}B(\dopp{R}^d)\|^{\bar{\varphi}} =
        \bigg(\sum\limits_{\bar{\ell}\in\dopp{N}_0^d}2^{\bar{r}\cdot\bar{\ell}q}
        \big\|[(\varphi^1\otimes \cdots \otimes\varphi^d)_{\bar{\ell}}(\xi)\widehat{f}]^{\vee}(x)|L_p(\R)
        \big\|^q\bigg)^{1/q}
      $$
      is finite (modification if $q=\infty$).
      \item(ii) Let $0<p<\infty$. Then $S^{\bar{r}}_{p,q}F(\dopp{R}^d)$ is the collection of all
      $f\in \mathcal{S}'(\dopp{R}^d)$ such that
      $$
        \|f|S^{\bar{r}}_{p,q}F(\dopp{R}^d)\|^{\bar{\varphi}} =
        \bigg\|\bigg(\sum\limits_{\bar{\ell}\in\dopp{N}_0^d}2^{\bar{r}\cdot\bar{\ell}q}
        |[(\varphi^1\otimes \cdots \otimes\varphi^d)_{\bar{\ell}}(\xi)\widehat{f}]^{\vee}(x)|^q\bigg)^{1/q}
        \bigg|L_p(\dopp{R}^d)\bigg\|
      $$
      is finite (modification if $q=\infty$).
  \end{description}
\end{definition}

The following theorem states the relation between previously defined coorbit spaces and the classical 
spaces with dominating mixed smoothness. 

\begin{Theorem}\label{coo:dom} Let $1\leq p,q \leq \infty$, $(p<\infty$ in the $F$-case), $\bar{r} \in \R$, and $\bar{a}>1/\min\{p,q\}$. Then 
we have in the sense of equivalent norms 
$$
    S^{\bar{r}}_{p,q}F(\R) = \Co (P^{\bar{r}+1/2-1/q}_{p,q,\bar{a}},\bar{\cf})
$$ 
and if $\bar{a}>1/p$
$$
   S^{\bar{r}}_{p,q}B(\R) = \Co (L^{\bar{r}+1/2-1/q}_{p,q,\bar{a}},\bar{\cf})\,.
$$
\end{Theorem}
\bproof We apply the continuous characterization in terms of Peetre maximal functions
of local means on the left-hand side, see \cite{T10} in connection with \cite{T08_2}. 
Then we apply Theorem \ref{prop} and get the equivalence. \eproof

 It is also possible to obtain a 
``semi-discrete'' characterization (in the sense of Definition \ref{sprd}) for the spaces on the right-hand side by using the abstract 
coorbit space theory from Section \ref{abstrth}.

\subsubsection*{Wavelet bases}
Below we state a multivariate version of Theorem \ref{classbtl} on wavelet basis characterizations using tensor product wavelet frames. Let us start with a scaling function $\psi^0$ and a corresponding wavelet $\psi^1 \in L_2(\re)$ satisfying $(D)$, $(M_{L-1})$, and $(S_K)$ for some $K$ and $L$. In the sequel we use the tensor product wavelet system $\{\psi_{\bar{j},\bar{k}}\}_{\bar{j},\bar{k}}$
defined in Appendix \ref{appdwavelet}.

We are interested in sufficient conditions on $K,L$ such that every $f\in S^{\bar{r}}_{p,q}F(\R)$ or
$S^{\bar{r}}_{p,q}B(\R)$, respectively, has the decomposition 
\begin{equation}\label{exp}
   \begin{split}
        f  =&  \sum\limits_{\bar{j}\in \n^d}\sum\limits_{\bar{k}\in \Z} \lambda_{\bar{j},\bar{k}}\psi_{\bar{j},\bar{k}}
   \end{split}
\end{equation}
and the sequence $\lambda = \lambda(f) = \{\lambda_{\bar{j},\bar{k}}\}_{\bar{j}\in \n^d,\bar{k}\in\Z}$ defined by 
$
   \lambda_{\bar{j},\bar{k}} = \langle f,\psi_{\bar{j},\bar{k}}\rangle,\bar{j}\in \n^d,\bar{k}\in \Z\,,
$
belongs to the sequence spaces
\begin{eqnarray}\nonumber
   \|\{\lambda_{\bar{j},\bar{k}}\}_{\bar{j},\bar{k}}|s^{\bar{r}}_{p,q}f\|
    &=& \Big\|\Big(\sum\limits_{\bar{j}\in \n^d}\Big[2^{|j|_1/2}\sum\limits_{\bar{k}\in \Z}
    2^{\bar{j}\bar{r}}|\lambda_{\bar{j},\bar{k}}|\chi_{\bar{j},\bar{k}}(x)\Big]^q\Big)^{1/q}|L_p(\R)\Big\|\\
   \|\{\lambda_{\bar{j},\bar{k}}\}_{\bar{j},\bar{k}}|s^{\bar{r}}_{p,q}b\| 
    &=&\Big(\sum\limits_{\bar{j}\in \n^d}\Big(\prod\limits_{i=1}^d 2^{j_i(r_i+1/2-1/p)q}\Big)\Big(\sum\limits_{\bar{k}\in \Z}
   |\lambda_{\bar{j},\bar{k}}|^p\Big)^{q/p}\Big)^{1/q}\nonumber
\end{eqnarray}
corresponding to $(P^{\bar{r}+1/2-1/q}_{p,q,\bar{a}})^{\sharp}$, see \eqref{mixedsequ1}, and $(L^{\bar{r}+1/2-1/q}_{p,q,\bar{a}})^{\sharp}$, see \eqref{mixedsequ2}, where we used $\overline{\mathcal{U}}^{1,2}$. 
Furthermore, we aim at the converse that an element $f\in (\mathcal{H}^1_{v})^{\sim}$ belongs to $S^{\bar{r}}_{p,q}F(\R)$ or $S^{\bar{r}}_{p,q}B(\R)$ if the sequence $\lambda(f)$ belongs 
to $s^{\bar{r}}_{p,q}f$ or $s^{\bar{r}}_{p,q}b$, respectively. 
The convergence in \eqref{exp} is required to be in the norm of $S^{\bar{r}}_{p,q}F(\R)$ or $S^{\bar{r}}_{p,q}B(\R)$
if the finite sequences are dense in $s^{\bar{r}}_{p,q}f$ or $s^{\bar{r}}_{p,q}b$, respectively. In general we require weak$^{\ast}$-convergence. 

The following theorem provides the corresponding wavelet basis characterization of spaces of mixed dominating smoothness.

\begin{Theorem}\label{wbasestensor} Let $1\leq p,q\leq \infty$ ($p<\infty$ in the $F$-case) and $\bar{r} \in \R$. Let further 
  $\psi^0, \psi^1 \in L_2(\re)$ be a scaling function and associated wavelet where $\psi^0$ satisfies $(D)$, $(S_K)$, and $\psi^1$ 
  satisfies $(D)$, $(S_K)$, $(M_{L-1})$.
  \begin{itemize}
  \item[(i)] If, for $i=1,\hdots,d$, 
   \begin{equation}\label{eq-71}
     \begin{split} 
        K,L>&\max\Big\{\max\Big\{r_i,\frac{1}{q}-\frac{1}{2}\Big\}+\max\Big\{\frac{1}{p},\frac{1}{q}\Big\},\\
        &~~~~~~~~\Big|\min\Big\{r_i,\frac{1}{q}-\frac{1}{2}\Big\}\Big|+\max\Big\{\frac{1}{p},1-\frac{1}{p}\Big\},\\
        &~~~~~~~-\min\Big\{r_i,\frac{1}{q}-\frac{1}{2}\Big\}+2\max\Big\{\frac{1}{p},\frac{1}{q}\Big\}\Big\}
     \end{split}
   \end{equation}
 then the inhomogeneous Lizorkin-Triebel space with dominating mixed smoothness $S^{\bar{r}}_{p,q}F(\R)$ can be discretized in the sense of \eqref{exp} using the sequence space $s^{\bar{r}}_{p,q}f$.
 \item[(ii)] If, for $i=1,\hdots,d$,
    \begin{equation}\label{eq-72}
      \begin{split} 
        K,L>&\max\Big\{\max\Big\{r_i,\frac{1}{q}-\frac{1}{2}\Big\}+\frac{1}{p},
        \Big|\min\Big\{r_i,\frac{1}{q}-\frac{1}{2}\Big\}\Big|+\max\Big\{\frac{1}{p},1-\frac{1}{p}\Big\},\\
        &~~~~~~~~-\min\Big\{r_i,\frac{1}{q}-\frac{1}{2}\Big\}+\frac{2}{p}\Big\}
     \end{split}   
   \end{equation}
then the inhomogeneous Besov spaces with dominating mixed smoothness $S^{\bar{r}}_{p,q}B(\R)$ can be discretized in the sense of \eqref{exp} using the sequence space $s^{\bar{r}}_{p,q}b$.
\end{itemize}
\end{Theorem}

\begin{corollary}\label{corsplines} The wavelet basis characterization of the previous theorem holds for the choice of an 
orthogonal spline wavelets system $(\varphi_m,\psi_m)$
of order $m$, see Appendix \ref{splines}. For $S^{\bar{r}}_{p,q}B(\R)$ we need for $i=1,\hdots,d$
$$
    m-1 > \mbox{rhs}\eqref{eq-72}\,,
$$
whereas in case $S^{\bar{r}}_{p,q}F(\R)$ we need for $i=1,\hdots,d$
$$
   m-1>\mbox{rhs}\eqref{eq-71}\,.
$$
\end{corollary}
\bproof We apply Theorem \ref{wbasestensor} and the fact that 
$\psi^1 = \psi_m$ satisfies $(M_{L-1})$ for $L=m$ and $\psi^0=\varphi_m, \psi^1 = \psi_m$ satisfy $(D)$ and $(S_{K})$ for $K<m-1$.
\eproof

\begin{remark} 
\begin{itemize}
\item[(i)] Atomic decompositions of spaces with dominating mixed smoothness were already given by Vyb\'iral \cite{Vy06}. He provides compactly supported atomic decompositions and in particular wavelet isomorphisms in terms of compactly supported Daubechies wavelet. Bazarkhanov \cite{Ba03} provided the $\varphi$-transform for dominating mixed spaces and obtained atomic decompositions in the sense of Frazier, Jawerth. 
\item[(ii)] Wavelet bases in terms of orthonormal spline wavelets with optimal conditions on the order $m$ were given in 
\cite{SiUl09} in case $p=q$. However, this restriction is due to the tensor product approach in \cite{SiUl09}, and is not needed in our result.
\end{itemize}
\end{remark}

\setcounter{section}{0}
\renewcommand{\thesection}{\Alph{section}}
\renewcommand{\theequation}{\Alph{section}.\arabic{equation}}
\renewcommand{\theTheorem}{\Alph{section}.\arabic{Theorem}}
\renewcommand{\thesubsection}{\Alph{section}.\arabic{subsection}}

\section{Appendix: Wavelets}
For the notion of multi-resolution analysis, 
scaling function and associated wavelet we refer to Wojtaszczyk \cite[2.2]{Wo97} and Daubechies \cite[Chapt.\ 5]{Dau92}.

\subsection{Spline wavelets}
\label{splines}
As a main example we use the spline wavelet system $(\varphi_m,\psi_m)$. Let us recall here the basic 
construction and refer to \cite[Chap.\ 3.3]{Wo97} for the properties listed below.
The normalized cardinal B-spline of order
$m+1$ is given by
\[
\cn_{m+1} (x):= \cn_m * \cx (x)\, , \qquad x \in \re\, , \quad m \in
\N\, ,
\]
beginning with $\cn_1 = \cx$, the characteristic function of the interval $(0,1)$. 
By
\[
\varphi_m (x):= \frac{1}{\sqrt{2\pi}} \,  \left[\frac{\widehat{\cn}_m
(\xi)}{\Big(\sum\limits_{k=-\infty}^\infty |\widehat{\cn}_m (\xi + 2\pi
k)|^2\Big)^{1/2}}\right]^{\vee}(x)\, , \qquad x \in \re\, ,
\]
we obtain an orthonormal scaling function which is again a spline of
order $m$. 
Finally, by
\[
\psi_m (x) := \sum_{k=-\infty}^\infty \langle \, \varphi_m (t/2),
\varphi_m (t-k)\rangle\, (-1)^k \, \varphi_m (2x+k+1)
\]
the generator of an orthonormal wavelet system is defined. For $m=1$
it is easily checked that $-\psi_1 (\cdot-1)$ is the Haar wavelet. In
general, these functions $\psi_m$ have the following properties:
\begin{itemize}
\item Restricted to intervals $[\frac{k}{2},\frac{k+1}{2}]$, $k\in \zz$,
$\psi_m$ is a polynomial of degree at most $m-1$;
\item $\psi_m \in C^{m-2} (\re)$ if $m\ge 2$;
\item The derivative $\psi_m^{(m-2)}$ is uniformly Lipschitz continuous on
$\re$ if $m \ge 2$.
\item
The function $\psi_m$ satisfies moment conditions of order up to $m-1$,
i.e.,
\[
\int_{-\infty}^\infty x^\ell \, \psi_m (x)\, dx =0\, , \qquad \ell
=0,1,\ldots, m-1\, .
\]
In particular, $\psi_m$ satisfies $(M_{L-1})$ for $L=m$ and $\varphi_m, \psi_m$ satisfy $(D)$ and $(S_{K})$ for $K<m-1$.
\end{itemize}

\subsection{Tensor product wavelet bases on $\mathbf{\R}$}
\label{appdwavelet}
%We are interested
%in wavelet bases on $\R$. 
There is a straightforward method to construct a wavelet basis on $\R$ from a wavelet basis on $\re$. Putting 
$$
      \psi_{j,k} = \left\{\begin{array}{rcl}
                             \psi^0(\cdot-k)&:& j=0\\
                             2^{j/2}\psi^1(2^j\cdot-k)&:&j\geq 1
                          \end{array}\right.\quad,\quad j\in \n, k\in \zz
$$
and
$$
   \psi_{\bar{j},\bar{k}} = \psi_{j_1,k_1}\otimes\cdots \otimes \psi_{j_d,k_d}\,\quad,\quad \bar{j}=(j_1,\hdots,j_d)\in \n^d,\,
   \bar{k} = (k_1,\hdots,k_d)\in \Z\,,
$$
we obtain the following
\begin{lemma} Let  $\psi^{0} \in L_2(\re)$ be an orthonormal scaling function with associated orthonormal wavelet 
$\psi^{1} \in L_2(\re)$. Then the system
$$
   \{\psi_{\bar{j},\bar{k}}~:~ \bar{j}\in \n^d,\bar{k}\in \Z \Big\}
$$ 
is an orthonormal basis in $L_2(\R)$.
\end{lemma}

The next construction is slightly more involved. The following lemma is taken from \cite[1.2.1]{Tr08}. 
\begin{lemma}\label{dwavelet} Suppose, that we have a multi-resolution analysis
in $L_2(\re)$ with scaling function $\psi^{0}$ and associated
wavelet $\psi^{1}$. Let $E = \{0,1\}^d$, $c = (c_1,...,c_d) \in E$, and $\psi^c =
\bigotimes_{j=1}^d \psi^{c_j}$. Then the system
\begin{equation}\nonumber
        \{\psi^0(x-k)\,:\,k\in \Z\} \cup \Big\{2^{\frac{jd}{2}}\psi^{c}(2^jx-k)\,:\,c\in E\setminus \{0\}, j\in \n,k\in \Z \Big\}
\end{equation} is an orthonormal basis in $L_2(\R)$.
\end{lemma}


\begin{thebibliography}{10}

\bibitem{alanga93}
S.~T. {A}li, J.-P. {A}ntoine, and J.-P. {G}azeau.
\newblock {C}ontinuous frames in {H}ilbert space.
\newblock {\em {A}nn. {P}hysics}, 222(1):1--37, 1993.

\bibitem{Ba03}
D.~B. {B}azarkhanov.
\newblock {C}haracterizations of the {N}ikol'skij-{B}esov and
  {L}izorkin-{T}riebel function spaces of mixed smoothness.
\newblock {\em Proc. Steklov Inst. Math.}, 243:53--65, 2003.

\bibitem{Bo95}
G.~{B}ourdaud.
\newblock Ondelettes et espaces de {B}esov.
\newblock {\em Rev. Matem. Iberoam.}, 11:477--512, 1995.

\bibitem{Bo07}
M.~Bownik.
\newblock {A}nisotropic {T}riebel-{L}izorkin spaces with doubling measures.
\newblock {\em Journ. of Geom. Anal.}, 17(3):387--424, 2007.

\bibitem{Bu82}
H.-Q. Bui.
\newblock Weighted {B}esov and {T}riebel spaces: interpolation by the real
  method.
\newblock {\em Hiroshima Math. J.}, 12(3):581--605, 1982.

\bibitem{Bu84}
H.-Q. Bui.
\newblock Characterizations of weighted {B}esov and {T}riebel-{L}izorkin spaces
  via temperatures.
\newblock {\em J. Funct. Anal.}, 55(1):39--62, 1984.

\bibitem{BuPaTa96}
H.-Q. {B}ui, M.~{P}aluszy\'nski, and M.~H. {T}aibleson.
\newblock A maximal function characterization of weighted {B}esov-{L}ipschitz
  and {T}riebel-{L}izorkin spaces.
\newblock {\em Stud. Math.}, 119(3):219--246, 1996.

\bibitem{BuPaTa97}
H.-Q. {B}ui, M.~{P}aluszy\'nski, and M.~H. {T}aibleson.
\newblock Characterization of the {B}esov-{L}ipschitz and {T}riebel-{L}izorkin
  spaces, the case $q<1$.
\newblock {\em J. Four. Anal. and Appl. (special issue)}, 3:837--846, 1997.

\bibitem{ChFr87}
F.~Chiarenza and M.~Frasca.
\newblock {M}orrey spaces and {H}ardy-{L}ittlewood maximal functions.
\newblock {\em Rend. Math.}, 7:273--279, 1987.

\bibitem{CoDaFe92}
A.~Cohen, I.~Daubechies, and J.-C. Feauveau.
\newblock Biorthogonal bases of compactly supported wavelets.
\newblock {\em Comm. Pure Appl. Math.}, 45(5):485--560, 1992.

\bibitem{CoMeSt85}
R.~R. {C}oifman, Y.~{M}eyer, and E.~M. {S}tein.
\newblock {S}ome new function spaces and their application to harmonic
  analysis.
\newblock {\em Journ. of Funct. Anal.}, 62:304--335, 1985.

\bibitem{daforastte08}
S.~{D}ahlke, M.~{F}ornasier, H.~{R}auhut, G.~{S}teidl, and G.~{T}eschke.
\newblock {G}eneralized coorbit theory, {B}anach frames, and the relation to
  alpha-modulation spaces.
\newblock {\em {P}roc. {L}ondon {M}ath. {S}oc. (3)}, 96:464--506, 2008.

\bibitem{dakustte09}
S.~{D}ahlke, G.~{K}utyniok, G.~{S}teidl, and G.~{T}eschke.
\newblock {S}hearlet coorbit spaces and associated {B}anach frames.
\newblock {\em {A}ppl. {C}omput. {H}armon. {A}nal.}, 27(2):195--214, 2009.

\bibitem{dastte04}
S.~{D}ahlke, G.~{S}teidl, and G.~{T}eschke.
\newblock {C}oorbit spaces and {B}anach frames on homogeneous spaces with
  applications to the sphere.
\newblock {\em {A}dv. {C}omput. {M}ath.}, 21(1-2):147--180, 2004.

\bibitem{dastte04-1}
S.~{D}ahlke, G.~{S}teidl, and G.~{T}eschke.
\newblock {W}eighted coorbit spaces and {B}anach frames on homogeneous spaces.
\newblock {\em {J}. {F}ourier {A}nal. {A}ppl.}, 10(5):507--539, 2004.

\bibitem{Dau92}
I.~Daubechies.
\newblock {\em {T}en {L}ectures on {W}avelets}, volume~61 of {\em CBMS-NSF
  Regional Conference Series in Applied Mathematics}.
\newblock Society for Industrial and Applied Mathematics (SIAM), Philadelphia,
  PA, 1992.

\bibitem{FeMu74}
C.~Fefferman and B.~Muckenhoupt.
\newblock Two nonequivalent conditions for weight functions.
\newblock {\em Proc. Amer. Math. Soc.}, 45:99--104, 1974.

\bibitem{fe83-4}
H.~G. {F}eichtinger.
\newblock {M}odulation spaces on locally compact {A}belian groups.
\newblock Technical report, {J}anuary 1983.

\bibitem{fegr85}
H.~G. {F}eichtinger and P.~{G}r{\"o}bner.
\newblock {B}anach spaces of distributions defined by decomposition methods.
  {I}.
\newblock {\em {M}ath. {N}achr.}, 123:97--120, 1985.

\bibitem{FeGr86}
H.~G. Feichtinger and K.~Gr{\"o}chenig.
\newblock A unified approach to atomic decompositions via integrable group
  representations.
\newblock In {\em Function spaces and applications (Lund, 1986)}, volume 1302
  of {\em Lecture Notes in Math.}, pages 52--73. Springer, Berlin, 1988.

\bibitem{FeGr89a}
H.~G. {F}eichtinger and K.~{G}r{\"o}chenig.
\newblock {B}anach spaces related to integrable group representations and their
  atomic decompositions, {I}.
\newblock {\em Journ. Funct. Anal.}, 21:307--340, 1989.

\bibitem{FeGr89b}
H.~G. {F}eichtinger and K.~{G}r{\"o}chenig.
\newblock {B}anach spaces related to integrable group representations and their
  atomic decompositions, {II}.
\newblock {\em {M}onatsh. {M}athem.}, 108:129--148, 1989.

\bibitem{fegr92-1}
H.~G. {F}eichtinger and K.~{G}r{\"o}chenig.
\newblock {G}abor wavelets and the {H}eisenberg group: {G}abor expansions and
  short time {F}ourier transform from the group theoretical point of view.
\newblock In C.~K. {C}hui, editor, {\em {W}avelets :a tutorial in theory and
  applications}, volume~2 of {\em {W}avelet {A}nal. {A}ppl.}, pages 359--397.
  {A}cademic {P}ress, 1992.

\bibitem{fora05}
M.~{F}ornasier and H.~Rauhut.
\newblock {C}ontinuous frames, function spaces, and the discretization problem.
\newblock {\em {J}. {F}ourier {A}nal. {A}ppl.}, 11(3):245--287, 2005.

\bibitem{FrJa90}
M.~{F}razier and B.~{J}awerth.
\newblock {A} discrete transform and decompositions of distribution spaces.
\newblock {\em Journ. of Funct. Anal.}, 93:34--170, 1990.

\bibitem{gr92-2}
P.~{G}r{\"o}bner.
\newblock {\em {B}anachr{\"a}ume glatter {F}unktionen und
  {Z}erlegungsmethoden}.
\newblock PhD thesis, 1992.

\bibitem{Gr88}
K.~Gr{\"o}chenig.
\newblock Unconditional bases in translation and dilation invariant function
  spaces on {$\mathbb{R}\sp n$}.
\newblock In {\em Constructive theory of functions (Varna, 1987)}, pages
  174--183. Publ. House Bulgar. Acad. Sci., Sofia, 1988.

\bibitem{Gr91}
K.~{G}r{\"o}chenig.
\newblock {D}escribing functions: atomic decompositions versus frames.
\newblock {\em {M}onatsh. {M}athem.}, 112:1--41, 1991.

\bibitem{gr01}
K.~{G}r{\"o}chenig.
\newblock {\em {F}oundations of {T}ime-{F}requency {A}nalysis}.
\newblock {A}ppl. {N}umer. {H}armon. {A}nal. {B}irkh{\"a}user {B}oston, 2001.

\bibitem{HaDiss10}
M.~{H}ansen.
\newblock {\em {N}onlinear {A}pproximation and {F}unction {S}paces of
  {D}ominating {M}ixed {S}moothness}.
\newblock PhD thesis, FSU Jena, Germany, 2010.

\bibitem{HaPi08}
D.~D. Haroske and I.~Piotrowska.
\newblock Atomic decompositions of function spaces with {M}uckenhoupt weights,
  and some relation to fractal analysis.
\newblock {\em Math. Nachr.}, 281(10):1476--1494, 2008.

\bibitem{HeNe07}
L.~I. {H}edberg and Y.~{N}etrusov.
\newblock {A}n axiomatic approach to function spaces, spectral synthesis, and
  {L}uzin approximation.
\newblock {\em Mem. Amer. Math. Soc.}, 188(882), 2007.

\bibitem{Ke09}
H.~{K}empka.
\newblock 2-microlocal spaces of variable integrability.
\newblock {\em {R}ev. {M}at. {C}omplut.}, 22(1):227--251, 2009.

\bibitem{Ke10}
H.~{K}empka.
\newblock Atomic, molecular and wavelet decomposition of generalized
  2-microlocal {B}esov spaces.
\newblock {\em Journ. Function Spaces Appl.}, 8:129--165, 2010.

\bibitem{Ke11}
H.~{K}empka.
\newblock Atomic, molecular and wavelet decomposition of 2-microlocal {B}esov
  and {T}riebel-{L}izorkin spaces with variable integrability.
\newblock {\em Funct. Appr. Comment. Math.}, 43(2):171--208, 2010.

\bibitem{Ky03}
G.~Kyriazis.
\newblock Decomposition systems for function spaces.
\newblock {\em Studia Math.}, 157:133--169, 2003.

\bibitem{Mo38}
C.~B. Morrey.
\newblock On the solutions of quasi linear elliptic partial differential
  equations.
\newblock {\em Trans. Amer. Math. Soc.}, 43:126--166, 1938.

\bibitem{Pe75}
J.~{P}eetre.
\newblock {O}n spaces of {T}riebel-{L}izorkin type.
\newblock {\em {A}rk. {M}at.}, 13:123--130, 1975.

\bibitem{ra05}
H.~{R}auhut.
\newblock {B}anach frames in coorbit spaces consisting of elements which are
  invariant under symmetry groups.
\newblock {\em {A}ppl. {C}omput. {H}armon. {A}nal.}, 18(1):94--122, 2005.

\bibitem{ra05-6}
H.~{R}auhut.
\newblock {\em {T}ime-frequency and wavelet analysis of functions with symmetry
  properties}.
\newblock {L}ogos-{V}erlag, 2005.
\newblock {P}h{D} thesis.

\bibitem{ra05-2}
H.~{R}auhut.
\newblock {R}adial time-frequency analysis and embeddings of radial modulation
  spaces.
\newblock {\em {S}ampl. {T}heory {S}ignal {I}mage {P}rocess.}, 5(2):201--224,
  2006.

\bibitem{ra05-3}
H.~{R}auhut.
\newblock {C}oorbit space theory for quasi-{B}anach spaces.
\newblock {\em {S}tudia {M}ath.}, 180(3):237--253, 2007.

\bibitem{Ry99a}
V.~S. {R}ychkov.
\newblock {O}n a theorem of {B}ui, {P}aluszy\'nski and {T}aibleson.
\newblock {\em {P}roc. {S}teklov {I}nst.}, 227:280--292, 1999.

\bibitem{Sa08}
Y.~Sawano.
\newblock Wavelet characterization of {B}esov-{M}orrey and
  {T}riebel-{L}izorkin-{M}orrey spaces.
\newblock {\em Funct. Approx. Comment. Math.}, 38:93--107, 2008.

\bibitem{sa08-1}
Y.~{S}awano.
\newblock {W}eighted modulation space $m_{p,q}^s(w)$ with $w \in
  a_\infty^{loc}$.
\newblock {\em {J}ournal of {M}athematical {A}nalysis and {A}pplications},
  345(2):615--627, 2008.

\bibitem{SaTa07}
Y.~Sawano and H.~Tanaka.
\newblock Decompositions of {B}esov-{M}orrey spaces and
  {T}riebel-{L}izorkin-{M}orrey spaces.
\newblock {\em Math. Zeitschr.}, 257:871--905, 2007.

\bibitem{SaYaYu10}
Y.~Sawano, D.~Yang, and W.~Yuan.
\newblock New applications of {B}esov-type and {T}riebel-{L}izorkin-type
  spaces.
\newblock {\em J. Math. Anal. Appl.}, 363:73--85, 2010.

\bibitem{ScTr87}
H.-J. {S}chmeisser and H.~{T}riebel.
\newblock {\em {T}opics in {F}ourier {A}nalysis and {F}unction {S}paces}.
\newblock Wiley, Chichester, 1987.

\bibitem{SiUl09}
W.~{S}ickel and T.~{U}llrich.
\newblock {T}ensor products of {S}obolev-{B}esov {S}paces and applications to
  approximation from the hyperbolic cross.
\newblock {\em Journ. Approx. Theory}, 161:748--786, 2009.

\bibitem{SiUl10}
W.~{S}ickel and T.~{U}llrich.
\newblock {S}pline interpolation on sparse grids.
\newblock {\em Applicable Analysis}, to appear.

\bibitem{StWe71}
E.~M. Stein and G.~{W}eiss.
\newblock {\em {I}ntroduction to {F}ourier {A}nalysis on {E}uclidean {S}paces}.
\newblock Princeton Univ. Press, 1971.

\bibitem{TaXu05}
L.~Tang and J.~Xu.
\newblock Some properties of {M}orrey type {B}esov-{T}riebel spaces.
\newblock {\em Math. Nachr.}, 278:904--917, 2005.

\bibitem{Tr83}
H.~{T}riebel.
\newblock {\em {T}heory of {F}unction {S}paces}.
\newblock {B}irkh{\"a}user, Basel, 1983.

\bibitem{Tr88}
H.~{T}riebel.
\newblock {C}haracterizations of {B}esov-{H}ardy-{S}obolev spaces: a unified
  approach.
\newblock {\em {J}ourn. of {A}pprox. {T}heory}, 52:162--203, 1988.

\bibitem{Tr92}
H.~{T}riebel.
\newblock {\em {T}heory of {F}unction {S}paces II}.
\newblock {B}irkh{\"a}user, Basel, 1992.

\bibitem{Tr06}
H.~{T}riebel.
\newblock {\em Theory of {F}unction {S}paces III}.
\newblock Birkh{\"a}user, Basel, 2006.

\bibitem{Tr08}
H.~{T}riebel.
\newblock {\em Function {S}paces and {W}avelets on {D}omains}.
\newblock EMS Publishing House, Z{\"u}rich, 2008.

\bibitem{Tr09}
H.~Triebel.
\newblock {\em {B}ases in {F}unction {S}paces, {S}ampling, {D}iscrepancy,
  {N}umerical {I}ntegration}.
\newblock EMS Publ. House, Z{\"u}rich, 2010.

\bibitem{TDiff06}
T.~{U}llrich.
\newblock {F}unction spaces with dominating mixed smoothness, characterization
  by differences.
\newblock {\em {J}enaer {S}chriften zur {M}athematik und {I}nformatik},
  Math/Inf/05/06:1--50, 2006.

\bibitem{T08_2}
T.~{U}llrich.
\newblock Local mean characterization of {B}esov-{T}riebel-{L}izorkin type
  spaces with dominating mixed smoothness on rectangular domains.
\newblock {\em Preprint}, pages 1--26, 2008.

\bibitem{T10}
T.~{U}llrich.
\newblock Continuous characterizations of {B}esov-{L}izorkin-{T}riebel spaces
  and new interpretations as coorbits.
\newblock {\em Journ. Funct. Spaces Appl.}, to appear.

\bibitem{Vy06}
J.~{V}yb\'iral.
\newblock {F}unction spaces with dominating mixed smoothness.
\newblock {\em {D}iss. {M}ath.}, 436:73~pp, 2006.

\bibitem{Wo97}
P.~{W}ojtaszczyk.
\newblock {\em {A} {M}athematical {I}ntroduction to {W}avelets}.
\newblock {C}ambridge {U}niversity {P}ress, 1997.

\bibitem{YuSiYa10}
W.~{Y}uan, W.~{S}ickel, and D.~{Y}ang.
\newblock {M}orrey and {C}ampanato meet {B}esov, {L}izorkin and {T}riebel.
\newblock volume 2005 of {\em Lecture Notes in Math.}, pages xi+281. Springer,
  Berlin, 2010.

\end{thebibliography}
\end{document}